\documentclass[12pt,twoside]{report}
\usepackage{appendix}
\usepackage{graphicx}
\usepackage{amsmath}
\usepackage{setspace}
\usepackage{color}
\usepackage[byname]{smartref}
\usepackage[unicode=true,
bookmarks=false,
breaklinks=false,pdfborder={0 0 1},backref=section,colorlinks=false]
{hyperref}
\hypersetup{
	colorlinks,bookmarksopen,bookmarksnumbered,citecolor=blue,urlcolor=blue}
\usepackage{cite}
\usepackage{txfonts}
\usepackage{tocloft}
\newcommand{\Lone}{$ \mathcal{L}_1 \,$}

\newtheorem{chapter_c}{Chapter}

\newtheorem{lemma}{Lemma}

\newtheorem{theorem}{Theorem}
\newtheorem{rem}{Remark}

\newtheorem{assumption}{Assumption}
\newtheorem{assum}{Assumption}

\makeatletter
\numberwithin{figure}{chapter}
{\clearemptydoublepage
 \begin{center}
  \section*{Acknowledgements}
 \end{center}
 \begingroup
}{\newpage\endgroup}

{\clearemptydoublepage 
 \begin{center}
  \section*{Dedication}
 \end{center}
 \begingroup
}{\newpage\endgroup}

\newenvironment{preliminary}%
{\pagestyle{plain}\pagenumbering{roman}}%
{\pagenumbering{arabic}}

\addtoreflist{chapter}

\newenvironment{remark}[1][Remark]{\begin{trivlist}
\item[\hskip \labelsep {\bfseries #1}]}{\end{trivlist}}

\newcommand{\qed}{\nobreak \ifvmode \relax \else
      \ifdim\lastskip<1.5em \hskip-\lastskip
      \hskip1.5em plus0em minus0.5em \fi \nobreak
      \vrule height0.75em width0.5em depth0.25em\fi}

\if@twoside 
\def\ps@thesis{\let\@mkboth\markboth
   \def\@oddfoot{}
   \let\@evenfoot\@oddfoot
   \def\@oddhead{
      {\sc\rightmark} \hfil \rm\thepage
      }
   \def\@evenhead{
      \rm\thepage \hfil {\sc\leftmark}
      }
   \def\chaptermark##1{\markboth{\ifnum \c@secnumdepth >\m@ne
      Chapter\ \thechapter. \ \fi ##1}{}}
   \def\sectionmark##1{\markright{\ifnum \c@secnumdepth >\z@
      \thesection. \ \fi ##1}}}
\else 
\def\ps@thesis{\let\@mkboth\markboth
   \def\@oddfoot{}
   \def\@oddhead{
      {\sc\rightmark} \hfil \rm\thepage
      }
   \def\chaptermark##1{\markright{\ifnum \c@secnumdepth >\m@ne
      Chapter\ \thechapter. \ \fi ##1}}}
\fi

\newcommand\isco[1]{%
  \edef\@tempa{#1}%
  \def\@tempb{}%
  \ifx\@tempa\@tempb
	\else \\\underline{Co-Supervisor:}\vspace{0.35in}\\\dots\dots\dots\dots\dots\dots\dots\\{#1}\\
  \fi
}

\newcommand\isjoint[1]{%
  \edef\@tempa{#1}%
  \def\@tempb{}%
  \ifx\@tempa\@tempb
	\else \\\underline{Joint Supervisor:}\vspace{0.35in}\\\dots\dots\dots\dots\dots\dots\dots\\{#1}\\
  \fi
}

\newcommand\isalt[1]{%
  \edef\@tempa{#1}%
  \def\@tempb{}%
  \ifx\@tempa\@tempb
	\else \\\underline{Alternate Supervisor:}\vspace{0.35in}\\\dots\dots\dots\dots\dots\dots\dots\\{#1}\\
  \fi
}

\newcommand\isdefinedsig[1]{%
  \edef\@tempa{#1}%
  \def\@tempb{}%
  \ifx\@tempa\@tempb
	\else \\ \dots\dots\dots\dots\dots\dots\dots\\{#1}\\
  \fi
}
\newcommand\isdefinedspinetitle[1]{%
  \edef\@tempa{#1}%
  \def\@tempb{}%
  \ifx\@tempa\@tempb
	\else (Spine title: #1)\\
  \fi
}
\newcommand\coauthor[1]{%
  \edef\@tempa{#1}%
  \def\@tempb{}%
  \ifx\@tempa\@tempb
	\else \newpage \Large Co-Authorship Statement\normalsize\\\indent\\#1\\
  \fi
}

\newcommand\acknowlege[1]{%
  \edef\@tempa{#1}%
  \def\@tempb{}%
  \ifx\@tempa\@tempb
	\else \newpage \Large Acknowlegements\normalsize\\\indent\\#1\newpage
  \fi
}

\newcommand{\department}{SYSTEMS Engineering}
\newcommand{\degree}{Masters of Science}
\newcommand{\firstname}{HASHIM}
\newcommand{\middlename}{ABDELLAH HASHIM}
\newcommand{\lastname}{MOHAMED}
\newcommand{\authorname}{{\firstname} {\middlename} {\lastname}}
\newcommand{\titl}{Improved Robust Adaptive Control of High-order Nonlinear Systems with Guaranteed Performance}

\newcommand{\listappendixname}{List of Appendices}
\newlistof{myappendices}{app}{\listappendixname}

\renewcommand{\maketitle}
{\begin{titlepage}
   \setcounter{page}{1}
   \begin{large}
   \begin{center}
      \mbox{}
      \vfill
      {\MakeUppercase{\titl}}\\
      \vfill
      by \\
      \vfill
      {\authorname} \\
      \vfill
      Graduate Program in {\department}\\
      \vfill
		A thesis submitted in partial fulfillment\\
		of the requirements for the degree of\\
		\degree\\
		\vfill
		SYSTEMS ENGINEERING\\
		KING FAHD UNIVERSITY OF PETROLEUM \& MINERALS\\
		\vfill
      {\copyright} {\authorname}   \\
      {December 2014}  \\
      \vspace*{.2in}
   \end{center}
   \end{large}
   \end{titlepage}

}

\newcommand{\makecert}{
   \setcounter{page}{2}
\vfill
}

\makeatother
\begin{document}


\begin{preliminary}

{ \bf The citation of this M.Sc thesis follows:}
\\

{ \bf \textcolor{red}{Hashim Abdellah Hashim Mohamed, Improved robust adaptive control of 
	 high-order nonlinear systems with guaranteed performance. M.Sc, King 
	 Fahd University Of Petroleum \& Minerals, 2014.}
}

\begin{center}
\vspace{500pt}\footnotesize { \bf
	Please contact us and provide details if you believe this document breaches copyrights. We will remove access to the work immediately and investigate your claim. }
\end{center}

\maketitle
\addcontentsline{toc}{chapter}{Certificate of Examination}
\makecert
\newpage
\addcontentsline{toc}{chapter}{Abstract}
\Large\begin{center}\textbf{Abstract}\end{center}\normalsize

This thesis presents fuzzy-\Lone adaptive controller and Model Reference Adaptive Control (MRAC) with Prescribed Performance Function (PPF) as two adaptive approaches for high nonlinear systems as two original contribution to the literature. Firstly, \Lone adaptive controller has a structure that allows decoupling between robustness and adaption owing to the use of a low pass filter with adjustable gain in the feedback loop. The trade-off between performance and robustness is a key factor in the tuning of the filter's parameters. In fuzzy-\Lone adaptive controller, we consider the class of strictly proper low pass filters with fixed structure but with the feedback gain as the only tunable parameter. A practical new fuzzy based approach for the tuning of the feedback filter of \Lone adaptive controller is proposed. The fuzzy controller is optimally tuned using Particle Swarm Optimization (PSO) to minimize the tracking error and the control signal range. The main function of the fuzzy logic controller is the on-line tuning of the feedback gain of the filter. Secondly, an adaptive control of multi-input multi-output uncertain high-order nonlinear system capable of guaranteeing a predetermined prescribed performance is presented as MRAC with PPF. In this work, prescribed performance is defined in terms of the tracking error converging to a smaller residual set at a rate no less than a predefined value and exhibiting a maximum overshoot/undershoot less than a sufficiently small fixed constant. The key step in such approach is to transform the constrained system into an equivalent unconstrained one through an adequate transformation of the output error. This will show that the robust stabilization of the transformed error, guaranties the stability and convergence of the constrained tracking error within the set of time varying constraints representing the performance limits. Finally, simulations are presented to illustrate the simplicity, the performance and the robustness of each new technique.

\vfill
\newpage
\addcontentsline{toc}{chapter}{Abstract}
\Large\begin{center}\textbf{Acknowledgments}\end{center}\normalsize 

This thesis is the result of my Master degree from January 2013 to December 2014. It has taken place at the Department of System Engineering at King Fahd University of Petroleum and Minerals under the counseling of my supervisor Dr. Sami El Ferik. First I would like to thank  Dr. Sami El Ferik. His exceptional motivational skills and ability to continuously guide me in the right directions have been much needed assets in my struggle towards finishing the Master Degree.I would like to thank also Dr. Mustafa El-Shafei, and Dr. Mohamed Abido for being my thesis committee members.

\clearpage  
\tableofcontents\newpage
\newpage
\addcontentsline{toc}{chapter}{List of Figures}
\listoffigures
\newpage
\addcontentsline{toc}{chapter}{List of Tables}
\listoftables\newpage
\end{preliminary}



\newpage


\chapter{INTRODUCTION}

\section{Introduction And Motivation}
The presence of uncertainties, nonlinearities, disturbances and lack in the precise modeling of nonlinear systems are common problems in dynamical applications. Over the last few decades, adaptive control has been developed to tackle the foregoing problems by providing fast adaption and ensure robustness. In this work, \Lone adaptive controller will be discussed briefly from different perspectives for different systems structures. \Lone adaptive controller has been inspired originally from  MRAC. Improving the feedback filter of \Lone adaptive control will enhance the performance of the controller and the robustness margin. Fuzzy filter will be proposed for \Lone adaptive controller in order to ensure fast closed loop dynamics with increasing the robustness margin. Neuro adaptive control with prescribed performance function will be investigated. Robust Model Reference Adaptive Control (MRAC) with Prescribed Performance Function (PPF) will be proposed to tackle problems of neuro-adaptive control and comparing the controller performance versus \Lone adaptive controller. Robust adaptive observer will be implemented with \Lone adaptive controller in order to check the performance of the controller in case of inaccessible states. These controllers will be applied on high nonlinear systems including Unmanned Vehicle Systems (UVS).

\section{Possible Applications of The Outcomes}
Unmanned Vehicle Systems (UVS) are important for different areas nowadays because they can be controlled and operated remotely without human interference. UVS is a research key because of the increase in demand of remote sensing and control in wide range of applications such as scientific surveys, traffic surveillance, transportation aids, and inspection in addition to operation in harsh environments. UVS have various configurations, characteristics, shapes and sizes which will be reflected on system dynamics. The development in miniaturization of UVS offers high potential effort for small size and low cost of UVS compared to manned applications especially in certain applications. Rapid growing of UVS comes with promising future because of its size, cost, construction simplicity and maneuverability.
 
UVS can be classified into two categories either remotely control vehicles, or autonomous vehicles. Each of these categories includes different types of UVS such as: Unmanned Aerial Vehicles, Underwater Vehicles, Unmanned Surface Vehicle, Unmanned Spacecraft and Unmanned Grounded Vehicle. Importance of UVS relies on performance and mission targets. Generally, each type is considered as a mechanical rigid body with different equations of motion. The majority of UVS can be represented by nonlinear dynamics. The dynamic of UVS have their own features as affine nonlinear systems with normal coupling or with strong coupling. Usually, the controller is required to drive the system to the desired trajectory with smooth transition and fast response. Smooth transition in both control signal and output response will contribute in protecting the life cycle of system rotors and other parts in the UVS.

Developing UVS in the absence of the operator is costly in the controller complexity for tracking and vision. The controller is demanded to overcome many drawbacks, starting with stabilizing the system, driving the system to the desired trajectory in the shortest possible time, adapt against any variations of system dynamics and finally be robust against any disturbances. All these requirements ended up making the control design as an important issue and an interested subject to be investigated.

In order to design a controller for UVS, accurate models are needed to reflect system dynamics either by precise modeling or real time identification. UVS have a framework of rigid body dynamics and can be described by a set of differential equations using Euler-Lagrange. The definition of exact model is a struggling problem because nominal model is usually defined under certain operating conditions with neglecting any uncertainties and disturbances that may exist during the control process. Classical controller will not be sufficient due to nonexact model represented by presence of uncertainties and/or disturbances. Other types of controllers have to be considered in order to overcome classical controller drawbacks.

In the literature, several control design approaches have been adopted for Euler-Lagrange systems like adaptive control, nonlinear control, robust control and so forth. The weakness of many control approaches resides in defining the appropriate model for nonlinearity cancellation. In nonlinear control, it is often difficult to use the approximated nonlinear Euler-Lagrange equations of the system without adding a robustifing term to ensure system operation in the stability region. Including a robustifing term in the control law introduces discontinuity and chattering on the control signal. On the other hand, estimation of system nonlinearities normally experienced with discontinuity or singularity in the estimation process which may take the system out of the stability region.

\section{Contribution to The Literature}
In our work, two robust adaptive control approaches will be proposed for high nonlinear systems with guaranteed performance. Firstly, A fuzzy logic feedback filter will be designed for \Lone adaptive controller mainly to improve the tracking capability and reduce the control signal range. The trade off between robustness range and fast closed loop dynamics will be averted and the proposed controller will contribute in solving this major problem. Next, robust MRAC-PPF will be proposed to tackle limitations of robust neuro-adaptive control with PPF. Also, it will be compared versus \Lone adaptive control to highlight merits of the new controller. The controller will be studied on affine and not-affine systems. Finally, the performance of \Lone adaptive controller with adaptive observers will be examined on Single-Input Single-Output (SISO) and Multi-Input Multi-Output (MIMO) systems.

The main features of the \Lone adaptive controller are:
\begin{itemize}
 \item Estimating the system to be controlled.
 \item For linear and nonlinear case without strong coupling, procedures consist of estimating uncertainties of the states, unmodelled input parameters and disturbances. For nonlinear case with strong coupling and/or unmatched uncertainties, it has same previous estimation process in addition to the estimate of unmatched part.
 \item The control law is based on Lyapunov function with compact set for previous item will be computed numerically.
\end{itemize}

The main features of robust neuro adaptive control with PPF are:
\begin{itemize}
\item Assign the prescribed function.
\item Derive the transformed error.
\item Estimating nonlinearities by neural network.
\item Computing the control signal based on Lyapunov function.
\end{itemize}

\subsection{Thesis Objectives and Contribution}
This thesis contributes to literature on several routs all aiming at improving \Lone adaptive controller in terms of adaptation and robustness. Therefore, there are several problems to be considered in this thesis:
\begin{enumerate}
\item We design a stabilizing controller based on fuzzy-\Lone adaptive controller and examine the controller performance for nonlinear systems.
\item We design a stabilizing controller based on MRAC with PPF and examine the controller performance for nonlinear systems.
\item We compare fuzzy-\Lone adaptive controller to \Lone adaptive controller.
\item We compare MRAC to PPF versus neuro adaptive conrol with PPF and \Lone adaptive controller.
\item Furthermore, we develop and implement adaptive observer with \Lone adaptive control for nonlinear systems.
\end{enumerate}

\section{Methodologies}

Developing thesis objective as mentioned in the previous section will go through several steps as following
\begin{enumerate}
\item Different UVS and nonlinear models have to be addressed as equation of motions.
\item Reproduce recent results upon literature of \Lone adaptive control for nonlinear systems including UVS.
\item Reproduce recent results upon the literature on robust neuro adaptive control with prescribed performance function for nonlinear systems.
\item Formulate fuzzy-\Lone adaptive controller and validate the new controller assuming complete unknown of nonlinear dynamics.
\item Formulate MRAC with PPF and validate the new controller assuming complete unknown of nonlinear dynamics.
\item Evaluating the performance of the controller by benchmarking the results to results in the literatures.
\item Develop and implement adaptive observer with \Lone adaptive controller and benchmarking the results to results of \Lone adaptive controller with accessible states.
\end{enumerate}
Out of this work I have succeeded to publish \cite{hashim1,hashim2,hashim3}

\section{Thesis Organization}
The thesis is organized as the following
\begin{chapter_c}
 includes introduction of the main work, motivation, thesis objective, methodology and finally thesis organization.
\end{chapter_c}

\begin{chapter_c}
 includes literature review of different control methods especially adaptive control for nonlinear systems. Literature review presents last research activities on \Lone adaptive control. Literature review of adaptive control with prescribed performance presents the main research activities over the last few years. Literature review of observer design shows the main research activities on this field.
\end{chapter_c}

\begin{chapter_c}
 includes \Lone adaptive controller for uncertain SISO systems, for uncertain MIMO systems and for MIMO systems in the presence of unmatched nonlinear uncertainties with strong coupling. Stability analysis, problem formulation and simulations will be validated for all foregoing cases.
\end{chapter_c}

\begin{chapter_c}
includes a brief review of \Lone adaptive controller. It proposes a design of fuzzy logic control to tune the feedback filter of \Lone adaptive controller. PSO is presented to design the output membership function of FLC. The controller will be examined on highly nonlinear system.
\end{chapter_c}

\begin{chapter_c}
 includes robust neuro adaptive controller for strict feedback MIMO system with PPF mainly functioned to capture the idea of PPF in addition to evaluate its performance by reproducing recent papers.
\end{chapter_c}

\begin{chapter_c}
proposes a design of MRAC with PPF for high uncertain nonlinear systems. \Lone adaptive controller and neuro-adaptive control with PPF are compared to the proposed controller.
\end{chapter_c}

\begin{chapter_c}
presents robust adaptive observer with \Lone adaptive controller for highly nonlinear systems with complete unknown dynamics.
\end{chapter_c}

\begin{chapter_c}
concludes the work and suggests possible future works.
\end{chapter_c}



\newpage


\chapter{LITERATURE REVIEW}

\section{Introduction}

This chapter summarizes the research activities of \Lone adaptive controller and adaptive control with PPF on different nonlinear systems with complete unknown dynamics. The first section include an introduction. The second section presents literature review of various control methods of UVS and a literature review of adaptive control techniques. The main contribution of this work is developed. Section three presents a brief review on \Lone adaptive control including the main recent research activities. The fourth section is a review on adaptive control with PPF including  including main research activities and recent works. Section five presents a study review on observer design. The last section is a conclusion.\\

\section{Feedback Control of UVS} 

   Adaptive control emerged in order to tackle time variant uncertainties, unmodeled dynamics and disturbances. Over the last few decades, various types of adaptive control has been proposed and modified to manipulate with aforementioned problems such as self-tuning regulators \cite{aastrom_self_1973,fortescue_implementation_1981,aastrom_automatic_1984,yesil_self_2004}, gain scheduling \cite{rugh_analytical_1991,apkarian_advanced_1998,rugh_research_2000}, model reference adaptive control system \cite{parks_liapunov_1966,chen_model_1995,yu_model_1996,patino_neural_2000} and adaptive neuro fuzzy control system \cite{kim_hyfis:_1999,orlowska-kowalska_adaptive_2010,kayacan_adaptive_2013}. In the recent few years, new adaptive control techniques were proposed rely on previous methods in terms of stability criteria and control law formulation. Immersion and Invariance adaptive control which is based on system immersion and manifold invariance was developed in order to reduce the control law and to ensure the asymptotic stability of the system \cite{astolfi_immersion_2003,liu_immersion_2010,hu_immersion_2013,bustan_immersion_2014}. Robust adaptive control with prescribed performance function mainly developed to force the error to start within large set and end within pre-assigned small set \cite{bechlioulis_robust_2008,bechlioulis_adaptive_2009,na_adaptive_2013}. \Lone adaptive control was developed to guarantee boundedness of transient and steady state performance in the absence knowledge of system nonlinearities, uncertainties and any disturbance \cite{cao_design_2008,luo_l_2010,xargay_l1_2010}.
   
   UVS control had been studied by many researchers trying to find a solution for improving the transient response and tracking trajectory. Sliding mode control for twin rotor MIMO system has been proposed in \cite{tao_novel_2010,mondal_adaptive_2012} where fuzzy control in \cite{tao_novel_2010} and adaptive rule technique in \cite{mondal_adaptive_2012} were used to cancel nonlinearities. Both techniques applied integral sliding mode for the vertical part with robust behavior against parameters variations and they showed great results. However, it has some intrinsic limitations due to design complexity, chattering on the sliding surface and manipulation of the controller only with strict feedback systems. Feedback linearization with sliding mode control for quadrotor has been implemented in \cite{lee_feedback_2009} and for micro unmanned automated vehicle was studied in \cite{voos_nonlinear_2009}. Limitations of feedback linearization is that the model should be in the strict feedback form and full knowledge of nonlinear model should be valid. In addition, uncertainties in model parameters should be within specific range. Backstepping control for quadrotor developed with neural nets mainly to estimate system dynamics in \cite{das_backstepping_2009}. Chattering in the control signal and complexity of developing control law are limitations of backstepping controller. Model Predictive Control (MPC) with friction compensation for mobile robot with inverse kinematics has been proposed in \cite{barreto_design_2014} and the work has been validated experimentally. The main drawback of MPC is the complexity of the optimization algorithm for linear and nonlinear case which takes more time for computations.

   In our work, \Lone adaptive controller will be studied on different classes of systems. Fuzzy-\Lone adaptive controller will be proposed to tackle problems of \Lone adaptive controller in terms of robustness margin and control signal range. Recent study of neuro-adptive control with PPF will be studied to evaluate the main role of PPF. MRAC with PPF will be proposed to tackle problems of neuro-adaptive control with PPF and \Lone adaptive controller in a proper way. Robust adaptive observer will be implemented with \Lone adaptive controller to examine the performance under inaccessible states. All foregoing tools will be applied on different classes of high nonlinear systems including UVS. Moreover, the nonlinearities will be assumed to be unknown with uncertainties in parameters.

\section{\Lone Adaptive Controller} 

    \Lone adaptive control was first inspired from MRAC. MRAC has been developed initially to control linear systems with uncertainty in parameters \cite{parks_liapunov_1966}. MRAC stability performance relies on Lyapunov function.

    \Lone adaptive controller has been built to enable fast adaption and ensuring robustness. \Lone adaptive controller ensures uniformly bounded in the transient response and  steady state tracking for both regulated output and control signal owing to the low pass filter in the feedback loop. Through the use of low pass filter in the feedback loop will increase the adaptation gain, \Lone adaptive control has been proposed to solve several issues that may exist in the control design. Output of the actual system will be compared to the output of the predicted system and the difference will be addressed into the projection function to help in estimating the uncertainties and disturbances. The output of the projection function will be used in building the required control signal. \Lone adaptive controller design could be adopted to control linear and nonlinear systems with uncertainties in both dynamics and input parameters in the presence of disturbances.

    Nonlinearities, uncertainties, disturbances and unmodelled input will be represented by compact regions and all these regions will give a complete view of system nonlinearities. The major advantage of \Lone adaptive controller is that the worst scenario of all previous unexact modeling can be represented by compact regions with upper and lower bounds without accurate knowledge of nonlinearities structure. \Lone adaptive controller can be defined as a robust controller for improving the transient and tracking response with appropriate assumptions of foregoing compact regions. All previous approximations have to be concerned to build approximated model allows us to build \Lone adaptive controller with satisfactory performance.

   \Lone adaptive controller has been proposed successfully for a simple SISO system in \cite{cao_design_2006}. In this work, the controller and stability analysis was mainly designed for an unstable linear system with constant uncertain parameters in the level of the states which assumed to be unknown. The output response shows a satisfactory transient and tracking performance with different values of a step input. In the following year, The work has been modified including control law and stability analysis in order to be able to deal with nonlinear time varying unknown uncertainties and disturbances for nonlinear SISO systems \cite{cao_guaranteed_2007}. The output performance of shows good results for both tracking, transient response and smooth control signal. Therefor, the controller has been tested on the same nonlinear system and with higher level of time varying uncertainties. Although, the output performance showed good results similar to previous case, the control signal included chattering in contrast to the first case. Finally, the work has been formulated in the following year as a journal paper \cite{cao_design_2008} considering the foregoing two cases SISO systems in \cite{cao_design_2006,cao_guaranteed_2007} in addition to the investigation of different feedback filter structures.

   \Lone adaptive control for nonlinear systems with unmatched uncertainties has been formulated in \cite{gregory_l1_2009} for NASA AIRSTAR flight. It was designed for single flight condition and data recorded during flight test and compared to simulated output data. The comparison study showed satisfactory results and good flight control although results were not very close due to insufficient representations of nonlinearities, disturbances and unmodeled input in the control law.

   \Lone adaptive controller was successfully designed for high nonlinear SISO systems \cite{luo_l_2010}. The control law formulation considered nonlinear time variant for each of uncertainties, system nonlinearities and disturbances in addition to unmodeled input parameters. The controller performance has been validated on high nonlinear SISO system including nonlinearities in the input signal. The transient and tracking performance showed great results with cosine reference input. The same procedure can be applied on MIMO nonlinear systems.

   \Lone adaptive controller for MIMO nonlinear systems in the presence of strong coupling and unmatched uncertainties has been proposed successfully in \cite{xargay_l1_2010}. The work in \cite{xargay_l1_2010} approximated the system into two parts where the first was matched and the second was unmatched part. The control law was developed successfully and stability analysis ensured the robustness of the proposed controller. The output performance showed impressive results for tracking capabilities.

   \Lone adaptive control has been tested for different applications and specifically for flight tests in \cite{gregory_l1_2009,leman_l1_2010,kaminer_coordinated_2007,michini_l1_2009,wang_novel_2008,kharisov_l1_2008} where it shows promising results with flight applications. It has been formulated for different aspects of control problems in \cite{hovakimyan_l1_2010}.  The structure of \Lone adaptive control theory depends on three features and one of them is the implementation of a low pass filter in order to limit the frequency range of the control signal and reduce the effect of uncertainties. The low pass filter should be selected such that the system output tracks properly the reference input and the undesirable uncertainties and frequencies are filtered \cite{cao_design_2006,hovakimyan_l1_2010}. Using the low pass filter, \Lone ensures decoupling between robustness, fast adaptation, infinity norm boundedness of the transient and steady state responses.
   
   The optimal structure of filter has been studied extensively in \cite{hovakimyan_l1_2010} by investigating different type of structures and identifying the optimal filter coefficients. Indeed, the determination of the appropriate parameters of the best filter within a certain class of predefined structure has attracted a particular attention and several attempts on identifying these optimal coefficients have been made. This includes convex optimization based on linear matrix inequality \cite{hovakimyan_l1_2010,li_filter_2008} and multi-objective optimization using MATLAB optimization solver \cite{li_optimization_2007}. Limitations of \Lone adaptive controller and the interconnection between adaptive estimates and the feedback filter were studied in \cite{kharisov_limiting_2011}, where Several filter designs were considered based on disturbance observer. More recent, Systematic approach was presented in \cite{kim_multi-criteria_2014} to determine the optimal feedback filter coefficients in order to increase the zone of robustness margin. The authors proposed the use of greedy randomized algorithms during the analysis of the system performance and robustness in the presence of uncertainties.  
     
   The trade-off between fast desired closed loop dynamics and filter parameters relies on error values. However, all previous studies assume constant coefficients of the feedback filter and the effort of tuning the filter's parameters is performed off-line. Increasing the bandwidth of the low pass filter will reduce robustness margin, which will require slowing the desired closed loop performance in order to regain the robustness. However, slower selection of desired closed loop performance will deteriorate the output performance especially during the transient period \cite{hovakimyan_l1_2010}. We argue that increasing the robustness with fast closed loop dynamics requires dynamic on-line tuning of the feedback filter gain. The method should practical and implementable. Therefore, in this thesis, we propose a fuzzy tuning of the filter coefficients function based on the rate and value of the tracking error between the model output and the system output. 

\section{Adaptive Control with Prescribed Performance Function} 

    Prescribed performance is considered as convergence the tracking error into an arbitrarily small residual set and the convergence error should be within range. Prescribed performance with robust adaptive control will provide a smooth control signal for soft tracking. It comes to solve the problem of accurate computation of the upper bounds for systematic convergence owing to nonexistence adaptive control nonlinear systems for error convergence into a predefined small set.

    The main function of the prescribed performance is the ability of tracking the error into a defined small set.  Prescribed performance should guarantee many factors
    \begin{itemize}
      \item The convergence has to be less than a prescribed value.
      \item Maximum overshot is sufficiently less than small prescribed value.
      \item Uniform ultimate boundedness property for the transformed output error.
      \item Adaptive and smooth tracking.
    \end{itemize}

   Several studies included in their design the use of PPF with linearly parameterized neural network as approximation model to handle unknown nonlinearities and disturbances with or without fuzzy techniques \cite{bechlioulis_robust_2008,bechlioulis_adaptive_2009,wang_verifiable_2010,kostarigka_adaptive_2012,na_adaptive_2013,na_adaptive_2014}. PPF has been applied in different applications and showed promising results. It was first introduced with neuro-adaptive control feedback for strict MIMO systems with unknown nonlinearities; linearly parameterized neural network has been used to approximate the model \cite{bechlioulis_robust_2008}. Although the control law prove robust performance and track the output performance into the desired trajectory, defining radial basis neural network weights offline by try and error is considered the main drawback in \cite{bechlioulis_robust_2008}. In addition, values of other constant parameters are sensitive. Overall, the output performance showed great results for 2-DOF planar robot.

   Robust adaptive controller with prescribed performance has been modified to deal with uncertain MIMO nonlinear systems \cite{kostarigka_adaptive_2012}. Linearly parameterized neural network has been used to compute the control signal and avoid the need of observer from the measured output. Although output performance proves robustness and control law refers to system stability, but limitations of \cite{bechlioulis_robust_2008} still exist in \cite{kostarigka_adaptive_2012}. Also, \cite{kostarigka_adaptive_2012} mentioned another flaw that even structure of each neuron in the neural network will be defined by try and error.

   SISO system with unknown nonlinearities for strict feedback systems studied in \cite{bechlioulis_prescribed_2010}. The work in \cite{bechlioulis_prescribed_2010} is mostly similar to that in \cite{bechlioulis_robust_2008} and the only difference was the way of developing control law. The output showed good performance and it had same limitations of \cite{bechlioulis_robust_2008}. Adaptive compensation control for uncertain nonlinear strict feedback systems with constrained input proposed in \cite{chen_adaptive_2011}. The control law mainly based on two adaptive backstepping controller with prescribed performance bound. Adaptive control with PPF has been proposed for nonlinear systems with unknown dead zone and in order to compensate nonlinearities and uncertainties in the system \cite{na_adaptive_2013}. In \cite{sun_fuzzy_2014}, A fuzzy adaptive prescribed performance control for MIMO uncertain chaotic systems is presented. The system is in a non-strict feedback form. A proportional integral adaptation law is proposed for updating the parameters of the fuzzy logic controller.

\section{Adaptive Observers} 
   Adaptive observer design is an active area of research and it was studied extensively for linear time invariant SISO systems in \cite{narendra_stable_2012,ioannou_stable_1995}. Robust observer for uncertain linear systems with solution provided by algebraic Riccati equation presented in \cite{gu_robust_2001}. Generally, sliding mode observers such as \cite{davila_second-order_2005,qiao_new_2013} are suitable with certain model structures. Neural network has been studied widely for observer design and showed efficacy in observing system states. Radial Basis Function (RBF) in \cite{lee_improved_2007,stepanyan_adaptive_2007} and Chebyshev neural network observer in \cite{shaik_real-time_2011} are designed as adaptive observers for nonlinear systems. Try and error are significant problem in adaptive Neural Network (NN) observer design in addition to the need of multi layers in certain cases. Adaptive observer design for nonlinear uncertain systems has been proposed in \cite{marine_robust_2001,liu_robust_2009}. The advantage of \cite{liu_robust_2009} is being effective for unmodeled dynamics in addition to the possibility of building the adaptation law of observer in the absence of control signal knowledge.

\section{Conclusions} 
   This chapter included overview of adaptive control research also included several research works on nonlinear systems especially UVS. The main work of research focused on \Lone adaptive controller and neuro-adaptive control with PPF. The main contribution in this work has been presented.

\clearpage

\newpage


\chapter{\Lone ADAPTIVE CONTROLLER}
\label{yyyyy3}

\section{Introduction}
  This chapter investigates the transient and tracking performance of \Lone adaptive controller on nonlinear systems with different structures. The control signal will be evaluated with respect to the foregoing features. The controller structure, stability analysis as well as simulations will be presented. The trade-off between fast closed loop dynamics and filter coefficients will be examined. The chapter consists of five sections with first section includes an introduction. The second section discusses \Lone adaptive controller for uncertain SISO systems. The third section handles \Lone adaptive controller for uncertain MIMO systems. The fourth section presents \Lone adaptive controller for uncertain MIMO systems in the presence of strong coupling and unmatched uncertainties. Finally, we conclude in the last section.
  
\section{\Lone Adaptive Controller for Uncertain SISO Systems}
\subsection{Problem Formulation}
  Consider the following class of systems:
  \begin{equation}
    \label{eq:ch2L1Act}
    \begin{aligned}
      &\dot{x}\left(t\right) = A_{m}x\left(t\right) +B(\omega u\left(t\right) + \theta^{\top} x\left(t\right)+\sigma\left(t\right) )\\
      & y\left(t\right) = Cx\left(t\right)
    \end{aligned}
  \end{equation}

   where $x\left(t\right) \in \mathbb{R}^n$ is the system state vector (measured); $u\left(t\right) \in \mathbb{R}$ the control input; $y\left(t\right) \in \mathbb{R}$ is the system output; $B$ and $C$ are constant matrices (known);  $A_m \in \mathbb{R}^{n \times n}$ is Hurwitz  matrix (known) and refers to the desired closed-loop dynamics; $\omega \in \mathbb{R}$ is an unknown matrix with known sign; $\theta\left(t\right) \in \mathbb{R}^n$ is a vector of time-varying unknown parameters; and $\sigma\left(t\right) \in \mathbb{R}$ models input disturbances.
  
   \begin{assumption}
       (The control input is partially known with known sign)
      Let the upper and lower input gain bounds be defined by $\omega_l$ and $\omega_u$ respectively, where
       \begin{equation*}
       \omega \in \Omega \triangleq [\omega_l , \omega_u],\hspace{10pt} |\dot{\omega}| < d_{\omega}
       \end{equation*}
      where $\Omega$ is assumed to be known convex compact set and $0<\omega_l<\omega_u$ are uniformly known conservative bounds.
      \label{assum31}
   \end{assumption}
   \begin{assumption} 
   (Unknown parameters are uniformly bounded)
   Let $\Theta$, $\Delta_0$ be known convex compact where $\Theta, \Delta_0 \in \mathbb{R}^+ $ are known (conservative) bound of $\theta$ and $\sigma$ where
    \begin{equation*}
      \theta\left(t\right) \in \Theta, \hspace{10pt} |\sigma\left(t\right)| \in \Delta_0, \hspace{10pt}\forall t \geq 0
    \end{equation*}
    \label{assum32}
   \end{assumption}
   \begin{assumption}
       (Partial derivatives are semiglobal uniformly bounded)
   Let $\theta\left(t\right)$ and $\sigma\left(t\right)$ be continuously differentiable with $\dot{\theta}$ and $\dot{\sigma}$ they are bounded by $d_\theta$ and $d_\sigma$ where
    \begin{equation*}
      ||\dot{\theta}|| \leq d_\theta< \infty, \hspace{10pt} ||\dot{\sigma}|| \leq d_\sigma< \infty, \hspace{10pt}\forall t \geq 0
    \end{equation*}
    \label{assum33}
   \end{assumption}
   The work in this section aims at designing a full-state feedback adaptive controller to ensure that $y\left(t\right)$ tracks a given bounded piecewise-continuous reference signal $r\left(t\right)$ with quantifiable performance bounds.We will apply the controller on many case studies to evaluate the output performance in terms of transient and tracking response and the control signal in terms of smoothness and boundedness.
   
  \subsection{$\mathcal{L}_1$ Adaptive Control Architecture}
    {\bf State predictor:} We consider the following state predictor:
     \begin{equation}
     \label{eq:ch2L1est}
     \begin{aligned}
       &\dot{\hat{x}}\left(t\right) = A_{m}\hat{x}\left(t\right) + B(\hat{\omega} u\left(t\right) + \hat{\theta}^{\top} x\left(t\right)+\hat{\sigma} )\\
       & \hat{y}\left(t\right) = C\hat{x}\left(t\right)
     \end{aligned}
   \end{equation}
   
\begin{figure}[h!]
   \centering
   \includegraphics[height=7cm, width=14cm]{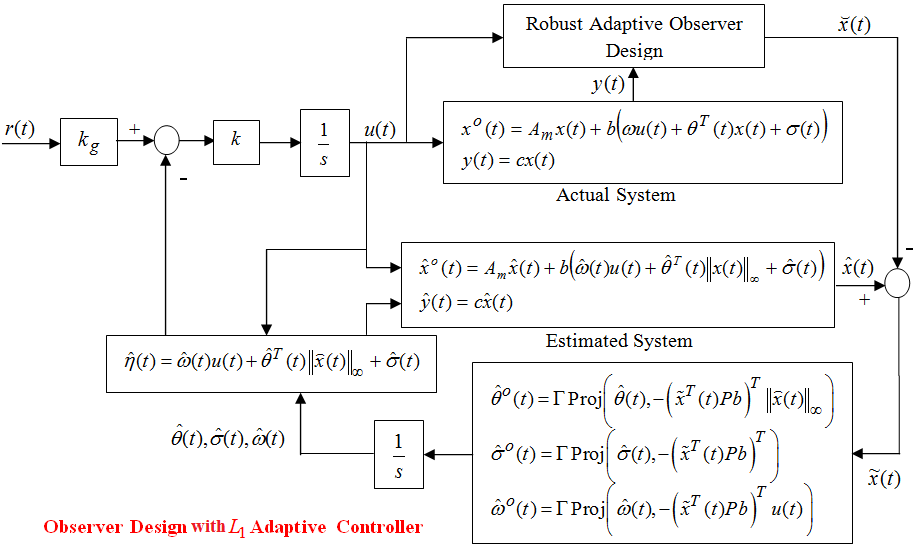}
   \caption{Closed loop \Lone adaptive control system.}
   \label{fig:ch2L1_Sys1Lin}
  \end{figure}

  The state predictor has the same structure as defined in \eqref{eq:ch2L1Act} except that the unknown parameters $\omega$, $\theta\left(t\right)$, and $\sigma\left(t\right)$ are being replaced by their adaptive estimates $\hat{\omega}$, $\hat{\theta}\left(t\right)$ and $\hat{\sigma}\left(t\right)$.
   \begin{equation}
     \label{eq:ch2L1Proj}
     \begin{aligned}
       &\dot{\hat{\omega}} = \Gamma Proj(\hat{\omega},-\tilde{x}^{\top}Pbu\left(t\right)), \hspace{10pt} \hat{\omega}(0) = \hat{\omega}_0  \\
       &\dot{\hat{\theta}} = \Gamma Proj(\hat{\theta},-\tilde{x}^{\top}Pbx\left(t\right)) \hspace{10pt} \hat{\theta}(0) = \hat{\theta}_0\\
       &\dot{\hat{\sigma}} = \Gamma Proj(\hat{\sigma},-\tilde{x}^{\top}Pb) \hspace{10pt} \hat{\sigma}(0) = \hat{\sigma}_0
     \end{aligned}
   \end{equation}
  where $\tilde{x} \triangleq \hat{x} - x\left(t\right)$, $\Gamma \in \mathbb{R}^{+}$ is the adaptation gain, and $P=P^{\top}>0$ is defined by solving the algebraic Lyapunov equation $A_m^{\top}P+PA_m=-Q$ for arbitrary symmetric $Q=Q^{\top}>0$. The projection operator ensures that $\hat{\omega} \in \Omega_0 \triangleq [\omega_l , \omega_u]$, $\hat{\theta} \in \Theta \triangleq [-\theta_b,\theta_b]$, $|\hat{\sigma}| \leq \Delta_0$,  while $\Omega_0$ and $\Delta_0$ are being replaced by $\Omega$ and $\Delta$ to satisfy
    \begin{equation*}
      \Omega_0<\Omega, \hspace{10pt} \Delta_0<\Delta,
    \end{equation*}
  {\bf Control Law:} Control signal can be calculated as follows
  \begin{equation}
    u(s) = -kD(s)(\hat{\eta}(s)-k_gr(s))
  \end{equation}
  where $r(s)$ and $\hat{\eta}(s)$ are the Laplace transforms of $r\left(t\right)$  and $\hat{\eta}\left(t\right) = \hat{\omega} u\left(t\right) + \hat{\theta} x\left(t\right)+\hat{\sigma}$ respectively; and the necessary feedforward gain in order to get unity steady state may be calculated  by $k_g \triangleq -1/(CA_{m}^{-1}B)$ ; Both of the feedback gain $k > 0$ and a strictly proper transfer function $D(s)$ will lead to a strictly proper stable closed loop system.
  \begin{equation}
    C(s)\triangleq \frac{\omega k D(s)}{1 + \omega k D(s)}, \hspace{10pt} \forall \omega \in \Omega_0
  \end{equation}
  with DC gain $C(0) = 1$. One simple choice is $D(s) = 1/s$, which yields a first-order strictly proper $C(s)$ of the form
  \begin{equation*}
    C(s)\triangleq \frac{\omega k }{s + \omega k }
  \end{equation*}
  Let
  \begin{equation}
    L \triangleq \max_{\theta\in\Theta} ||\theta\left(t\right)||_{\mathcal{L}_1}, \hspace{10pt} H(s)=(sI-A_m)^{-1}b,\hspace{10pt} G(s) \triangleq H(s)(1-C(s))
  \end{equation}
  Then the \Lone norm of \Lone adaptive controller will be
  \begin{equation*}
    ||G(s)||_{\mathcal{L}_1} L \leq 1
  \end{equation*}  
  Figure ~\ref{fig:ch2L1_Sys1Lin} shows the structure of closed loop \Lone adaptive controller for uncertain SISO systems.

\subsection{\Lone Adaptive Control Stability Analysis}
  {\bf Transient and Steady-State Performance:} The error dynamics between system dynamics in \eqref{eq:ch2L1Act} and state predictor in \eqref{eq:ch2L1est} can be written as
  \begin{equation}
    \label{eq:ch2L1err1}
      \dot{\tilde{x}}\left(t\right) = A_{m}\tilde{x}\left(t\right) + b(\tilde{\omega} u\left(t\right) + \tilde{\theta}^{\top}x\left(t\right) + \tilde{\sigma}\left(t\right) ) = A_{m}\tilde{x}\left(t\right) + b\tilde{\eta}\left(t\right)\\
  \end{equation}
  Where $\tilde{x} = \hat{x} - x$, $\tilde{\theta} = \hat{\theta} - \theta$, $\tilde{\omega} = \hat{\omega} - \omega$ and $\tilde{\sigma} = \hat{\sigma} - \sigma$. The nonlinear part is $\tilde{\eta}\left(t\right)$ and its Laplace transform  $\tilde{\eta}(s)$ where $\tilde{\eta}\left(t\right) \triangleq \tilde{\omega} u\left(t\right) + \tilde{\theta}^{\top}x\left(t\right) + \tilde{\sigma}\left(t\right)$. The Laplace transform of the error dynamics in \eqref{eq:ch2L1err1} can be rewritten as
  \begin{equation}
      \tilde{x}\left(t\right) = (sI-A_{m})^{-1}B\tilde{\eta}(s)=H(s)\tilde{\eta}(s)
  \end{equation}
   \begin{lemma}
   The prediction error $\tilde{x}\left(t\right)$ is uniformly bounded,
   \end{lemma}
  \begin{equation}
    ||\tilde{x}||_{\infty} \leq \sqrt{\frac{\theta_m}{\lambda_{min}(P)\Gamma}}
  \end{equation} 
  where
  \begin{equation}
  \label{eq:ch2thetam1}
    \theta_m \triangleq \max_{\theta\in\Theta}||\theta||^2 + 4\Delta^{2} + (\omega_l-\omega_u)^2 + 4\frac{\lambda_{max}(P)}{\lambda_{min}(Q)}(d_{\theta}\max_{\theta\in\Theta}||\theta|| + d_{\sigma}\Delta)
  \end{equation}
  which will be verified as follows.\\
  {\bf Stability proof:} Consider the Lyapunov function candidate
  \begin{equation}
    \label{eq:ch2VPart1}
    V(\tilde{x},\tilde{\theta},\tilde{\omega},\tilde{\sigma}) = \tilde{x}^{\top}P\tilde{x} + \frac{1}{\Gamma}(\tilde{\theta}^{\top}\tilde{\theta}+\tilde{\omega}^{\top}\tilde{\omega} +\tilde{\sigma}^{\top}\tilde{\sigma}) 
  \end{equation}
  Since $\hat{x}(0) = x(0)$ then we can verify that
  \begin{equation*}
    V(0) \leq  \max_{\theta\in\Theta}||\theta||^2 + 4\Delta^{2} + (\omega_l-\omega_u)^2 \leq \frac{\theta_m}{\Gamma}
  \end{equation*}
  \begin{equation*}
    \begin{split}
    \dot{V} \leq & \tilde{x}^{\top}P\dot{\tilde{x}} +\dot{\tilde{x}} ^{\top}P\tilde{x} + \frac{1}{\Gamma}(\tilde{\theta}^{\top}\dot{\hat{\theta}} +\dot{\hat{\theta}}^{\top}\tilde{\theta}+\tilde{\sigma}^{\top}\dot{\hat{\sigma}} +\dot{\hat{\sigma}}^{\top}\tilde{\sigma} +\tilde{\omega}^{\top}\dot{\tilde{\omega}}+\\
     & \dot{\tilde{\omega}}^{\top}\tilde{\omega} - \tilde{\theta}^{\top}\dot{\theta} -\dot{\theta}^{\top}\tilde{\theta} - \tilde{\sigma}^{\top}\dot{\sigma} - \dot{\sigma}^{\top}\tilde{\sigma})
    \end{split}
  \end{equation*}
  \begin{equation*}
    \begin{split}
    \dot{V} \leq & \tilde{x}^{\top}Q\tilde{x} + \frac{2}{\Gamma}(\dot{\hat{\theta}} + \tilde{x}^{\top}PBx) +\frac{2}{\Gamma}(\dot{\hat{\sigma}} + \tilde{x}^{\top}PB) +\frac{2}{\Gamma}(\dot{\hat{\omega}} + \tilde{x}^{\top}PBu)\\
    &  - \frac{2}{\Gamma}(\tilde{\theta}^{\top}\dot{\theta} + \tilde{\sigma}^{\top}\dot{\sigma})
    \end{split}
  \end{equation*}
  \begin{equation}
  \label{eq:ch2_vodt1_1}
    \dot{V} \leq -\tilde{x}^{\top}Q\tilde{x} + \frac{2}{\Gamma}\big(|\dot{\theta}^{\top}\tilde{\theta}|  +|\dot{\sigma}^{\top}\tilde{\sigma}|)
  \end{equation}
  As mentioned in Assumption \ref{assum31}, \ref{assum32} and \ref{assum33}, the projection operator ensures that $\theta\left(t\right)\in\Theta$, $|\sigma\left(t\right)|\in\Delta$ for all $t \geq 0$, and therefore, the upper bounds in assumption \ref{assum32} lead to the following upper bound:
  \begin{equation}
   \tilde{\theta}^{\top}\dot{\theta} + \tilde{\sigma}^{\top}\dot{\sigma} \leq 2(d_{\theta}\max_{\theta\in\Theta}||\theta|| + d_{\sigma}\Delta)
  \end{equation}
  Moreover, the projection operator also ensures that
  \begin{equation}
   \label{eq:ch2_v1_1}
    \max_{t \geq 0} \big(\frac{1}{\Gamma}(\tilde{\theta}^{\top}\tilde{\theta} + \tilde{\omega}^{\top}\tilde{\omega} + \tilde{\sigma}^{\top}\tilde{\sigma})\big) \leq  \frac{1}{\Gamma}(\max_{\theta\in\Theta}||\theta||^2 + 4\Delta^{2} + (\omega_l-\omega_u)^2)
  \end{equation}
  which holds for all $t \geq 0$.
  If at any time $t_1>0$, one has $V(t_1) \geq \theta_m/\Gamma$, then it follows from \eqref{eq:ch2thetam1} and \eqref{eq:ch2VPart1} that
  \begin{equation}
   \label{eq:ch2_v1_2}
    \tilde{x}^{\top}(t_1)P\tilde{x}(t_1) > 4\frac{\lambda_{max}(P)}{\lambda_{min}(Q)}(d_{\theta}\max_{\theta\in\Theta}||\theta|| + d_{\sigma}\Delta)
  \end{equation}
  and thus
  \begin{equation}
  \label{eq:ch2_vodt1_2}
    \tilde{x}(t_1)^{\top}Q\tilde{x}(t_1) \geq \tilde{x}^{\top}(t_1)P\tilde{x}(t_1) >\frac{4}{\Gamma}(d_{\theta}\max_{\theta\in\Theta}||\theta|| + d_{\sigma}\Delta)
  \end{equation}
  Hence, if $V(t_1) \geq \theta_m/\Gamma$, then from \eqref{eq:ch2_vodt1_1} and \eqref{eq:ch2_vodt1_2} we have
  \begin{equation}
    \dot{V} \leq 0
  \end{equation}
\subsection{Problem Formulation and Simulation}
\begin{figure}[!]
   \centering
   \includegraphics[height=5cm, width=13cm]{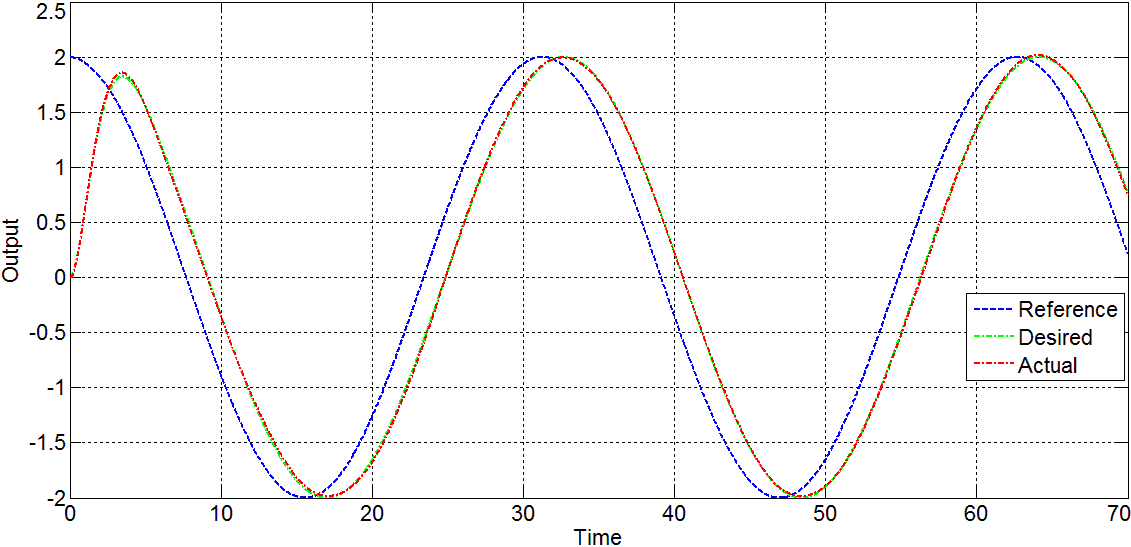}
   \caption{The output performance of \Lone adaptive controller for unknown nonlinear SISO system.}
   \label{fig:Chap_L1_Ex1_1}
  \end{figure}

\begin{figure}[!]
   \centering
   \includegraphics[height=5cm, width=13cm]{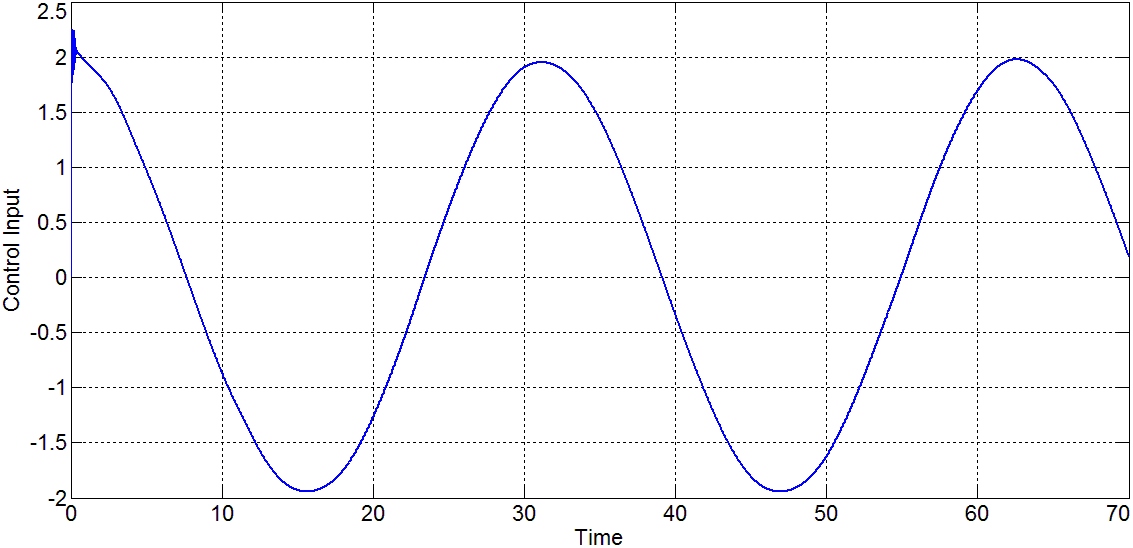}
   \caption{Control signal of \Lone Adaptive controller for unknown nonlinear SISO system.}
   \label{fig:Chap_L1_Ex1_2}
  \end{figure}

\begin{figure}[!]
   \centering
   \includegraphics[height=5cm, width=13cm]{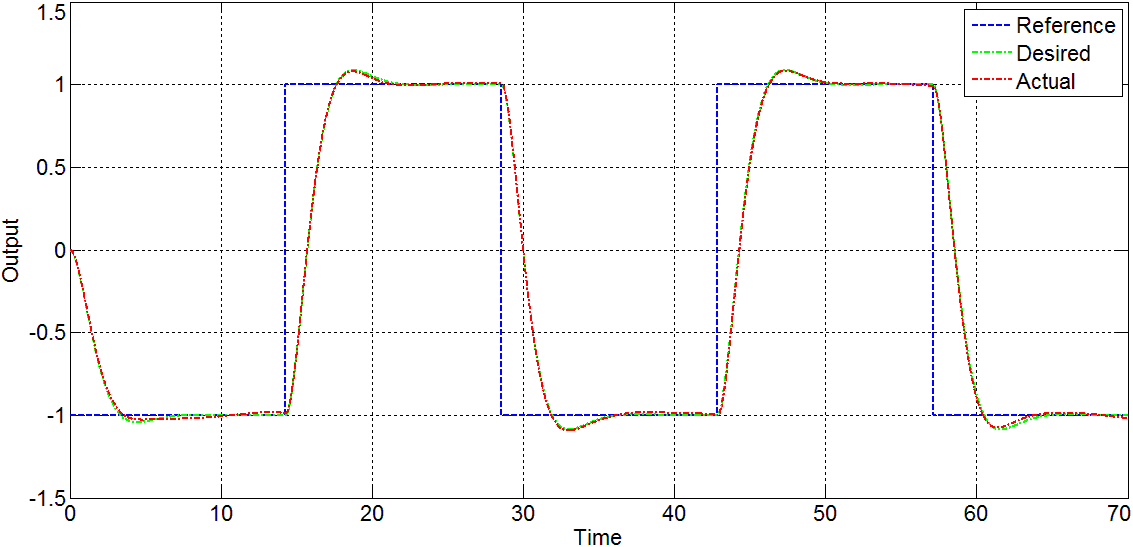}
   \caption{The output performance of \Lone adaptive controller for unknown nonlinear SISO system.}
   \label{fig:Chap_L1_Ex1_3}
  \end{figure}

\begin{figure}[!]
   \centering
   \includegraphics[height=5cm, width=13cm]{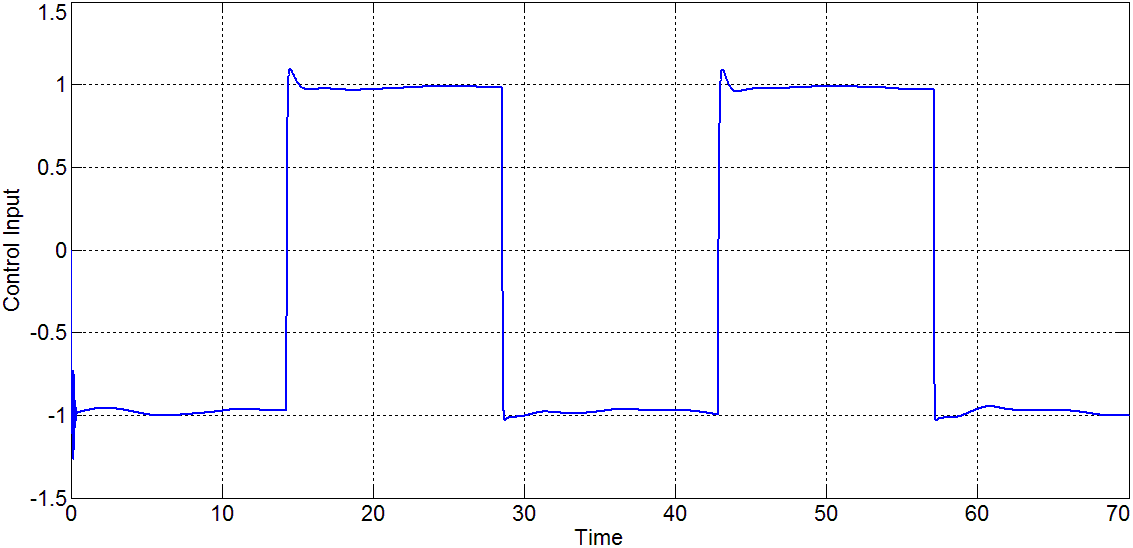}
   \caption{Control signal of \Lone Adaptive controller for unknown nonlinear SISO system.}
   \label{fig:Chap_L1_Ex1_4}
  \end{figure}

    {\bf Example 3.2.1}
      Consider the following unknown nonlinear system \cite{luo_l_2010}
    \begin{equation*}
     \begin{aligned}
       &\dot{x}\left(t\right) = A_{m}x\left(t\right) + b(\omega u\left(t\right) + f(x\left(t\right),u\left(t\right),t))\\
       & y\left(t\right) = cx\left(t\right)
     \end{aligned}
    \end{equation*}
  where $x\left(t\right)=[x_1\left(t\right),x_2\left(t\right)]^{\top}$ are the system states, $u\left(t\right)$is the system control input, $f(x\left(t\right),u\left(t\right),t)$ is assumed to be unknown nonlinear function, $y\left(t\right)$ is the output of the system and the system parameters are presented as following
  \begin{equation*}
    A_m=
    \begin{bmatrix}
     0 & 1\\
     -1 & -1.4
    \end{bmatrix}, \hspace{10pt}
    b=
    \begin{bmatrix}
     0\\
     1
    \end{bmatrix}, \hspace{10pt}
    c=
    \begin{bmatrix}
     1 & 0
    \end{bmatrix}
  \end{equation*}
  \begin{equation*}
  \begin{split}
    f(x\left(t\right),u\left(t\right),t) = & x_1\left(t\right) + 1.4x_2\left(t\right) + (2+0.2sin\left(t\right))u\left(t\right) + sin(u\left(t\right))sin(x_1\left(t\right))\\ 
    & + x_1^2\left(t\right) +  x_2^2\left(t\right) + sin(0.5t)
  \end{split} 
  \end{equation*}
  Parameters of \Lone can be computed numerically and they are chosen to be $\omega_l=0.5$, $\omega_u=3$, $\theta_b = 10$, $\sigma_b = 10$ and the adaptation gain $\Gamma = 100000$. \Lone adaptive control parameters are defined as $Q = \big(\begin{smallmatrix} 1 & 0 \\ 0 & 1 \end{smallmatrix} \big)$, $k = 20$, hence $P = \big(\begin{smallmatrix} 1.4144 & 0.5001 \\ 0.5001 & 0.7144 \end{smallmatrix} \big)$. Figure (\ref{fig:Chap_L1_Ex1_1}) and (\ref{fig:Chap_L1_Ex1_2}) are the output response and control signal respectively with reference input $r\left(t\right) = 2cos(0.2t)$ while figure (\ref{fig:Chap_L1_Ex1_3}) and (\ref{fig:Chap_L1_Ex1_4}) are the output response and control signal respectively with 0.23Hz square wave reference input for the same problem
  
  
\section{\Lone Adaptive Controller for Uncertain MIMO Systems}
\subsection{Problem Formulation}
  Consider in the following class of systems:
    \begin{equation}
     \begin{aligned}
      \label{eq:ch2L1Actf2}
       &\dot{x}\left(t\right) = A_{m}x\left(t\right) + B(\omega u\left(t\right) + f(x\left(t\right),t))\\
       & y\left(t\right) = Cx\left(t\right)
     \end{aligned}
    \end{equation}
   where $x\left(t\right) \in \mathbb{R}^n$ is the system state vector (measured); $u\left(t\right) \in \mathbb{R}^m$ the control input; $y\left(t\right) \in \mathbb{R}^m$is the system output;$b \in \mathbb{R}^{n \times m}$ and $C \in \mathbb{R}^{m \times n}$ are constant matrices (known); $A_m \in \mathbb{R}^{n \times n}$ is Hurwitz matrix (known) and refers to the desired closed-loop dynamics; $\omega \in \mathbb{R}^{m \times m}$ is an unknown constant matrix with known sign; $f(x\left(t\right),t)$ is an unknown nonlinearity.
   
   \begin{assumption} (Uniform boundedness of $f(0,t)$)
      Let $B>0$ such that $ f(0,t) \leq B_l $ for all  $ t \geq 0 $
      \label{assum34}
   \end{assumption}
   \begin{assumption} (Partial derivatives are semiglobal uniformly bounded) 
   For any $\delta >0$, there exist $d_{f_{x}} (\delta) >0 $ and $d_{f_{t}} (\delta) >0 $  such that for arbitrary $||x||_{\infty} \leq \delta$ and any $u$, the partial derivatives of $ f(x\left(t\right),t)) $ is piecewise-continuous and bounded,
       \begin{equation*}
       ||\frac{\partial f(x\left(t\right),t) }{\partial x}|| \leq d_{f_{x}} (\delta),\hspace{10pt} ||\frac{\partial f(x\left(t\right),t) }{\partial t}|| \leq d_{f_{t}} (\delta)
       \end{equation*}
       \label{assum35}
   \end{assumption}
   \begin{assumption}  (Asymptotically stable of initial conditions)
   The system is assumed to start initially with $x_0$ inside an arbitrarily known set $\rho_0$ i.e., $||x_0||_{\infty} \leq \rho_0 < \infty$.
   \label{assum36}
   \end{assumption}
   \begin{assumption} (The control inputs are partially known with known sign) 
   Let upper and lower input gain bounds are defined by $\omega_l$ and $\omega_u$ respectively, where
    \begin{equation*}
    \omega \in \Omega \triangleq [\omega_l , \omega_u],\hspace{10pt} |\dot{\omega}| < \omega
    \end{equation*}
    \label{assum37}
   \end{assumption}
\subsection{\Lone Adaptive Control Architecture}
   {\bf State Predictor:} We consider the following state predictor:
   \begin{equation}
     \label{eq:ch2L1est2}
     \begin{aligned}
       &\dot{\hat{x}}\left(t\right) = A_{m}\hat{x}\left(t\right) + b(\hat{\omega} u\left(t\right) + \hat{\theta} ||x\left(t\right)||_{\infty}+\hat{\sigma} )\\
       & \hat{y}\left(t\right) = c\hat{x}\left(t\right)
     \end{aligned}
   \end{equation}
  The system in \eqref{eq:ch2L1Actf2} can be defined as
     \begin{equation}
       \label{eq:ch2L1Act2}
       \begin{aligned}
         &\dot{\hat{x}}\left(t\right) = A_{m}x\left(t\right) + b(\omega u\left(t\right) + \theta ||x\left(t\right)||_{\infty}+\sigma )\\
         & \hat{y}\left(t\right) = c\hat{x}\left(t\right)
       \end{aligned}
     \end{equation}
  which is similar to \eqref{eq:ch2L1est2} except that the unknown parameters $\omega$, $\theta\left(t\right)$, and $\sigma\left(t\right)$ are being replaced by their adaptive estimates $\hat{\omega}\left(t\right)$, $\hat{\theta}\left(t\right)$ and $\hat{\sigma}\left(t\right)$.
   \begin{equation}
     \label{eq:ch2L1Proj2}
     \begin{aligned}
       &\dot{\hat{\omega}} = \Gamma Proj(\hat{\omega},-\tilde{x}^{\top}Pbu\left(t\right)), \hspace{10pt} \hat{\omega}(0) = \hat{\omega}_0  \\
       &\dot{\hat{\theta}} = \Gamma Proj(\hat{\theta},-\tilde{x}^{\top}Pbx\left(t\right)) \hspace{10pt} \hat{\theta}(0) = \hat{\theta}_0\\
       &\dot{\hat{\sigma}} = \Gamma Proj(\hat{\sigma},-\tilde{x}^{\top}Pb) \hspace{10pt} \hat{\sigma}(0) = \hat{\sigma}_0
     \end{aligned}
   \end{equation}
   where $\tilde{x} \triangleq \hat{x} - x\left(t\right)$, $\Gamma \in \mathbb{R}^{+}$ is the adaptation gain, and $P=P^{\top}>0$ is defined by solving the algebraic Lyapunov equation $A_m^{\top}P+PA_m=-Q$ for arbitrary symmetric $Q=Q^{\top}>0$. The projection operator ensures that $\hat{\omega} \in \Omega_0 \triangleq [\omega_l , \omega_u]$, $\hat{\theta} \in \Theta \triangleq [-\theta_b,\theta_b]$, $|\hat{\sigma}| \leq \Delta_0$,  while $\Omega_0$ and $\Delta_0$ are being replaced by $\Omega$ and $\Delta$ to satisfy
    \begin{equation*}
      \Omega_0 \subset \Omega, \hspace{10pt} \Delta_0 \subset \Delta,
    \end{equation*}
     {\bf Control Law:} Control signal can be calculated as following
     \begin{equation}
       u(s) = -kD(s)(\hat{\eta}(s)-k_gr(s))
     \end{equation}
     where $r(s)$ and $\hat{\eta}(s)$ are the Laplace transforms of $r\left(t\right)$  and $\hat{\eta}\left(t\right) = \hat{\omega} u\left(t\right) + \hat{\theta} x\left(t\right)+\hat{\sigma}$ respectively; and the necessary feedforward gain in order to get unity steady state is calculated  by $k_g \triangleq -1/(CA_{m}^{-1}B)$ ; $k > 0$ is a feedback gain and and $D(s)$ is a strictly proper transfer function such that both of them lead to a strictly proper stable closed loop system.
     \begin{equation}
       C(s)\triangleq \frac{\omega k D(s)}{1 + \omega k D(s)}, \hspace{10pt} \forall \omega \in \Omega_0
     \end{equation}
     with DC gain $C(0) = 1$. One simple choice is $D(s) = 1/s$, which yields a first-order strictly proper $C(s)$ of the form
     \begin{equation*}
       C(s)\triangleq \frac{\omega k }{s + \omega k }
     \end{equation*}
     Let
     \begin{equation}
       L \triangleq \max_{\theta\in\Theta} ||\theta\left(t\right)||_{\mathcal{L}_1}, \hspace{10pt} H(s)=(sI-A_m)^{-1}b,\hspace{10pt} G(s) \triangleq H(s)(1-C(s))
     \end{equation}
     Then the \Lone norm of \Lone adaptive controller will be
     \begin{equation*}
       ||G(s)||_{\mathcal{L}_1} L \leq 1
     \end{equation*}  
     Now, for a given $\rho_0$ as in assumption \ref{assum36}, $k$ and $D(s)$ should be chosen such that there exist $\rho_r>\rho_{in}$ such that the following \Lone norm condition verified
     \begin{equation}
     \label{eq:ch2rhor2}
       ||G(s)||_{\mathcal{L}_1} < \frac{\rho_{r} - ||H(s)C(s)k_g||_{\mathcal{L}_1}||r||_{\mathcal{L}_{\infty}} - \rho_{in}}{L_{\rho_{r}}\rho_{r} + B}
     \end{equation}
     let
     \begin{equation}
       \gamma_{1} \triangleq \frac{||C(s)||_{\mathcal{L}_1}}{1-||C(s)||_{\mathcal{L}_1}L_{\rho_{r}}}\gamma_{0}+ \beta
     \end{equation}
     where $\gamma_{0}$ and $\beta$ are arbitrarily small positive constants.\\
     let
     \begin{equation}
       \rho_{u} \triangleq \rho_{ur} + \gamma_{2}
     \end{equation}
     where $\rho_{ur}$ and $\gamma_{2}$ are defined as following
     \begin{equation}
       \label{eq:ch2rhour2}
       \rho_{ur} \triangleq ||\omega^{-1}C(s)||_{\mathcal{L}_1}(|k_g|||r||_{\mathcal{L}_{\infty}} + L_{\rho_{r}}\rho_{r} + B )
     \end{equation}  
     \begin{equation}
       \gamma_{2} \triangleq ||\omega^{-1}C(s)||_{\mathcal{L}_1}L_{\rho_{r}} \gamma_{1} + ||\omega^{-1}C(s)(c_{0}^{\top}H(s))^{-1}c_{0}^{\top}||_{\mathcal{L}_1}\gamma{0}
     \end{equation}
     and finally let
      \begin{equation}
        \label{eq:ch2L1Max2}
        \theta_b \triangleq d_{f_{x}}(\delta),\hspace{10pt} \Delta \triangleq B + \epsilon
      \end{equation}
    where $\epsilon$ is an arbitrary positive constant.
        
\subsection{\Lone Adaptive Control Stability Analysis}
     {\bf Transient and Steady-State Performance:} The error dynamics between system dynamics in \eqref{eq:ch2L1Act2} and state predictor in \eqref{eq:ch2L1est2} can be written as
     \begin{equation}
       \label{eq:ch2L1err2}
       \dot{\tilde{x}}\left(t\right) = A_{m}\tilde{x}\left(t\right) + b(\tilde{\omega} u\left(t\right) + \tilde{\theta} ||x\left(t\right)||_{\infty}+\tilde{\sigma}\left(t\right) ) = A_{m}\tilde{x}\left(t\right) + b\tilde{\eta}\left(t\right)\\
     \end{equation}
     Where $\tilde{x} = \hat{x} - x$, $\tilde{\theta} = \hat{\theta} - \theta$, $\tilde{\omega} = \hat{\omega} - \omega$ and $\tilde{\sigma} = \hat{\sigma} - \sigma$.The nonlinear part is $\tilde{\eta}\left(t\right)$ and its Laplace transform  $\tilde{\eta}(s)$ where $\tilde{\eta}\left(t\right) \triangleq \tilde{\omega} u\left(t\right) + \tilde{\theta}^{\top}x\left(t\right) + \tilde{\sigma}\left(t\right)$. The Laplace transform of the error dynamics in \eqref{eq:ch2L1err1} can be rewritten as
     \begin{equation}
         \tilde{x}\left(t\right) = (sI-A_{m})^{-1}B\tilde{\eta}(s)=H(s)\tilde{\eta}(s)
     \end{equation}
     Assume
     \begin{equation}
     \label{eq:ch2xrho2}
         ||x\left(t\right)||_{\infty} \leq \rho
     \end{equation}
     \begin{equation}
     \label{eq:ch2urho2}
         ||u\left(t\right)||_{\infty} \leq \rho_u
     \end{equation}
   \begin{lemma}
      The prediction error $\tilde{x}\left(t\right)$ is uniformly bounded,
      \label{lemma32}
   \end{lemma}
     from Lemma \ref{lemma32} and equations \ref{eq:ch2xrho2} and \ref{eq:ch2urho2}, the derivatives of $\omega$, $\theta$ and $\sigma$ are bounded:
     \begin{equation}
     \label{eq:ch2omegadot2}
         |\dot{\omega}| \leq d_{\omega} < \infty
     \end{equation}
     \begin{equation}
     \label{eq:ch2thetadot2}
         |\dot{\theta}| \leq d_{\theta} < \infty
     \end{equation}
     \begin{equation}
     \label{eq:ch2sigmadot2}
         |\dot{\sigma}| \leq d_{\sigma} < \infty
     \end{equation}
     Then we have
     \begin{equation}
       ||\tilde{x}||_{\infty} \leq \sqrt{\frac{\theta_m}{\lambda_{min}(P)\Gamma}}
     \end{equation} 
     where
     \begin{equation}
     \label{eq:ch2thetam2}
       \theta_m \triangleq 4\bigg(\theta_b^2m + \Delta^{2}m + \max_{\omega\in\Omega}tr(\omega^{\top}\omega) + m\frac{\lambda_{max}(P)}{\lambda_{min}(Q)}(d_{\theta}\max_{\theta\in\Theta}||\theta|| + d_{\sigma}\Delta)\bigg)
     \end{equation}
     which will be verified as follows
     {\bf Stability proof:} Consider the Lyapunov function candidate
     \begin{equation}
       \label{eq:ch2VPart2}
       V(\tilde{x},\tilde{\theta},\tilde{\omega},\tilde{\sigma}) = \tilde{x}^{\top}P\tilde{x} + \frac{1}{\Gamma}(\tilde{\theta}^{\top}\tilde{\theta}+tr(\tilde{\omega}^{\top}\tilde{\omega}) +\tilde{\sigma}^{\top}\tilde{\sigma}) 
     \end{equation}
     Since $\hat{x}(0) = x(0)$ then we can verify that
     \begin{equation*}
       V(0) \leq  \frac{4}{\Gamma}\big(\theta_b^2m + \Delta^{2}m + \max_{\omega\in\Omega}tr(\omega^{\top}\omega)\big) \leq \frac{\theta_m}{\Gamma}
     \end{equation*}
     \begin{equation*}
       \begin{split}
       \dot{V} \leq & \tilde{x}^{\top}Q\tilde{x} + \frac{2}{\Gamma}(\dot{\hat{\theta}} + \tilde{x}^{\top}PB||x||_{\infty}) +\frac{2}{\Gamma}(\dot{\hat{\sigma}} + \tilde{x}^{\top}PB) +\frac{2}{\Gamma}(\dot{\hat{\omega}} + \tilde{x}^{\top}PBu)\\
       &  - \frac{2}{\Gamma}(\tilde{\theta}^{\top}\dot{\theta} + \tilde{\sigma}^{\top}\dot{\sigma})
       \end{split}
     \end{equation*}
     \begin{equation}
       \dot{V} = -\tilde{x}^{\top}Q\tilde{x} + \frac{2}{\Gamma}\big(|\dot{\theta}^{\top}\tilde{\theta}|  +|\dot{\sigma}^{\top}\tilde{\sigma}|)
     \end{equation}
     \begin{equation}
     \label{eq:ch2_vodt2_1}
       \dot{V} \leq -\tilde{x}^{\top}Q\tilde{x} + \frac{4}{\Gamma}(d_{\theta}\theta_b + d_{\sigma}\Delta)
     \end{equation}
     As mentioned in Assumption \ref{assum34}, \ref{assum35}, \ref{assum36} and \ref{assum37}, the projection operator ensures that $\theta\left(t\right)\in\Theta$, $|\sigma\left(t\right)|\in\Delta$ for all $t \geq 0$, and therefore, the upper bound in assumption \ref{assum37} lead to the following upper bound:
     \begin{equation}
      \tilde{\theta}^{\top}\dot{\theta} + \tilde{\sigma}^{\top}\dot{\sigma} \leq d_{\theta}\theta_b + d_{\sigma}\Delta
     \end{equation}
     Moreover, the projection operator ensures that
     \begin{equation}
      \label{eq:ch2_v2_1}
       \max_{t \geq 0} \big(\frac{1}{\Gamma}(\tilde{\theta}^{\top}\tilde{\theta}+tr(\tilde{\omega}^{\top}\tilde{\omega}) +\tilde{\sigma}^{\top}\tilde{\sigma}\big) \leq  \frac{1}{\Gamma}(\theta_b^2m + \Delta^{2}m + \max_{\omega\in\Omega}tr(\omega^{\top}\omega))
     \end{equation}
     which holds for all $t \geq 0$.
     If at any time $t_1>0$, one has $V(t_1) \geq \theta_m/\Gamma$, then from \eqref{eq:ch2thetam2} and \eqref{eq:ch2VPart2}, one has
     \begin{equation}
      \label{eq:ch2_v2_2}
       \tilde{x}^{\top}(t_1)P\tilde{x}(t_1) > 4\frac{\lambda_{max}(P)}{\lambda_{min}(Q)}(d_{\theta}\theta_b + d_{\sigma}\Delta)
     \end{equation}
     and thus
     \begin{equation}
     \label{eq:ch2_vodt2_2}
       \tilde{x}(t_1)^{\top}Q\tilde{x}(t_1) \geq \frac{\lambda_{min}(Q)}{\lambda_{max}(P)}\tilde{x}^{\top}(t_1)P\tilde{x}(t_1) >\frac{4}{\Gamma}\frac{\lambda_{max}(P)}{\lambda_{min}(Q)}(d_{\theta}\theta_b + d_{\sigma}\Delta)
     \end{equation}
     Hence, if $V(t_1) \geq \theta_m/\Gamma$, then from \eqref{eq:ch2_vodt2_1} and \eqref{eq:ch2_vodt2_2},
     \begin{equation}
       \dot{V} \leq 0
     \end{equation}

\subsection{Problem Formulation and Simulation}
     {\bf Example 3.3.1} Simulation Problem of Two Link Planar Robot \cite{bechlioulis_robust_2008}
       \begin{equation*}
          M(q)\ddot{q} + C(\dot{q},q)\dot{q} + G_0(q) = \tau
       \end{equation*}
       where $q = [q_1 \hspace{10pt} q_2]^{\top}$ are the angular position and $\tau = [ \tau_1 \hspace{10pt} \tau_2]^{\top}$  are representing the applied torques.\\
       The inertia matrix is represented by 
       \begin{equation*}
          M(q)=
          \begin{bmatrix} 
          M_{11} & M_{12}\\
          M_{21} & M_{22}\\
          \end{bmatrix}
       \end{equation*}
       with\\
       $M_{11} = I_{z_1} + I_{z_2} + \frac{m_1l_1^{2}}{2} + m_2\big(l_1^{2} + \frac{l_2^{2}}{4} + l_1l_2c_2 \big)$\\
       $M_{12} = M_{21} = I_{z_2} + m_2\big( \frac{l_2^{2}}{4} + \frac{1}{2}l_1l_2c_2 \big)$\\
       $M_{22} = I_{z_2} + m_2 \frac{l_2^{2}}{4} $\\
       $C(\dot{q},q)$ is the Coriolis and centrifugal torques matrix, $\dot{q}$ is angular speed and $C(\dot{q},q)\dot{q}$ is actuator joint friction forces where
       \begin{equation*}
          C(\dot{q},q)\dot{q}=
          \begin{bmatrix} 
            c\dot{q_2}+k_1 & -c(\dot{q_1}+\dot{q_2})\\
            c\dot{q_1} & k_2\\
          \end{bmatrix}
          \begin{bmatrix} 
            \dot{q_1}\\
            \dot{q_2}\\
          \end{bmatrix}
       \end{equation*}
       with $c=\frac{1}{2}m_2l_1l_2s_2$.
       and $G_0(q)$  is the vector of gravitational torques
       \begin{equation*}
          G_0(q)=
          \begin{bmatrix} 
            \frac{1}{2}m_1gl_1c_1+m_2g(l_1c_1+\frac{1}{2}l_2c_{12})\\
            \frac{1}{2}m_2gl_2c_{12}\\
          \end{bmatrix}
       \end{equation*}
       with $c1=cos(q_1)$, $c12=cos(q_1+q_2)$, $s1=sin(q_1)$ and $c2=cos(q_2)$.
       Table (\ref{tab:PPFdesc}) and (\ref{tab:PPFparm}) defines the necessary symbols, description and their associated values.
       \begin{table}[h]
       \caption {Description of symbols and their units}
       \label{tab:PPFdesc} 
       \begin{center}
           \begin{tabular}{ l | l | l }
             \hline\hline
             {\bf Symbol} & {\bf Description}  & {\bf Unit} \\ \hline \hline
             $q_i$ & Angular position of joint-$i$ & $rad$\\ \hline
             $\dot{q}_i$ & Angular velocity of joint-$i$ & $rad/sec$\\ \hline
             $\tau_i$ & Applied torque at joint-$i$ & $N/m$\\ \hline
             $m_i$ & Mass of link-$i$ & $kg$\\ \hline
             $l_i$ & Length of link-$i$ & $m$\\ \hline
             $I_{Z_i}$ & Moment Inertia of link-$i$ & $kg.m^2$\\ \hline
             $k_i$ & Friction coefficient of joint-$i$ & $kg.m^2/s$\\ \hline
             $g$ & Gravity acceleration & $m/s^2$\\ \hline \hline
           \end{tabular}
         \end{center}
       \end{table}
     
       \begin{table}[h]
         \caption {System parameters}
         \label{tab:PPFparm} 
         \begin{center}
           \begin{tabular}{ l | l | l | l | l | l | l | l | l }
             \hline\hline
             $m_1$ & $l_1$  & $I_{Z_1}$ & $k_1$ & $m_2$ & $l_2$  & $I_{Z_2}$ & $k_2$ & $g$ \\ \hline
             3.2 & 0.5 & 0.96 & 1 & 2.0 & 0.4 & 0.841 & 1 & 9.81 \\ \hline \hline
           \end{tabular}
         \end{center}
       \end{table}
     
\begin{figure}[!]
     \centering
     \includegraphics[height=5cm, width=13cm]{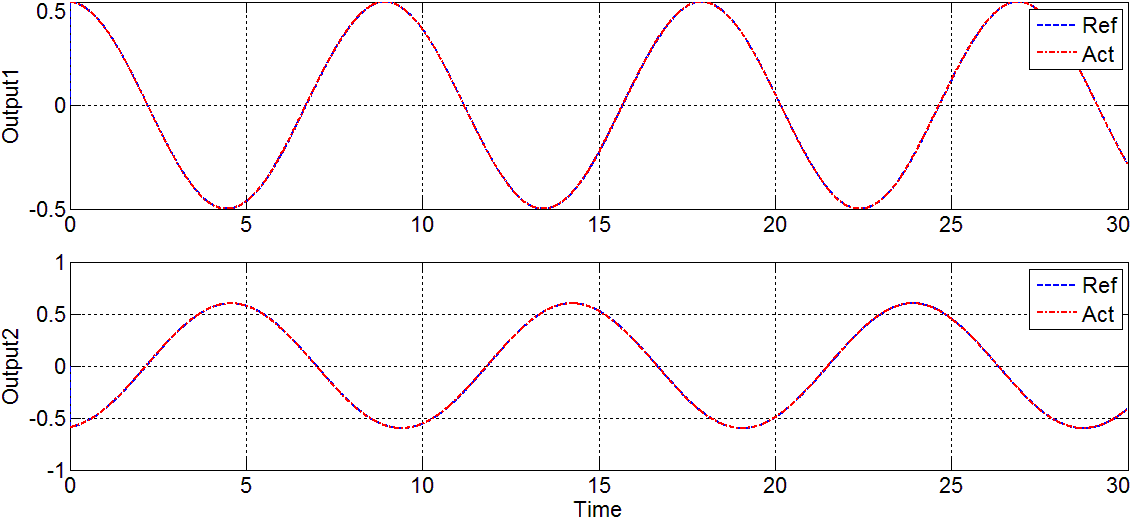}
     \caption{\Lone adaptive control of two link planar robot with reference and actual tracking }
     \label{fig:Chap_L1_Ex2_1}
    \end{figure}

\begin{figure}[!]
     \centering
     \includegraphics[height=5cm, width=13cm]{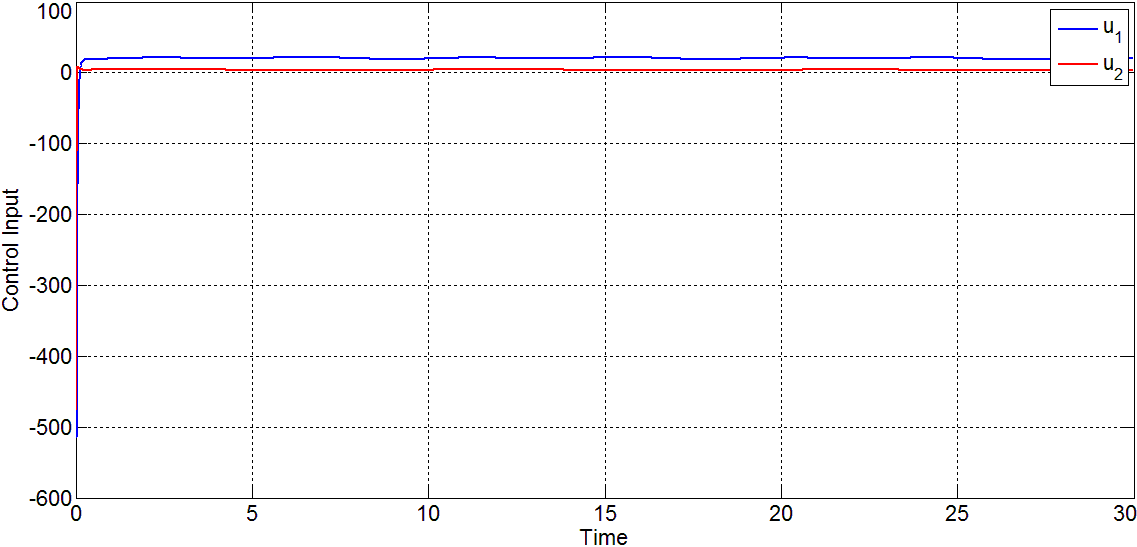}
     \caption{Control signal of \Lone adaptive control for two link planar robot}
     \label{fig:Chap_L1_Ex2_2}
    \end{figure}

       The equation of motion of the nonlinear plant can be represented as following
       \begin{equation*}
          \ddot{q} = -M^{-1}(q)( C(\dot{q},q)\dot{q} + G_0(q) ) + M^{-1}(q)\tau
       \end{equation*}
     {\bf Case 1:} Parameters of \Lone can be computed numerically where their bounds were chosen to be $\omega_l=0.5$, $\omega_u=10$, $\theta_b = 100$, $\sigma_b = 10$ and the adaptation gain $\Gamma = 100000$. Assuming the desired poles are $-300 \pm j5$ and $-400 \pm j5$. The feedback controller was set to be $30diag(4)$. The simulated response will be demonstrated in figure (\ref{fig:Chap_L1_Ex2_1}) and (\ref{fig:Chap_L1_Ex2_2}) for \Lone output performance and control signal respectively.
     
\begin{figure}[!]
     \centering
     \includegraphics[height=5cm, width=13cm]{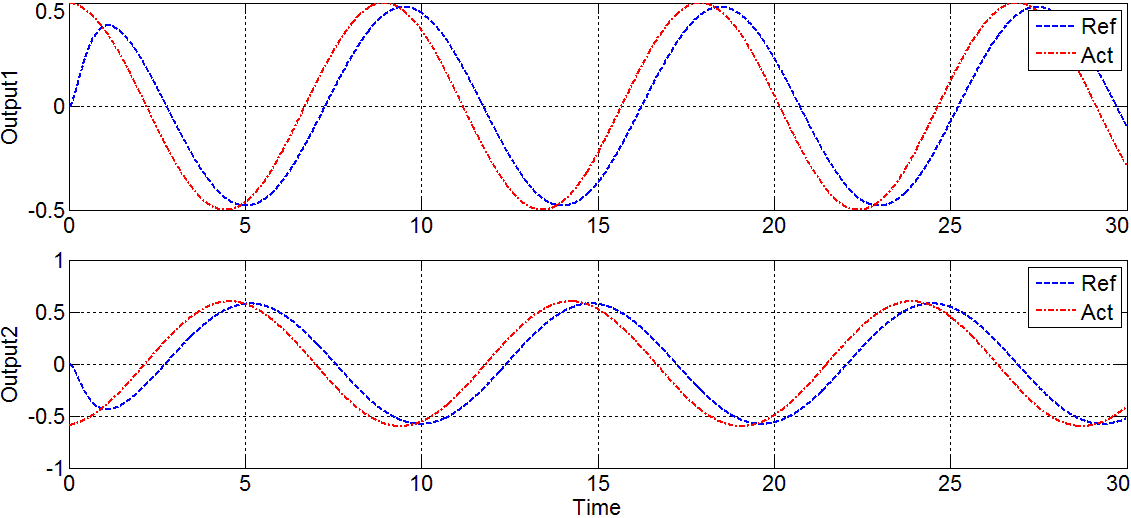}
     \caption{\Lone adaptive control of two link planar robot with reference and actual tracking }
     \label{fig:Chap_L1_Ex2_3}
    \end{figure}

\begin{figure}[!]
     \centering
     \includegraphics[height=5cm, width=13cm]{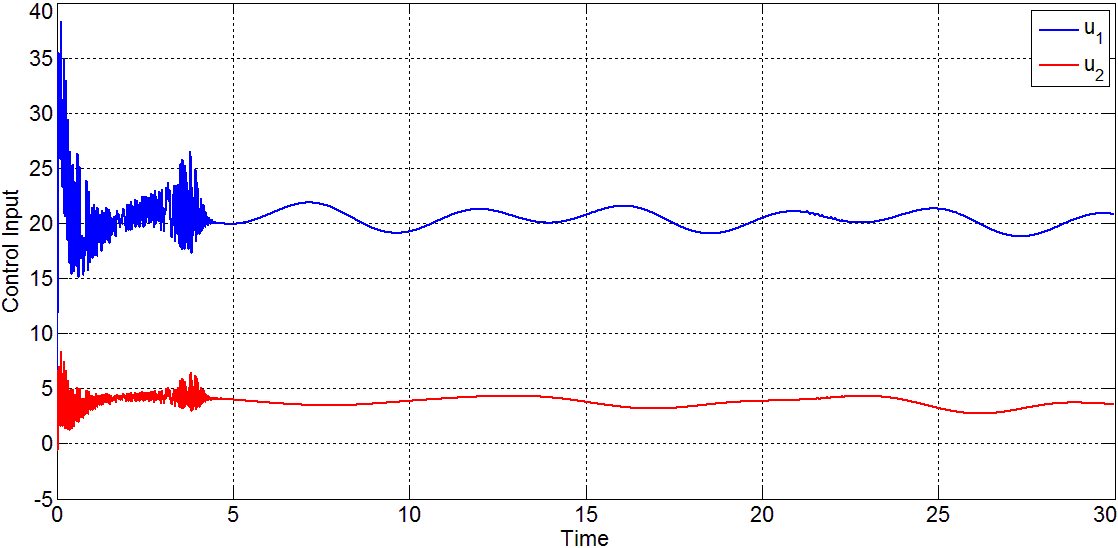}
     \caption{Control signal of \Lone adaptive control for two link planar robot}
     \label{fig:Chap_L1_Ex2_4}
    \end{figure}

    {\bf Case 2:} Figure ~\ref{fig:Chap_L1_Ex2_3}  and ~\ref{fig:Chap_L1_Ex2_4} present the outputs of \Lone adaptive control and control signals respectively considering same assumptions as in case 1 except setting desired poles  $-30 \pm j0.5$ and $-40 \pm j0.5$ in order to investigate the relation between fast and slow desired dynamics with respect to the control signal and tracking performance.
    
    
    Figures (\ref{fig:Chap_L1_Ex2_1}), (\ref{fig:Chap_L1_Ex2_2}), (\ref{fig:Chap_L1_Ex2_3}) and (\ref{fig:Chap_L1_Ex2_4}) describe the relation between robustness and fast tracking response from one hand and control signal range from the other hand. Increasing the speed of transient and tracking performance has a direct relation with how far the desired poles can be located in the left hand side from the origin of $(\sigma-j\omega)$ axis. However, it reduces the robustness of the zone wish demand reducing the feedback gain value. On the other hand, the narrow range of control signal has adverse relation with transient speed.
    
    {\bf Example 3.3.2} Simulation Problem of Quadrotor\\
    Consider the quadrotor model in \cite{das_backstepping_2009} with model parameters presented in \cite{freddi_feedback_2011}
     \begin{equation*}
       \ddot{\eta}_1 = \frac{1}{m}R(\eta_2)\begin{bmatrix} 0 & 0 & \tau_z \end{bmatrix}^{\top} - g\begin{bmatrix} 0 & 0 & 1 \end{bmatrix}^{\top}
     \end{equation*}
     \begin{equation*}
       \ddot{\eta}_2 = f(\eta_2) + G(\eta_2)\begin{bmatrix} \tau_p & \tau_q & \tau_r \end{bmatrix}^{\top}
     \end{equation*}
    Where $R$ is the Euler transformation angle matrix, $\eta_2$ is the Euler angles, $f(\eta_2) \in \mathbb{R}^{3 \times 1}$ is the nonlinear function and $G(\eta_2) \in \mathbb{R}^{3 \times 3}$ is the inverse of the inertia matrix.
    
\begin{figure}[!]
     \centering
     \includegraphics[height=5cm, width=13cm]{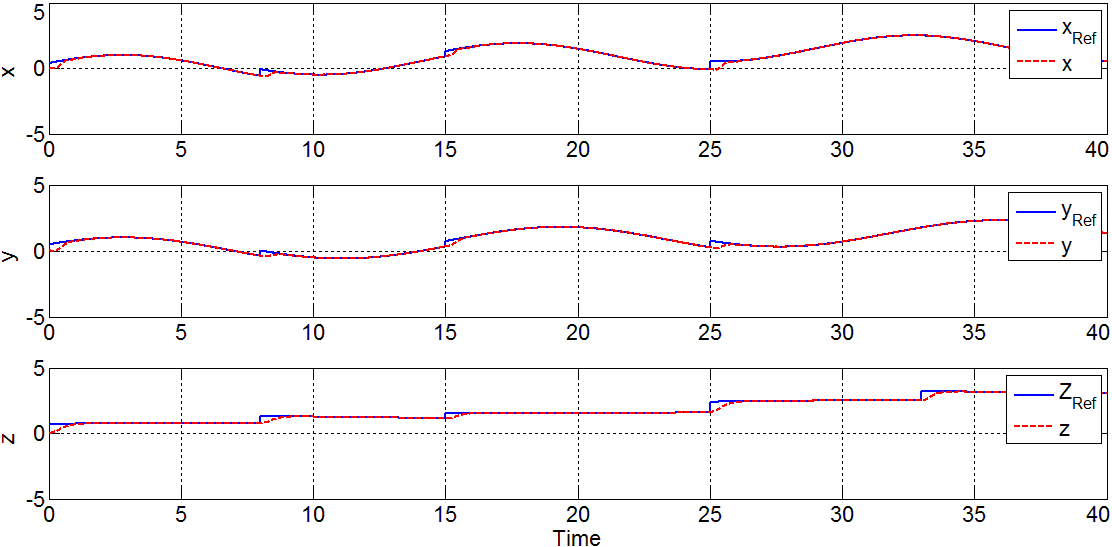}
     \caption{\Lone adaptive controller with reference and actual tracking positions for quadrotor.}
     \label{fig:Chap_L1_Ex3_1pos}
    \end{figure}

\begin{figure}[!]
     \centering
     \includegraphics[height=5cm, width=13cm]{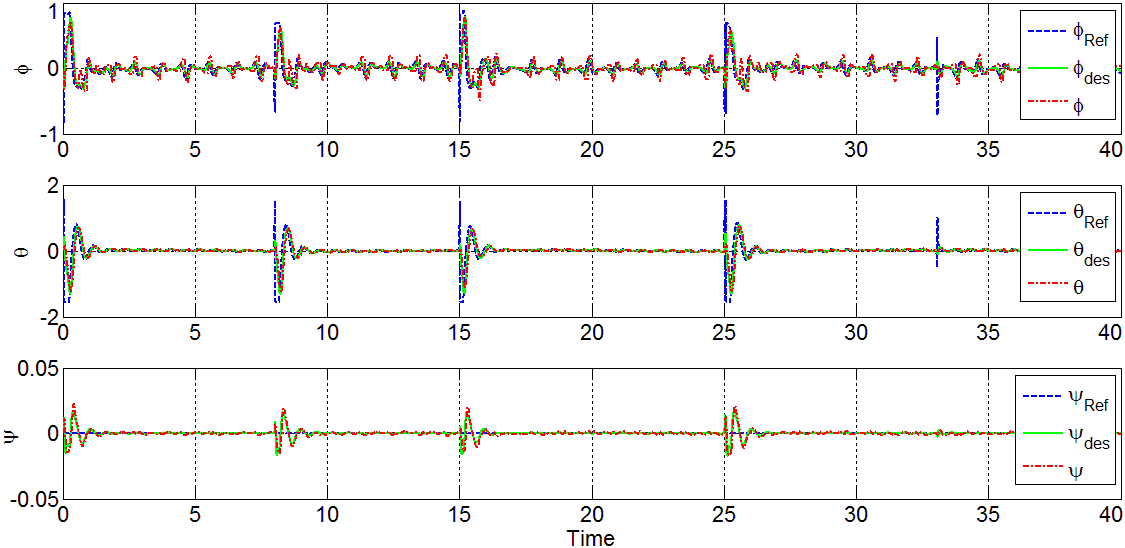}
     \caption{\Lone adaptive controller with reference, desired and actual tracking angles of a quadrotor system.}
     \label{fig:Chap_L1_Ex3_1ang}
    \end{figure}

\begin{figure}[!]
     \centering
     \includegraphics[height=5cm, width=13cm]{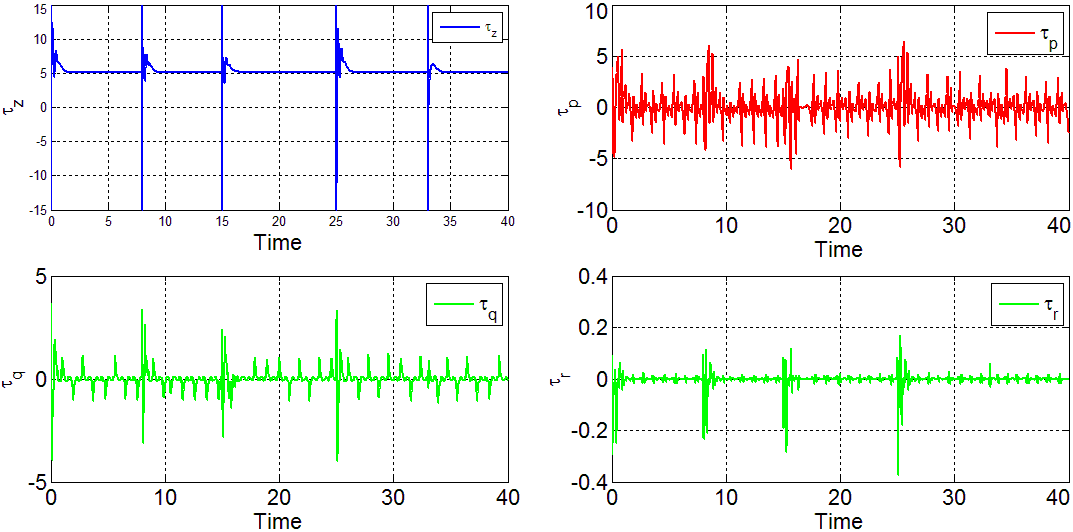}
     \caption{Control input of \Lone adaptive controller of a quadrotor system.}
     \label{fig:Chap_L1_Ex3_1cont}
    \end{figure}

\begin{figure}[!]
     \centering
     \includegraphics[height=5cm, width=13cm]{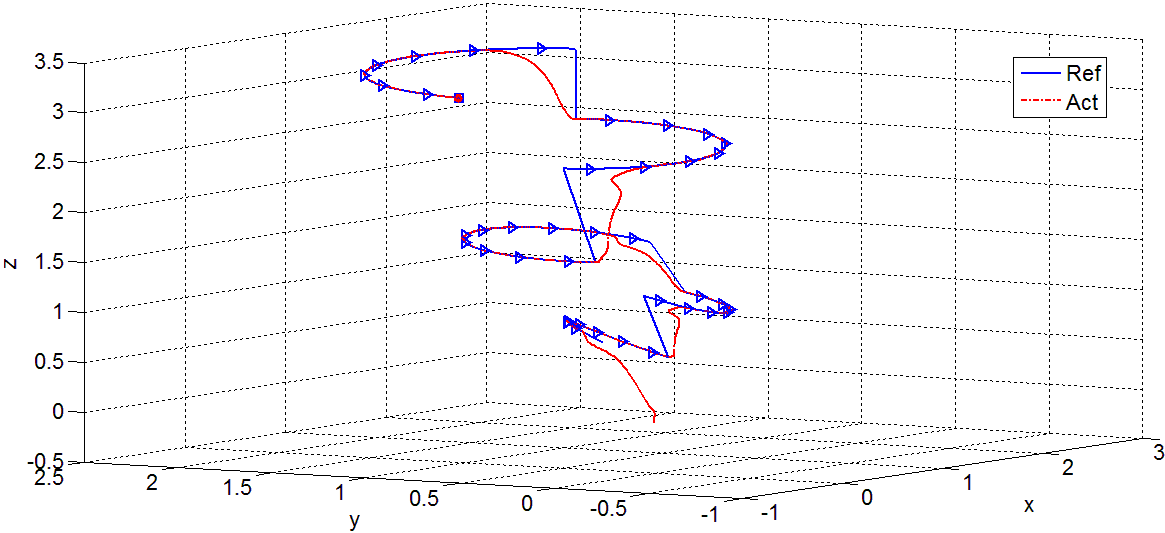}
     \caption{The 3D space tracking trajectory for both reference and actual output of a quadrotor system.}
     \label{fig:Chap_L1_Ex3_13D}
    \end{figure}

    {\bf Case 1:} We assume exact modeling and system with free disturbances where projection bounds of adaptation laws were defined numerically. Parameters of \Lone can be computed numerically where their bounds were chosen such as $\omega_l=0.5$, $\omega_u=10$, $\theta_b = 100$, $\sigma_b = 100$ and the adaptation gain $\Gamma = 100000$. The control input is constrained to $\tau_z = 15$ while other control signals are set free. The desired poles were set to $-30 \pm j0.5$, $-35 \pm j0.5$ and $-40 \pm j0.5$ and the feedback gain were set to diag(30,30,30). Figures (\ref{fig:Chap_L1_Ex3_1pos}), (\ref{fig:Chap_L1_Ex3_1ang}), (\ref{fig:Chap_L1_Ex3_1cont}) and (\ref{fig:Chap_L1_Ex3_13D}) represent the output positions, angles, control signals and 3D trajectory of quadrotor system by \Lone adaptive control respectively.
    
\begin{figure}[!]
     \centering
     \includegraphics[height=5cm, width=13cm]{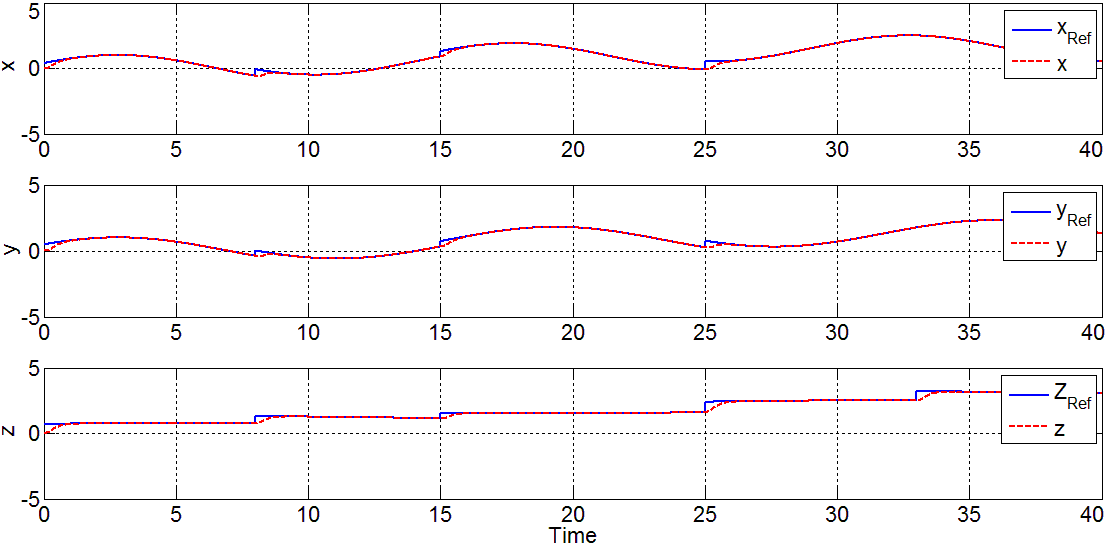}
     \caption{\Lone adaptive controller with reference and actual tracking positions for quadrotor.}
     \label{fig:Chap_L1_Ex3_2pos}
    \end{figure}

\begin{figure}[!]
     \centering
     \includegraphics[height=5cm, width=13cm]{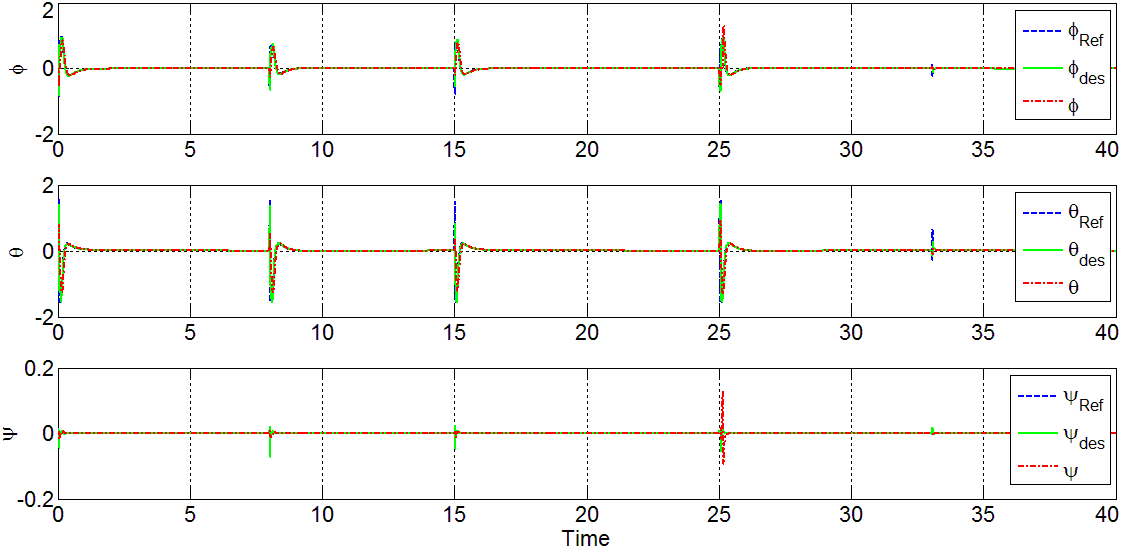}
     \caption{\Lone adaptive controller with reference, desired and actual tracking angles of a quadrotor syste.m}
     \label{fig:Chap_L1_Ex3_2ang}
    \end{figure}

\begin{figure}[!]
     \centering
     \includegraphics[height=5cm, width=13cm]{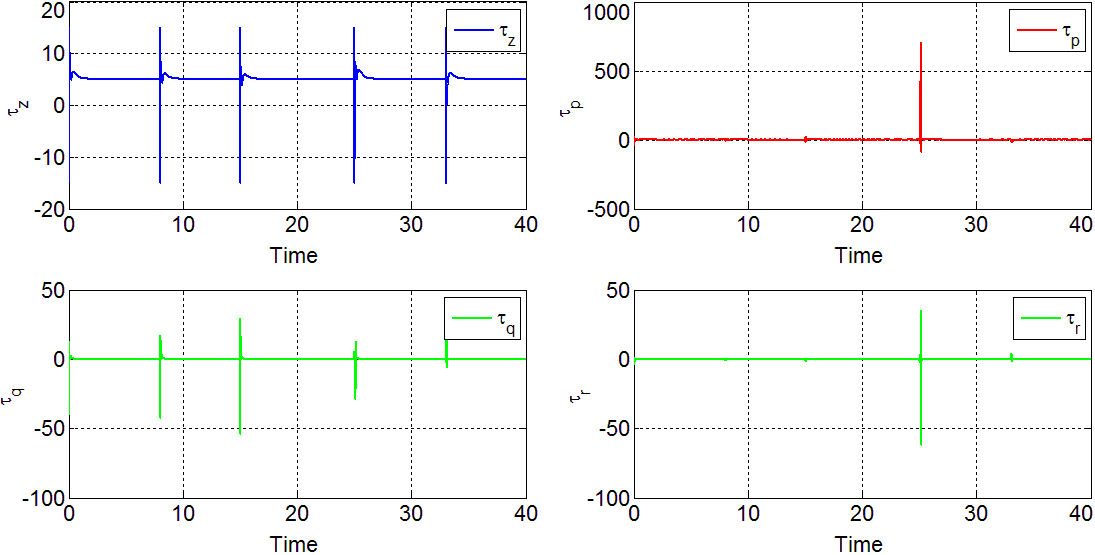}
     \caption{Control input of \Lone adaptive controller of a quadrotor system}
     \label{fig:Chap_L1_Ex3_2cont}
    \end{figure}

\begin{figure}[!]
     \centering
     \includegraphics[height=5cm, width=13cm]{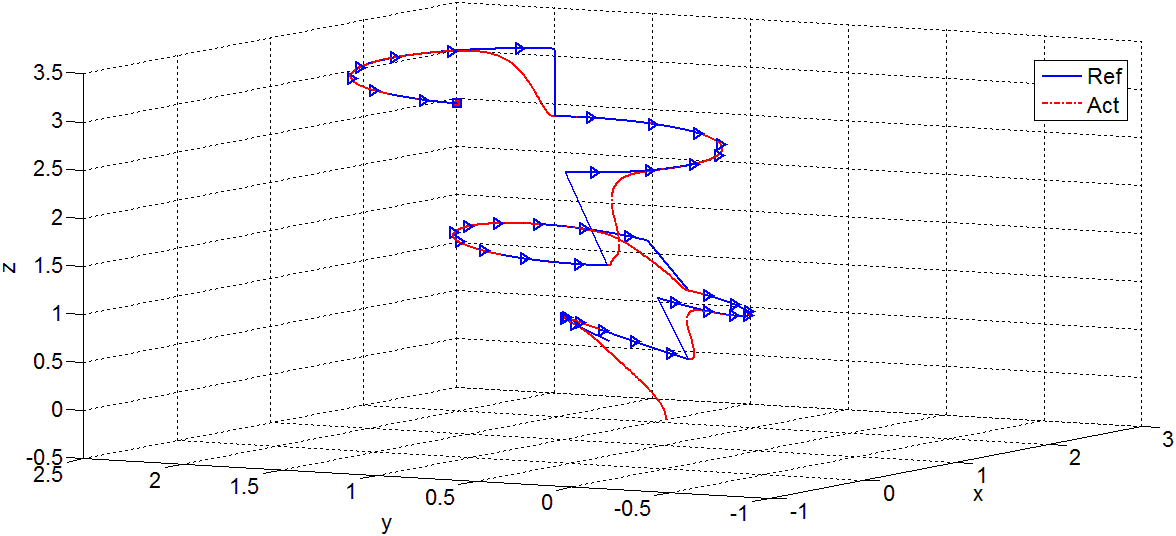}
     \caption{The 3D space tracking trajectory for both reference and actual output of a quadrotor system.}
     \label{fig:Chap_L1_Ex3_23D}
    \end{figure}

     {\bf Case 2:} Same assumptions and given data as mentioned in part 1 are considered here except the model is no longer exact. Uncertainties in the level of the states, disturbances and unmodeled input represented will be addressed into the system.
     \begin{equation*}
       \ddot{\eta}_1 = \frac{1}{m}R(\eta_2)\begin{bmatrix} 0 & 0 & \tau_z \end{bmatrix}^{\top} - g\begin{bmatrix} 0 & 0 & 1 \end{bmatrix}^{\top}
     \end{equation*}
     \begin{equation*}
       \ddot{\eta}_2 = f(\eta_2) + f_{\Delta}(\eta_2) + G_{\Delta}(\eta_2)G(\eta_2)\begin{bmatrix} \tau_p & \tau_q & \tau_r \end{bmatrix}^{\top} + D(s)
     \end{equation*}
     \begin{equation*}
      f_{\Delta}(\eta_2)=
      \begin{bmatrix} 
        0.2cos(\phi)sin(\theta)+0.2\phi\psi\\
        0.2cos(\phi)sin(\psi)+0.2\phi\psi^2\\
        0.2cos(\theta)sin(\phi)+0.2\phi\theta\psi
      \end{bmatrix},
      \hspace{10pt}
      D(s) =
      \begin{bmatrix} 
        \frac{0.2}{s+1}u_{d1}(s)\\
        \frac{0.24}{s^2+2s+3}u_{d2}(s)\\
        \frac{0.15}{s^2+3s+2}u_{d3}(s)
      \end{bmatrix}
     \end{equation*}
     \begin{equation*}
      G_{\Delta}(\eta_2) =
      \begin{bmatrix} 
        1.6 & 0 & 0\\
        0 & 0.7 & 0\\
        0 & 0 & 1.23
      \end{bmatrix}
     \end{equation*}
     \begin{equation*}
        u_{d1}\left(t\right) = sin(0.4t), \hspace{10pt} u_{d2}\left(t\right) = sin(0.6t), \hspace{10pt} u_{d3}\left(t\right) = sin(0.5t),
     \end{equation*}
    Figures (\ref{fig:Chap_L1_Ex3_2pos}), (\ref{fig:Chap_L1_Ex3_2ang}), (\ref{fig:Chap_L1_Ex3_2cont}) and (\ref{fig:Chap_L1_Ex3_23D}) are describing the output positions, angles, control signals and 3D trajectory of quadrotor system by \Lone adaptive control after admitting uncertainties, unmodeled input and disturbances.
    
    {\bf Example 3.3.3} Simulation Problem of Fully Actuated MARES Autonomous Underwater Vehicle\\
    MARES underwater vehicle model and parameters were defined in \cite{fossen_guidance_1994,ferreira_control_2009,ferreira_hydrodynamic_2009,ferreira_modeling_2010}. The submarine model can be represented as following 
     \begin{equation*}
       \tau_{\eta}(\eta) = M_{\eta}(\eta)\ddot{\eta} + C_{\eta}(\eta,\nu)\dot{\eta}+ D_{\eta}(\eta,\nu)\dot{\eta} + G_{\eta}(\eta)
     \end{equation*}
    Where $\eta$ is the earth coordinate frame, $G_{\eta}(\eta)$ is vector of gravitational/buoyancy forces and moments, $D_{\eta}(\eta,\nu)$ is damping matrix, $C_{\eta}(\eta,\nu)$ is coriolis-centripetal matrix (including added mass), $M_{\eta}(\eta)$ is system inertia matrix (including added mass) and $\tau_{\eta}(\eta)$ is the control input vector.
    
\begin{figure}[!]
      \centering
      \includegraphics[height=5cm, width=13cm]{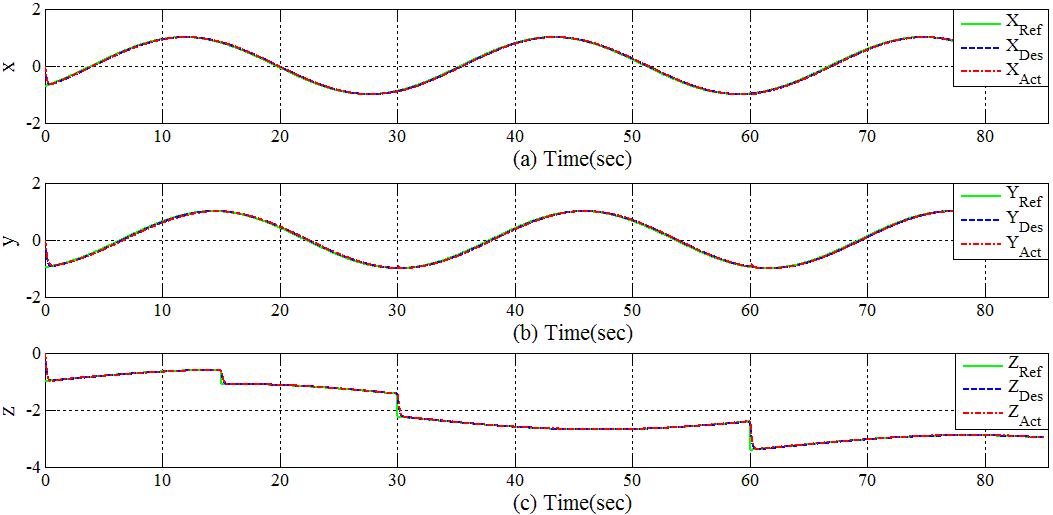}
      \caption{\Lone adaptive controller with reference and actual tracking positions of MARES.}
      \label{fig:Chap_L1_Ex4_1pos}
     \end{figure}

\begin{figure}[!]
      \centering
      \includegraphics[height=5cm, width=13cm]{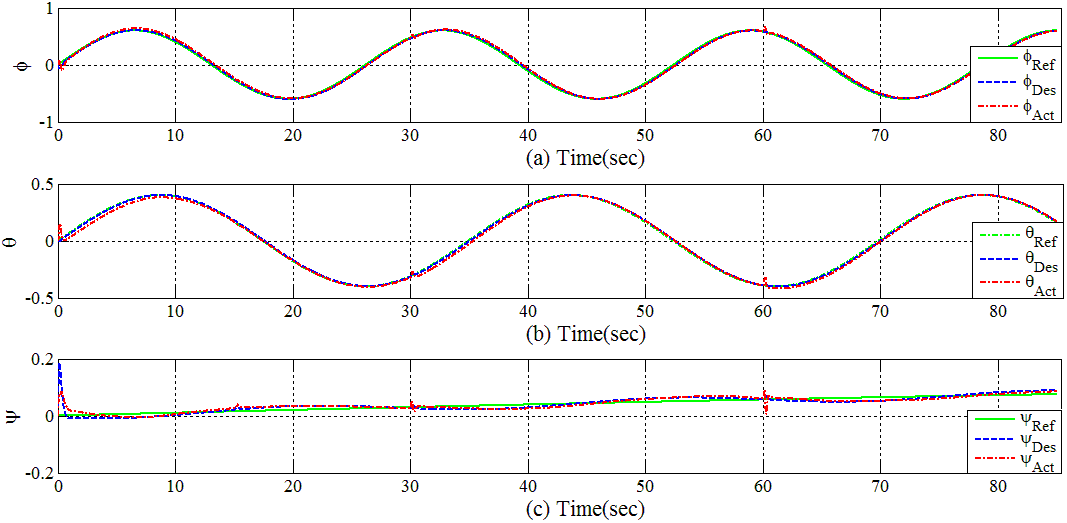}
      \caption{\Lone adaptive controller with reference, desired and actual tracking angles of MARES}
      \label{fig:Chap_L1_Ex4_1ang}
     \end{figure}

\begin{figure}[!]
      \centering
      \includegraphics[height=5cm, width=13cm]{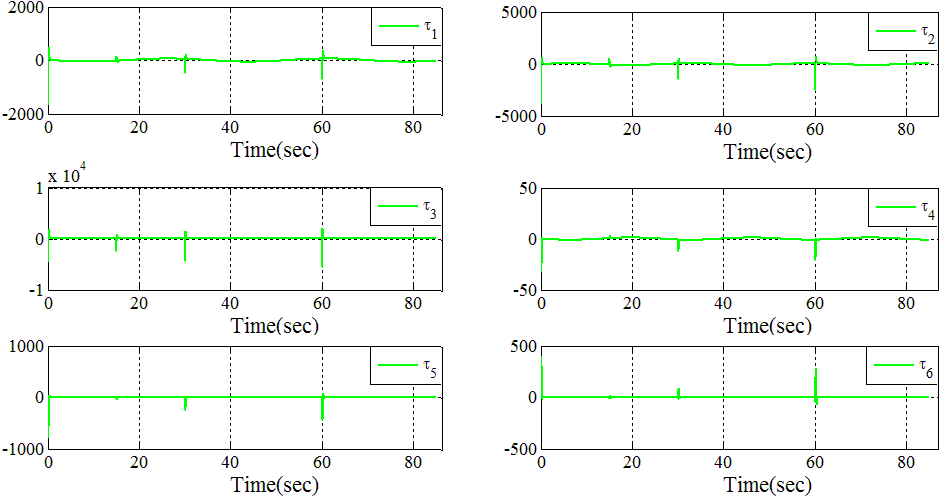}
      \caption{Control input of \Lone adaptive controller of MARES.}
      \label{fig:Chap_L1_Ex4_1cont}
     \end{figure}

\begin{figure}[!]
      \centering
      \includegraphics[height=5cm, width=13cm]{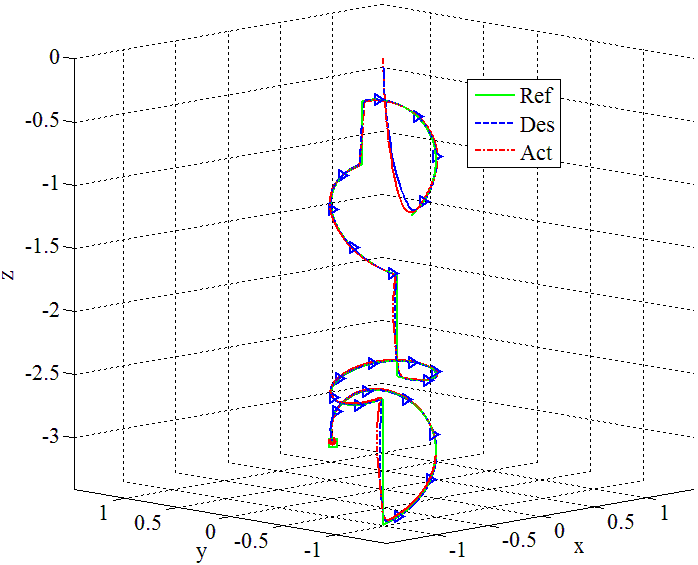}
      \caption{The 3D space tracking trajectory for both reference and actual output of MARES.}
      \label{fig:Chap_L1_Ex4_13D}
     \end{figure}

     Parameters of \Lone can be computed numerically where their bounds were chosen to $\omega_l=0.5$, $\omega_u=20$, $\theta_b = 100$, $\sigma_b = 100$ and the adaptation gain $\Gamma = 100000$. The desired poles are  $-9 \pm j0.1$, $-10.5 \pm j0.1$, $-12 \pm j0.1$, $-13.5 \pm j0.1$, $-15 \pm j0.1$ and $-16.5 \pm j0.1$. Finally, the feedback gain is diag(30,30,30,30,30,30).
    Figures (\ref{fig:Chap_L1_Ex4_1pos}), (\ref{fig:Chap_L1_Ex4_1ang}), (\ref{fig:Chap_L1_Ex4_1cont}) and (\ref{fig:Chap_L1_Ex4_13D}) are describing the output positions, angles, control signals and 3D trajectory respectively of MARES submarine using \Lone adaptive control. 
    
\section{\Lone Adaptive Controller for MIMO Systems in the Presence of Unmatched Nonlinear Uncertainties and Strong Coupling Effect}
    This section presents \Lone adaptive control for MIMO uncertain system in the presence of high nonlinearities with strong coupling effect. All uncertainties and nonlinearities are assumed  unknown.  
\subsection{Problem Formulation}
  Consider in the following class of systems:
    \begin{equation}
     \begin{aligned}
      \label{eq:ch2L1Actf3}
       & \dot{x}\left(t\right) = A_{m}x\left(t\right) + B_m\omega u\left(t\right) + f(x\left(t\right),z\left(t\right),t),\hspace{10pt}x(0)=x_0\\
       & \dot{x}_z = g(x\left(t\right),x_z\left(t\right),t),\hspace{10pt}x_z(0)=x_{z_0}\\
       & z\left(t\right) = g_{0}(x_z\left(t\right),t)\\
       & y\left(t\right) = Cx\left(t\right)
     \end{aligned}
    \end{equation}
   where $x\left(t\right) \in \mathbb{R}^n$ is the system state vector (measured); $u\left(t\right) \in \mathbb{R}^m$ is the control input; $y\left(t\right) \in \mathbb{R}^m$ is the system output;$B_m \in \mathbb{R}^{n \times m}$  is a known full column-rank constant matrix where $(A_m,B_m)$ is controllable and $C \in \mathbb{R}^{m \times n}$  is a known full-row rank constant matrix(known) where $(A_m,C_m)$ is observable; $A_m \in \mathbb{R}^{n \times n}$ is a known Hurwitz matrix that defines the desired dynamics for the closed-loop system; $\omega \in \mathbb{R}^{m \times m}$ is a gain matrix defines uncertain system input, $x_z\left(t\right)$ are the output and the state vector of internal unmodeled dynamics; also $f:\mathbb{R}\times\mathbb{R}^n\times\mathbb{R}^p\to\mathbb{R}^n$, $g_0:\mathbb{R}^l\times\mathbb{R}\to\mathbb{R}^p$ and $g:\mathbb{R}\times\mathbb{R}^l\times\mathbb{R}^n\to\mathbb{R}^l$ are unknown nonlinear continuous functions.
   
   The system in \eqref{eq:ch2L1Actf3} can also be written in the form
    \begin{equation}
     \begin{aligned}
      \label{eq:ch2L1Actf4}
       & \begin{split}
       \dot{x}\left(t\right) = & A_{m}x\left(t\right) + B_m(\omega u\left(t\right) + f_1(x\left(t\right),z\left(t\right),t))\\ 
       & + B_{um}(\omega u\left(t\right) + f_2(x\left(t\right),z\left(t\right),t)) ,\hspace{10pt}x(0)=x_0
       \end{split}\\
       & \dot{x}_z = g(x\left(t\right),x_z\left(t\right),t),\hspace{10pt}x_z(0)=x_{z_0}\\
       & z\left(t\right) = g_{0}(x_z\left(t\right),t)\\
       & y\left(t\right) = Cx\left(t\right)
     \end{aligned}
    \end{equation}
    Where $B_{um} \in \mathbb{R}^{n \times (n-m)}$ is a constant matrix has the property of $B_m \times B_{um}=0$ with $rank([B_m,B_{um}])=n$, while $f:\mathbb{R}\times\mathbb{R}^n\times\mathbb{R}^p\to\mathbb{R}^m$ and $f:\mathbb{R}\times\mathbb{R}^n\times\mathbb{R}^p\to\mathbb{R}^{n-m}$ are unknown nonlinear functions that verify
    \begin{equation}
     \begin{bmatrix}
      f_1(x\left(t\right),z\left(t\right),t)\\
      f_2(x\left(t\right),z\left(t\right),t)
     \end{bmatrix}
     =
     \begin{bmatrix}
        B_m & B_{um}
     \end{bmatrix}^{-1}
     f(x\left(t\right),z\left(t\right),t)
    \end{equation}
    where $f_1(\cdot)$ represents the matched component of the unknown nonlinearities and $f_2(\cdot)$ represents the unmatched uncertainties.
    Let $X \triangleq [x^{\top},z^{\top}]^{\top}$ and let $F(X\left(t\right),t) \triangleq f_i(x\left(t\right),z\left(t\right),t),i=1,2$. The system in \eqref{eq:ch2L1Actf4} verifies the following assumptions: 
     \begin{assumption} (The control input are partially known with known sign)
   The system input gain matrix $\omega$ assumed to be nonsingular and unknown with known diagonal sign with strictly row-diagonally dominant matrix.
    \begin{equation*}
     \omega \in \Omega \subset \mathbb{R}^{m \times m}
    \end{equation*}
    where $\Omega$ is assumed to be known convex compact set.
    \label{assum38}
     \end{assumption}
     \begin{assumption}  (Uniform boundedness of $f(0,t)$) 
      Let $B>0$ such that $ f_i(0,t) \leq B $ for all  $ t \geq 0 $.
      \label{assum39}
     \end{assumption}
     \begin{assumption} (Partial derivatives are semiglobal uniform bounded)
   For any $\delta >0$, there exist $d_{f_{x}} (\delta) >0 $ and $d_{f_{t}} (\delta) >0 $  such that for arbitrary $||x||_{\infty} \leq \delta$ and any $u$, the partial derivatives of $ f(x\left(t\right),t)) $ is piecewise-continuous and bounded,
    \begin{equation*}
     ||\frac{\partial f_i(x\left(t\right),t) }{\partial x}|| \leq d_{f_{xi}} (\delta),\hspace{10pt} ||\frac{\partial f_i(x\left(t\right),t) }{\partial t}|| \leq d_{f_{ti}} (\delta) \hspace{10pt} where \: i=1,2
    \end{equation*}
    \label{assum310}
     \end{assumption}
     \begin{assumption} (Asymptotically stability of initial conditions)
   The system assumed to start initially with $x_0$ inside an arbitrarily known set $\rho_0$ i.e., $||x_0||_{\infty} \leq \rho_0 < \infty$.\label{assum311}
     \end{assumption}
     \begin{assumption} (BIBO stability of internal dynamics) The states $x_z$ of internal dynamics are BIBO stable with respect to $x_{z0}$ and $x\left(t\right)$ and there exist $L_z,B_z>0$ such that for all $t \geq 0$
   \begin{equation*}
    ||z_t||_{\mathcal{L}_{\infty}} \leq L_z||x\left(t\right)||_{\mathcal{L}_{\infty}}+B_z
   \end{equation*} 
   \label{assum312}
     \end{assumption}
     \begin{assumption} (Stability of Transmission zeros)The transmission zeros of the transfer matrix $H_m(s) = C(sI-A_m)^{-1}B_m$ lie in the open left half complex plane.\label{assum313}
     \end{assumption}
   The objective in this section aims at designing a full-state feedback adaptive controller to ensure that $y\left(t\right)$ tracks a given bounded piecewise-continuous reference signal $r\left(t\right)$ with quantifiable performance bounds given $M(s)$.
   \begin{equation*}
    M(s) = C(sI-A_m)^{-1}B_mK_g(s)
   \end{equation*} 
   where $K_g(s)$ is a feedforward pre-filter,

\subsection{Definitions and \Lone-Norm Sufficient Condition for Stability}
   Let
   \begin{equation*}
     H_{xm}(s) \triangleq (sI-A_m)^{-1}B_m,
   \end{equation*}
   \begin{equation*}
     H_{xum}(s) \triangleq (sI-A_m)^{-1}B_{um},
   \end{equation*}
   \begin{equation*}
     H_{m}(s) \triangleq C(sI-A_m)^{-1}B_{m},
   \end{equation*}
   \begin{equation*}
     H_{um}(s) \triangleq C(sI-A_m)^{-1}B_{um},
   \end{equation*}
   and let $x_{in}\left(t\right)$ be the signal with Laplace transform $x_{in}(s) \triangleq (sI-A_m)^{-1} x_0$ and $\rho_{in} \triangleq ||s(sI-A_m)^{-1}||_{\mathcal{L}_{1}}\rho_{0}$. Since $A_m$ is Hurwitz and $x_0$ is bounded, then $||x_{in}||_{\mathcal{L_{\infty}}} \leq \rho_{in}$.\\
   \begin{equation}
    \label{eq:ch2Ldelta}
    L_{i_{\delta}}\triangleq \frac{\bar{\delta}(\delta)}{\delta} d_{f_{xi}}(\bar{\delta}(\delta)), \hspace{10pt} \bar{\delta}(\delta)\triangleq max \{\delta + \bar{\gamma},L_z(\delta+\bar{\gamma})B_z\}
   \end{equation}
   where $\bar{\gamma}$ is a small positive constant assigned arbitrarily.
   The objective of the adaptive controller aims in achieving DC gain $C(0)\triangleq \mathbb{I}_m$. $K$ is a feedback gain matrix and $D(s)$ is strictly proper transfer matrix and both of them aim to strictly proper transfer function as follows
   \begin{equation}
     C(s) \triangleq \omega K D(s)(\mathbb{I}_m +\omega K D(s))^{-1}
   \end{equation}
   The choice of $D(s)$ needs to ensure also that $C(s)H^{-1}(s)$ is a proper stable transfer matrix. For a particular class of systems, a possible choice for $D(s)$ might be $D(s)=1/s \cdot \mathbb{I}_m$, which yields a strictly proper $C(s)$ of the form
   \begin{equation}
     C(s) \triangleq \omega K(s\mathbb{I}_m +\omega K)^{-1}
   \end{equation}
  Now, for a given $\rho_0$, $k$ and $D(s)$ should be chosen such that there exist $\rho_{r} > \rho_{in}$ and the following  $\mathcal{L}_1$ norm condition verified
  \begin{equation}
  \label{eq:ch2rhor}
    ||G_{m}(s)||_{\mathcal{L}_1}+||G_{um}(s)||_{\mathcal{L}_1}\ell_0 < \frac{\rho_{r} - ||H_{xm}(s)C(s)K_g||_{\mathcal{L}_1}||r||_{\mathcal{L}_{\infty}} - \rho_{in}}{L_{1\rho_{r}}\rho_{r} + B_0}
  \end{equation}
  \begin{equation*}
    G_{m}(s) = H_{xm}(s)(\mathbb{I}_m - C(s))
  \end{equation*}
  \begin{equation*}
    G_{um}(s)= (\mathbb{I}_m - H_{xm}(s)C(s)H^{-1}_{m}(s))H_{xum}(s)
  \end{equation*}
  \begin{equation*}
    \ell_0 = \frac{L_{1\rho_{r}}}{L_{2\rho_{r}}}, \hspace{10pt} B_0 = max\{B_{10},B_{20}/\ell_0\}
  \end{equation*}
  let
  \begin{equation}
    \gamma_{1} \triangleq \frac{||H_{xm}(s)C(s)H^{-1}_{m}(s)||_{\mathcal{L}_1}}{1-||G_{m}(s))||_{\mathcal{L}_1}L_{1\rho_{r}}-||G_{um}(s))||_{\mathcal{L}_1}L_{2\rho_{r}}}\gamma_{0}+ \beta
  \end{equation}
  where $\gamma_{0}$ and $\beta$ are arbitrarily small positive constants.\\
  let
  \begin{equation}
    \rho_{u} \triangleq \rho_{ur} + \gamma_{2}
  \end{equation}
  where $\rho_{ur}$ and $\gamma_{2}$ are defined as following
  \begin{equation}
    \label{eq:ch2rhour}
    \begin{split}
    \rho_{ur} \triangleq &||\omega^{-1}C(s)||_{\mathcal{L}_1}\big(||K_g||_{\mathcal{L}_1}||r||_{\mathcal{L}_{\infty}} + L_{1\rho_{r}}\rho_{r} + B_{10}\\
    & + ||H^{-1}_{m}(s)H_{um}(s)||(L_{2\rho_{r}}\rho_{r} + B_{20})\big)
    \end{split}
  \end{equation}
  \begin{equation}
    \begin{split}
    \gamma_{2} \triangleq 
    &||\omega^{-1}C(s)||_{\mathcal{L}_1}L_{1\rho_{r}} \gamma_{1} + ||H^{-1}_{m}(s)H_{um}(s)||_{\mathcal{L}_1} L_{2\rho_{r}} \gamma_{1}\\
    &+||H^{-1}_{m}(s)C(s)||\gamma_{0}
    \end{split}
  \end{equation}
  and finally let
  \begin{equation}
    \label{eq:ch2L1Max3}
    \theta_{b_i} \triangleq L_{i_{\rho}},\hspace{10pt} \sigma_{b_i} \triangleq L_{i_{\rho}}B_z + B_i + \epsilon_i, \hspace{10pt} i = 1,2
  \end{equation}
  where $\epsilon$ is an arbitrary positive constant.
\subsection{\Lone Adaptive Control Architecture}
   {\bf State Predictor:} We consider the following state predictor:
   \begin{equation}
     \label{eq:ch2L1est3}
     \begin{aligned} 
        &\begin{split}
        \dot{\hat{x}}\left(t\right) = & A_{m}\hat{x}\left(t\right) + B_m(\hat{\omega} u\left(t\right) + \hat{\theta}_1 ||x\left(t\right)||_{\infty}+\hat{\sigma}_1 )\\
        &+ B_{um}(\hat{\theta}_2 ||x\left(t\right)||_{\infty}+\hat{\sigma}_2 ),\hspace{10pt}\hat{x}(0) = x(0)
     \end{split}\\
       & \hat{y}\left(t\right) = c\hat{x}\left(t\right)
     \end{aligned}
   \end{equation}
  where $\hat{\omega} \in \mathbb{R}^{m \times m}$, $\hat{\theta}_1\left(t\right)\in \mathbb{R}^{m}$, $\hat{\theta}_2\left(t\right)\in \mathbb{R}^{n-m}$, $\hat{\sigma}_1\left(t\right)\in\mathbb{R}^{m}$ and $\hat{\sigma}_2\left(t\right)\in\mathbb{R}^{n-m}$ are the adaptive estimates defined by the following adaptive laws.
   \begin{equation}
     \label{eq:ch2L1Proj3}
     \begin{aligned}
       &\dot{\hat{\omega}} = \Gamma Proj(\hat{\omega},-(\tilde{x}^{\top}PB_m)^{\top}(u\left(t\right)^{\top})), \hspace{10pt} \hat{\omega}(0) = \hat{\omega}_0  \\
       &\dot{\hat{\theta}}_1 = \Gamma Proj(\hat{\theta}_1,-(\tilde{x}^{\top}PB_m)^{\top}||x\left(t\right)||_{\infty}), \hspace{10pt} \hat{\theta}_1(0) = \hat{\theta}_{10}\\
       &\dot{\hat{\theta_2}} = \Gamma Proj(\hat{\theta}_2,-(\tilde{x}^{\top}PB_{um})^{\top}||x\left(t\right)||_{\infty}), \hspace{10pt} \hat{\theta}_2(0) = \hat{\theta}_{20}\\
       &\dot{\hat{\sigma}}_1 = \Gamma Proj(\hat{\sigma}_1,-(\tilde{x}^{\top}PB_m)^{\top}), \hspace{10pt} \hat{\sigma}_1(0) = \hat{\sigma}_{10}\\
       &\dot{\hat{\sigma}}_2 = \Gamma Proj(\hat{\sigma}_1,-(\tilde{x}^{\top}PB_m)^{\top}), \hspace{10pt} \hat{\sigma}_2(0) = \hat{\sigma}_{20}
     \end{aligned}
   \end{equation}
    where $\tilde{x} \triangleq \hat{x} - x\left(t\right)$, $\Gamma \in \mathbb{R}^{+}$ is the adaptation gain, and $P=P^{\top}>0$ is defined by solving the algebraic Lyapunov equation $A_m^{\top}P+PA_m=-Q$ for arbitrary symmetric $Q=Q^{\top}>0$. The projection operator ensures that $\hat{\omega} \in \Omega$, $||\hat{\theta}_i||_{\infty} \in \Theta_i$, $||\hat{\sigma}_i|| \leq \Delta_i$,  with  $\theta_{bi}$ and $\delta_{bi}$ are being defined by \eqref{eq:ch2L1Max3}
    
     {\bf Control Law:} Control signal can be calculated as following
     \begin{equation}
       u(s) = -kD(s)\hat{\eta}(s)
     \end{equation}
     where $r(s)$ and $\hat{\eta}(s)$ are the Laplace transforms of $r\left(t\right)$  and $\hat{\eta}\left(t\right) = \hat{\omega} u\left(t\right) + \hat{\eta}_1 + \hat{\eta}_2 - K_gr\left(t\right)$ respectively; and the necessary feedforward gain in order to get unity steady state is calculated by $K_g \triangleq -(CA_{m}^{-1}B)^{-1}$ ; the feedback gain $k$ is positive constant and $D(s)$ is a strictly proper transfer function where both of them lead to a strictly proper stable closed loop system.
     \begin{equation*}
       \hat{\eta}_1 \triangleq \hat{\theta}_1 ||x\left(t\right)||_{\infty}+\hat{\sigma}_1
     \end{equation*}
     \begin{equation*}
       \hat{\eta}_2 \triangleq \hat{\theta}_2 ||x\left(t\right)||_{\infty}+\hat{\sigma}_2
     \end{equation*}
 
\subsection{\Lone Adaptive Control Stability Analysis}
     {\bf Transient and Steady-State Performance:} The error between system dynamics in \eqref{eq:ch2L1Actf4} and state predictor in \eqref{eq:ch2L1est3} can be written as
     \begin{equation}
       \label{eq:ch2L1err3}
       \dot{\tilde{x}}\left(t\right) = A_{m}\tilde{x}\left(t\right) + B_m(\tilde{\omega} u\left(t\right) + \tilde{\eta}_1\left(t\right)) + B_{um}\tilde{\eta}_2\left(t\right)\\
     \end{equation}
     Where $\tilde{x} = \hat{x} - x$, $\tilde{\theta}_i = \hat{\theta}_i - \theta_i$, $\tilde{\omega} = \hat{\omega} - \omega$, $\tilde{\sigma}_i = \hat{\sigma}_i - \sigma_i$ and $\tilde{\eta}_i = \hat{\eta}_i - \eta_i$ where $i=1,2$.
     
     \begin{lemma}
      The prediction error $\tilde{x}\left(t\right)$ is uniformly bounded,
      \label{lemma33}
     \end{lemma}
     from Lemma \ref{lemma33} and equations (\ref{eq:ch2xrho2}) and (\ref{eq:ch2urho2}), the derivatives of $\omega$, $\theta$ and $\sigma$ are bounded:
     \begin{equation}
     \label{eq:ch2thetadot3}
         ||\theta_i||_{\infty} \leq \theta_{b_i}(\rho_r) < \infty, \hspace{10pt}  ||\dot{\theta}_i||_{\infty} \leq d_{\theta_i}(\rho_r) < \infty
     \end{equation}
     \begin{equation}
     \label{eq:ch2sigmadot3}
         ||\sigma_i||_{\infty} \leq \sigma_{b_i}(\rho_r) < \infty, \hspace{10pt}  ||\dot{\sigma}_i||_{\infty} \leq d_{\sigma_{i}}(\rho_r) < \infty
     \end{equation}
     Then
     \begin{equation}
       ||\tilde{x}||_{\infty} \leq \sqrt{\frac{\theta_m}{\lambda_{min}(P)\Gamma}}
     \end{equation} 
     where
     \begin{equation}
     \label{eq:ch2thetam3}
     \begin{split}
       \theta_m \triangleq & 4\bigg((\theta_{b_1}^2 +\sigma_{b_1}^2)m +(\theta_{b_2}^2 +\sigma_{b_2}^2)(n-m) + \max_{\omega\in\Omega}tr(\omega^{\top}\omega) + \\ 
       & 4\frac{\lambda_{max}(P)}{\lambda_{min}(Q)}\big((d_{\theta_1}\theta_{b_1} + d_{\sigma_1}\sigma_{b_1})m +(d_{\theta_2}\theta_{b_2} + d_{\sigma_2}\sigma_{b_2})(n-m)\big) \bigg)
     \end{split}
     \end{equation}
     which will be verified as following\\
     {\bf Stability proof:} Consider the Lyapunov function candidate
     \begin{equation}
       \label{eq:ch2VPart3}
       V(\tilde{x},\tilde{\omega},\tilde{\theta}_i,\tilde{\sigma}_i) = \tilde{x}^{\top}P\tilde{x} + \frac{1}{\Gamma}\big(tr(\tilde{\omega}^{\top}\tilde{\omega}) + \tilde{\theta}^{\top}\tilde{\theta} + \tilde{\sigma}^{\top}\tilde{\sigma}\big) 
     \end{equation}
     Since $\hat{x}(0) = x(0)$ then we can verify that
     \begin{equation*}
       V(0) \leq  \frac{4}{\Gamma}\big(\max_{\omega\in\Omega}tr(\omega^{\top}\omega) + \theta_{b_1}^2 +\sigma_{b_1}^2)m +(\theta_{b_2}^2 +\sigma_{b_2}^2)(n-m)\big) \leq \frac{\theta_m}{\Gamma}
     \end{equation*}
     \begin{equation*}
       \begin{split}
       \dot{V} \leq & \tilde{x}^{\top}Q\tilde{x} + \frac{2}{\Gamma}(\dot{\hat{\theta}} + \tilde{x}^{\top}PB||x||_{\infty}) +\frac{2}{\Gamma}(\dot{\hat{\sigma}} + \tilde{x}^{\top}PB) +\frac{2}{\Gamma}(\dot{\hat{\omega}} + \tilde{x}^{\top}PBu)\\
       &  - \frac{2}{\Gamma}\sum\limits_{i=1}^2(\tilde{\theta}_i^{\top}\dot{\theta}_i + \tilde{\sigma}_i^{\top}\dot{\sigma}_i)
       \end{split}
     \end{equation*}
     \begin{equation}
       \dot{V} = -\tilde{x}^{\top}Q\tilde{x} + \frac{2}{\Gamma}\sum\limits_{i=1}^2\big(|\tilde{\theta}_i^{\top}\dot{\theta}_i| + |\tilde{\sigma}_i^{\top}\dot{\sigma}_i|)
     \end{equation}
     \begin{equation}
     \label{eq:ch2_vodt3_1}
       \dot{V} \leq -\tilde{x}^{\top}Q\tilde{x} + \frac{4}{\Gamma}\big((d_{\theta_1}\theta_{b_1} + d_{\sigma_1}\sigma_{b_1})m + (d_{\theta_2}\theta_{b_2} + d_{\sigma_2}\sigma_{b_2})(n-m)\big)
     \end{equation}
     Now we can say
     \begin{equation}
      \sum\limits_{i=1}^2\big(\tilde{\theta}_i^{\top}\dot{\theta}_i + \tilde{\sigma}_i^{\top}\dot{\sigma}_i) \leq \sum\limits_{i=1}^2(d_{\theta_i}\theta_{b_i} + d_{\sigma_i}\sigma_{b_i})
     \end{equation}
     Moreover, the projection operator also ensures that
     \begin{equation}
      \label{eq:ch2_v3_1}
       \max_{t \geq 0} \big(\frac{1}{\Gamma}(\tilde{\theta}^{\top}\tilde{\theta}+tr(\tilde{\omega}^{\top}\tilde{\omega}) +\tilde{\sigma}^{\top}\tilde{\sigma}\big) \leq  \frac{1}{\Gamma}(\theta_b^2m + \Delta^{2}m + \max_{\omega\in\Omega}tr(\omega^{\top}\omega))
     \end{equation}
     which holds for all $t \geq 0$.
     If at any time $t_1>0$, one has $V(t_1) \geq \theta_m/\Gamma$, then it follows from \eqref{eq:ch2thetam3} and \eqref{eq:ch2VPart3} that
     \begin{equation}
      \label{eq:ch2_v3_2}
       \tilde{x}^{\top}(t_1)P\tilde{x}(t_1) > \frac{4}{\Gamma}\frac{\lambda_{max}(P)}{\lambda_{min}(Q)}\big((d_{\theta_1}\theta_{b_1} + d_{\sigma_1}\sigma_{b_1})m +(d_{\theta_2}\theta_{b_2} + d_{\sigma_2}\sigma_{b_2})(n-m)\big)
     \end{equation}
     thus
     \begin{equation}
     \label{eq:ch2_vodt3_2}
     \begin{split}
      & \tilde{x}(t_1)^{\top}Q\tilde{x}(t_1) \geq \frac{\lambda_{min}(Q)}{\lambda_{max}(P)}\tilde{x}^{\top}(t_1)P\tilde{x}(t_1) \\ 
      &\hspace{40pt} >\frac{4}{\Gamma}\sum\limits_{i=1}^2\big((d_{\theta_1}\theta_{b_1} + d_{\sigma_1}\sigma_{b_1})m + (d_{\theta_2}\theta_{b_2} + d_{\sigma_2}\sigma_{b_2})(n-m)\big)
     \end{split}
     \end{equation}
     Hence, if $V(t_1) \geq \theta_m/\Gamma$, then from \eqref{eq:ch2_vodt3_1} and \eqref{eq:ch2_vodt3_2} we have
     \begin{equation}
       \dot{V} \leq 0
     \end{equation}
\subsection{Problem Formulation and Simulation}
\begin{figure}[!]
      \centering
      \includegraphics[height=5cm, width=13cm]{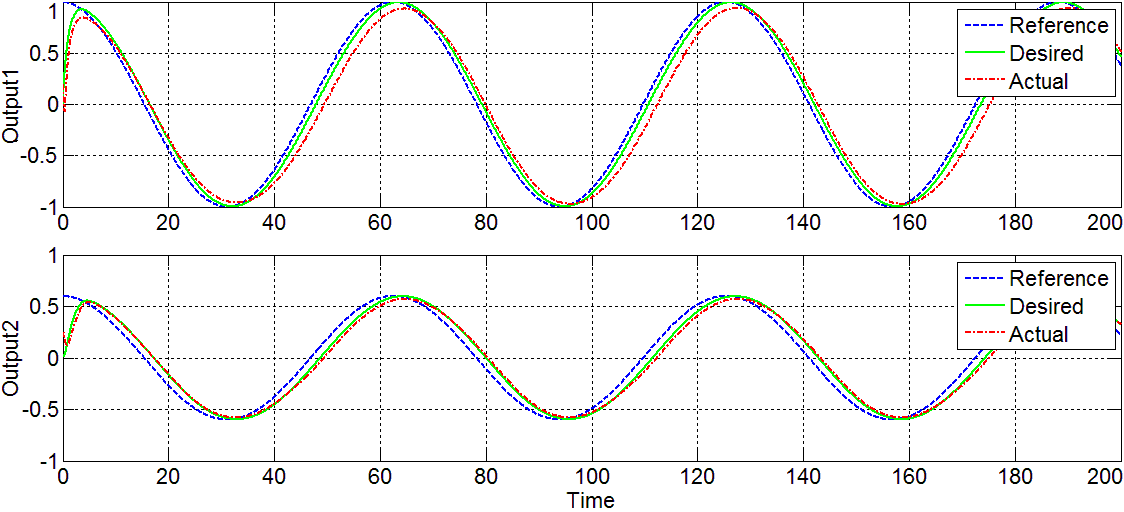}
      \caption{Tracking output of \Lone adaptive control with reference and desired outputs for unmatched MIMO uncertain system.}
      \label{fig:Chap_L1_Ex5_1}
     \end{figure}

\begin{figure}[!]
      \centering
      \includegraphics[height=5cm, width=13cm]{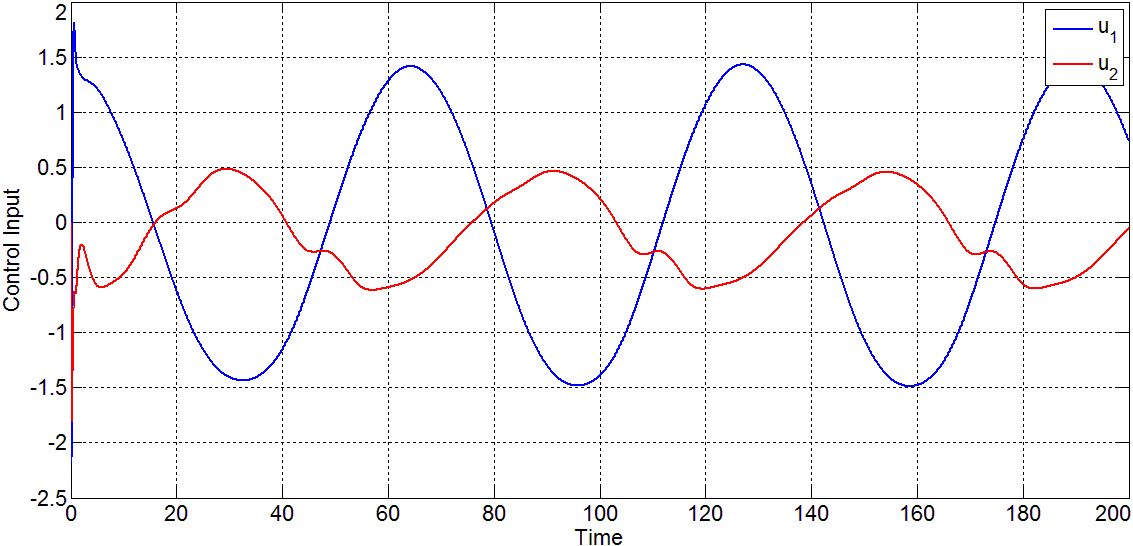}
      \caption{Control signal of \Lone adaptive control for unmatched MIMO uncertain system.}
      \label{fig:Chap_L1_Ex5_2}
     \end{figure}

    {\bf Example 3.4.1}
    MIMO System with Nonlinear Unmatched Uncertainties.\\
    \Lone adaptive control will be implemented to high nonlinear system with unmatched uncertainties in order to investigate output performance and control signals. Consider the system in \cite{xargay_l1_2010}.
     \begin{equation*}
     \begin{aligned}
       & \dot{x}\left(t\right) = (A_{m} + A_{\Delta})x\left(t\right) + B_m\omega u\left(t\right) + f_{\Delta}(x\left(t\right),z\left(t\right),t)\\ 
       & y\left(t\right) = Cx\left(t\right)
     \end{aligned}
     \end{equation*}
    where
     \begin{equation*}
       A_{m} =
       \begin{bmatrix}
       -1 & 0 & 0\\
       0 & 0 & 1\\
       0 & -1 & -1.8
       \end{bmatrix}, \hspace{10pt}
     B_{m} =
       \begin{bmatrix}
        1 & 0 \\
        0 & 0 \\
        1 & 1
       \end{bmatrix}, \hspace{10pt}
       C =
       \begin{bmatrix}
       1 & 0 & 0\\
       0 & 1 & 0
       \end{bmatrix}
     \end{equation*}
    while $A_{\Delta} \in \mathbb{R}^{3 \times 3}$ and $\omega_{\Delta} \in \mathbb{R}^{2 \times 2}$ are unknown constant matrices satisfying
     \begin{equation*}
      \omega \in 
      \begin{bmatrix}
       [0.6,1.2] & [-0.2,0.2]\\
       [-0.2,0.2] & [0.6,1.2]
      \end{bmatrix} = \Omega
     \end{equation*}
    and $f_{\Delta}$ is the (unknown) nonlinear function
     \begin{equation*}
      f_{\Delta}(x\left(t\right),z\left(t\right),t) = 
      \begin{bmatrix}
       \frac{k_1}{3}x^{\top}x + \tanh(\frac{k_2}{2}x_1)x_1 + k_3z\\
       \frac{k_4}{2}sec(x_2)x_2 + \frac{k_5}{5}x_3^2 + k_6(1-e^{-\lambda t}) + \frac{k_7}{2}z\\
       k_8x_3cos(\omega_ut) + k_9z^2\\
      \end{bmatrix}
     \end{equation*}
     where $k_1 = -1$, $k_2 = 1$, $k_3 = 0$, $k_4 = 1$, $k_5 = 0$, $k_6 = 0.2$, $k_7 = 1$, $k_8 = 0.6$, $k_9 = -0.7$, $\lambda = 0.3$ and $\omega_u = 5$.
     The internal unmodeled dynamics are given by
     \begin{equation*}
        \begin{aligned}
        & \dot{x}_{z1} = x_{z2}\left(t\right)\\
        & \dot{x}_{z2} = -x_{z1}\left(t\right) + 0.8(1-x_{z1}^2\left(t\right))x_{z2}\left(t\right)\\
        & z\left(t\right) = 0.1(x_{z1}\left(t\right) - x_{z2}\left(t\right)) + z_u\left(t\right)\\
        & z(s) = \frac{-s+1}{100s^2 + 8s + 1}\begin{bmatrix}1 & -2 & 1 \end{bmatrix}x(s)
        \end{aligned}
     \end{equation*}
     Desired poles are chosen as $p=-1,-0.9 \pm j0.4359$, $\Gamma = 80000$ and 
     \begin{equation*}
     Q = 
      \begin{bmatrix}
       1  &  0  &  0\\
       0  &  1  &  0\\
       0  &  0  &  1\\
      \end{bmatrix}, \hspace{10pt}
     K = 
      \begin{bmatrix}
       8  &  0\\
       0  &  8\\
      \end{bmatrix}
     \end{equation*}
     \begin{equation*}
        D(s) = \frac{1}{s(s/25+1)(s/70+1)(s^2/40^2 + 1.8s/40 + 1)}\mathbb{I}_2
     \end{equation*}
     Adaptive estimates belong to the following bounds $\hat{\theta}_1\left(t\right) \in [-40,40]\mathbb{I}_2$, $\hat{\theta}_2\left(t\right) \in [-40,40]$, $\hat{\sigma}_1\left(t\right) \in [-5,5]\mathbb{I}_2$, $\hat{\sigma}_2\left(t\right) \in [-5,5]$, $\hat{\omega}_{11}\left(t\right),\hat{\omega}_{22}\left(t\right) \in [0.25,3]$, and $\hat{\omega}_{12}\left(t\right),\hat{\omega}_{21}\left(t\right) \in [-0.2,0.2]$.
     Also other uncertainities and modeled input parameters will be defined by
     \begin{equation*}
        A_{\Delta} =
      \begin{bmatrix}
       0.2  &  -0.2  &  -0.3\\
       -0.2  &  -0.2 &  0.6\\
       -0.1  &  0  &  -0.9
      \end{bmatrix}, \hspace{10pt} \omega = 
      \begin{bmatrix}
       0.6  &  -0.2\\
       0.2  &  1.2\\
      \end{bmatrix}
     \end{equation*}
     Figure ~\ref{fig:Chap_L1_Ex5_1} and ~\ref{fig:Chap_L1_Ex5_2} show output response and control signals of \Lone adaptive control.

\begin{figure}[!]
      \centering
      \includegraphics[height=5cm, width=13cm]{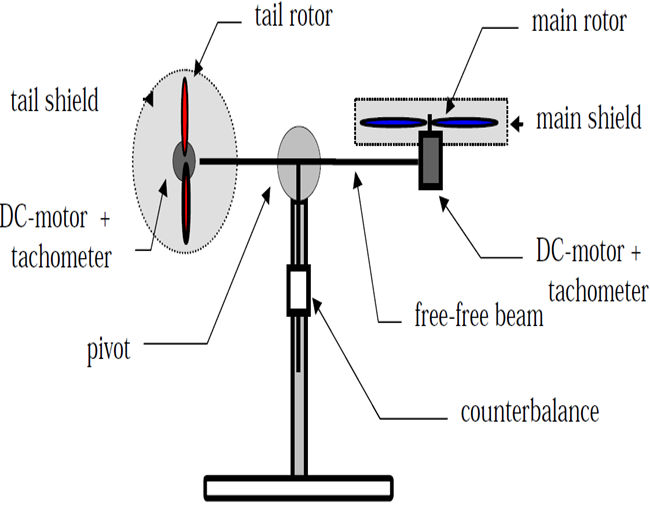}
      \caption{Laboratory set-up of TRMS.}
      \label{fig:Fig281}
     \end{figure}

    {\bf Example 3.4.2}
     Nonlinear Twin Rotor MIMO System (TRMS) with Strong Coupling.\\
     Twin rotor was designed for training high nonlinear control applications to mimic the behavior of the helicopter dynamics in terms of angle orientation \cite{_twin_1998}. The model and parameters of the system are defined in \cite{pratap_sliding_2010}. Complexity of the twin rotor comes from high nonlinearities in addition to strong coupling between control signals. Figure ~\ref{fig:Fig281} demonstrates TRMS set up. \Lone adaptive control will be implemented on high nonlinear TRMS with strong coupling effect in order to evaluate the control performance on output response and control signals.
     Adaptive estimates were defined as $\hat{\theta}_1\left(t\right) \in [-50,50]\mathbb{I}_2$, $\hat{\theta}_2\left(t\right) \in [-50,50]$, $\hat{\sigma}_1\left(t\right) \in [-15,15]\mathbb{I}_2$, $\hat{\sigma}_2\left(t\right) \in [-15,15]$, $\hat{\omega}_{11}\left(t\right),\hat{\omega}_{22}\left(t\right) \in [0.25,5]$, $\Gamma = 100000$ and the desired poles are assigned to $-15 \pm 0.3i$, $-17 \pm 0.5i$ and $-20 \pm 0.5i$ and finally the feedback gain = $5\big (\begin{smallmatrix} 1 & 0\\ 0 & 1 \end{smallmatrix}\big)$.

\begin{figure}[!]
      \centering
      \includegraphics[height=5cm, width=13cm]{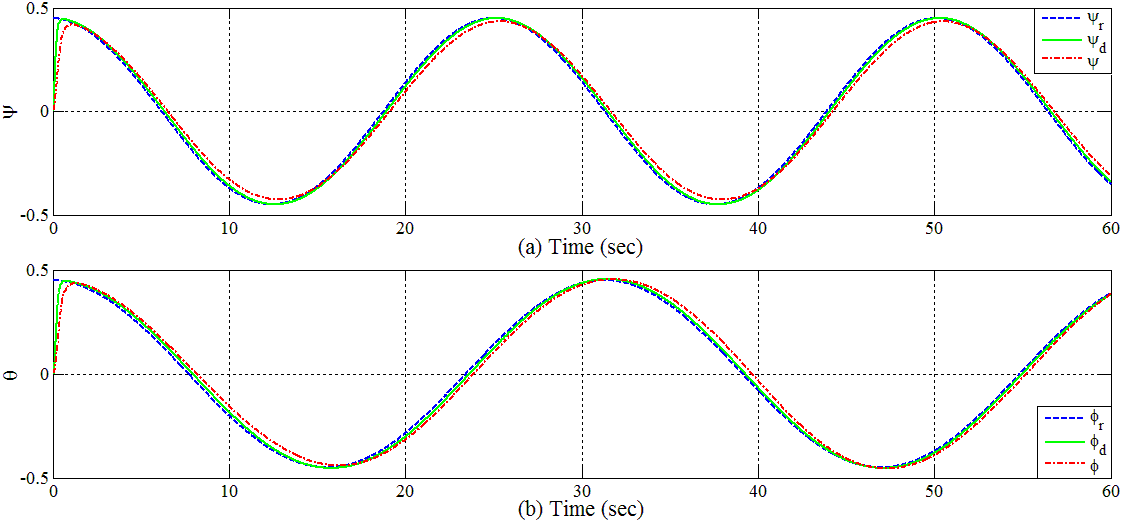}
      \caption{Tracking output of \Lone adaptive control with reference and desired outputs for TRMS.}
      \label{fig:Chap_L1_Ex6_1}
     \end{figure}

\begin{figure}[!]
      \centering
      \includegraphics[height=5cm, width=13cm]{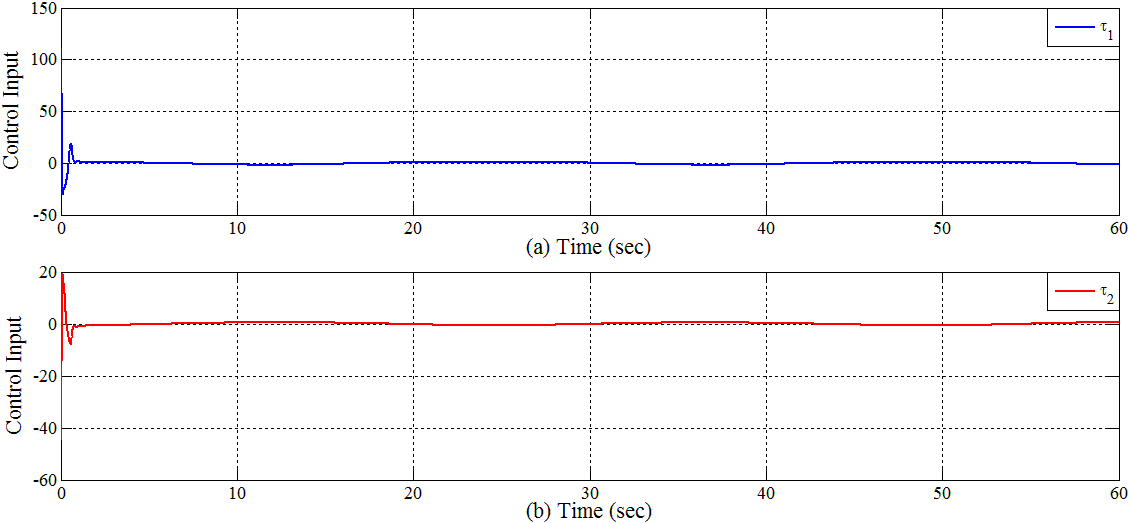}
      \caption{Control signal of \Lone adaptive control for TRMS.}
      \label{fig:Chap_L1_Ex6_2}
     \end{figure}

     Figure ~\ref{fig:Chap_L1_Ex6_1} and ~\ref{fig:Chap_L1_Ex6_2} show output response and control signals of \Lone adaptive control for TRMS.
     
\section{Conclusion}
  This chapter mainly handled \Lone adaptive controller from different perspectives and for different classes of nonlinear systems. The robustness, transient performance and tracking trajectory are prominent features of \Lone adaptive controller. All previous features have been validated through different cases of studies including reproducing recent results. From the literature, the relation between improving robustness, enhancing transient performance and control signal range have been demonstrated. In conclusion, improving robustness and enhancing the transient performance have a direct effect on the control signal range. We will present a satisfactory solution will be studied in subsequent chapters.

\clearpage

\newpage


\chapter{A Fuzzy Logic Feedback Filter Design Tuned with PSO for $ \mathcal{L}_1 $ Adaptive Controller}

\section{Introduction}   
\begin{figure*}[ht]
                    \centering
                    \includegraphics[scale=0.8]{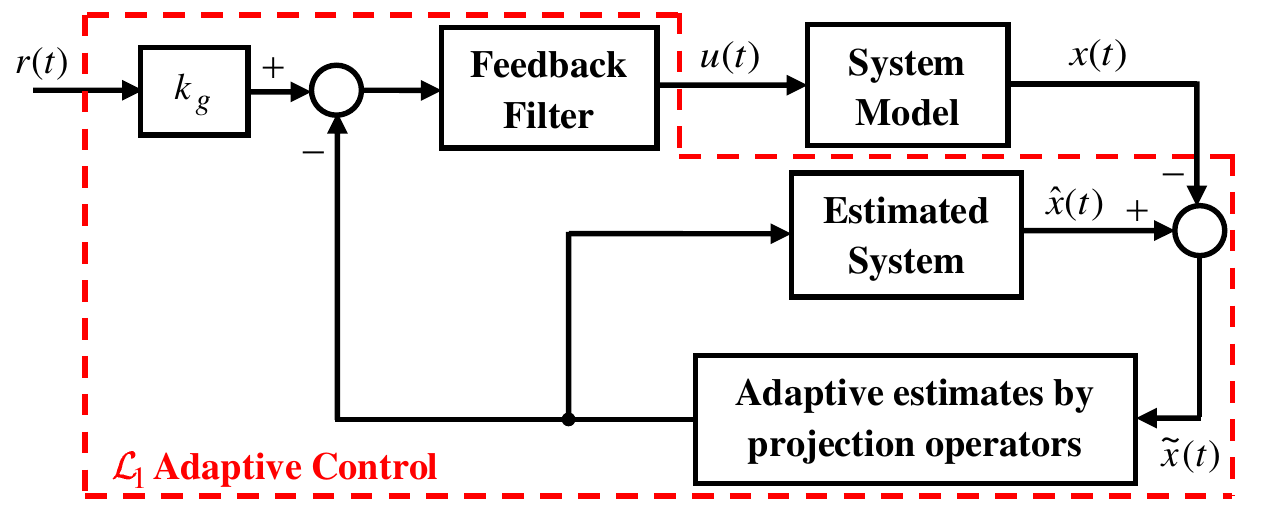}
                    \caption{The general structure of L1 adaptive controller.}\label{Fuzzy_L1_Gen}
          \end{figure*}

     The structure of \Lone adaptive controller offers three features including the implementation of a low pass filter in order to limit the frequency range of the control signal and reduce the effect of the uncertainties (see Figure (\ref{Fuzzy_L1_Gen})). 
        The structure allows decoupling of the adaption and robustness using high-gain for fast adaption.
        The filter is selected such that the system's output tracks properly the reference input and the undesirable uncertainties and frequencies are filtered ( see \cite{cao_design_2006} or \cite{hovakimyan_l1_2010}). Using the low pass filter, \Lone controller reduces the coupling between robustness and fast adaptation and provides infinity norm boundedness of the transient and steady state responses.
     \Lone adaptive control was first introduced by \cite{cao_design_2006}.  It has been applied successfully to uncertain linear systems \cite{cao_design_2008}, uncertain nonlinear single-input-single-output (SISO) systems \cite{cao_guaranteed_2007}, \cite{luo_l_2010}, and nonlinear system multi-input-multi-output (MIMO) with unmatched uncertainties \cite{xargay_l1_2010-1}. And, the control approach showed satisfactory results on experimental flight tests \cite{gregory_l1_2009},  \cite{xargay_l1_2010}. 
     The optimal structure of \Lone filter has been studied extensively in \cite{hovakimyan_l1_2010}. The trade-off between fast desired closed loop dynamics and filter parameters has been debated for long \cite{cao_design_2006,hovakimyan_l1_2010,li_filter_2008,li_optimization_2007,kharisov_limiting_2011,kim_multi-criteria_2014}.  Increasing the bandwidth of the low pass filter will reduce robustness margin, which will require slowing the desired closed loop performance in order to regain the robustness. However, slower selection of desired closed loop performance will deteriorate the output performance especially during the transient period \cite{hovakimyan_l1_2010}.  Limitations of \Lone adaptive controller and the interconnection between adaptive estimates and the feedback filter were studied in \cite{kharisov_limiting_2011}, where several filter designs were considered based on the use of disturbance observer. The authors showed that it is crucial to select the appropriate coefficients for a given filter to achieve the desired performance. Several attempts on identifying these optimal coefficients have been made in the literature.  This includes convex optimization based on linear matrix inequality (LMI) \cite{hovakimyan_l1_2010}, \cite{li_filter_2008} and multi-objective optimization using MATLAB optimization solver \cite{li_optimization_2007}. More recently, a systematic approach was presented in \cite{kim_multi-criteria_2014} to determine the optimal feedback filter coefficients in order to increase the zone of robustness margin. The authors proposed the use of greedy randomized algorithms. \\
       
\begin{figure*}[ht]
               \centering
               \includegraphics[scale=0.8]{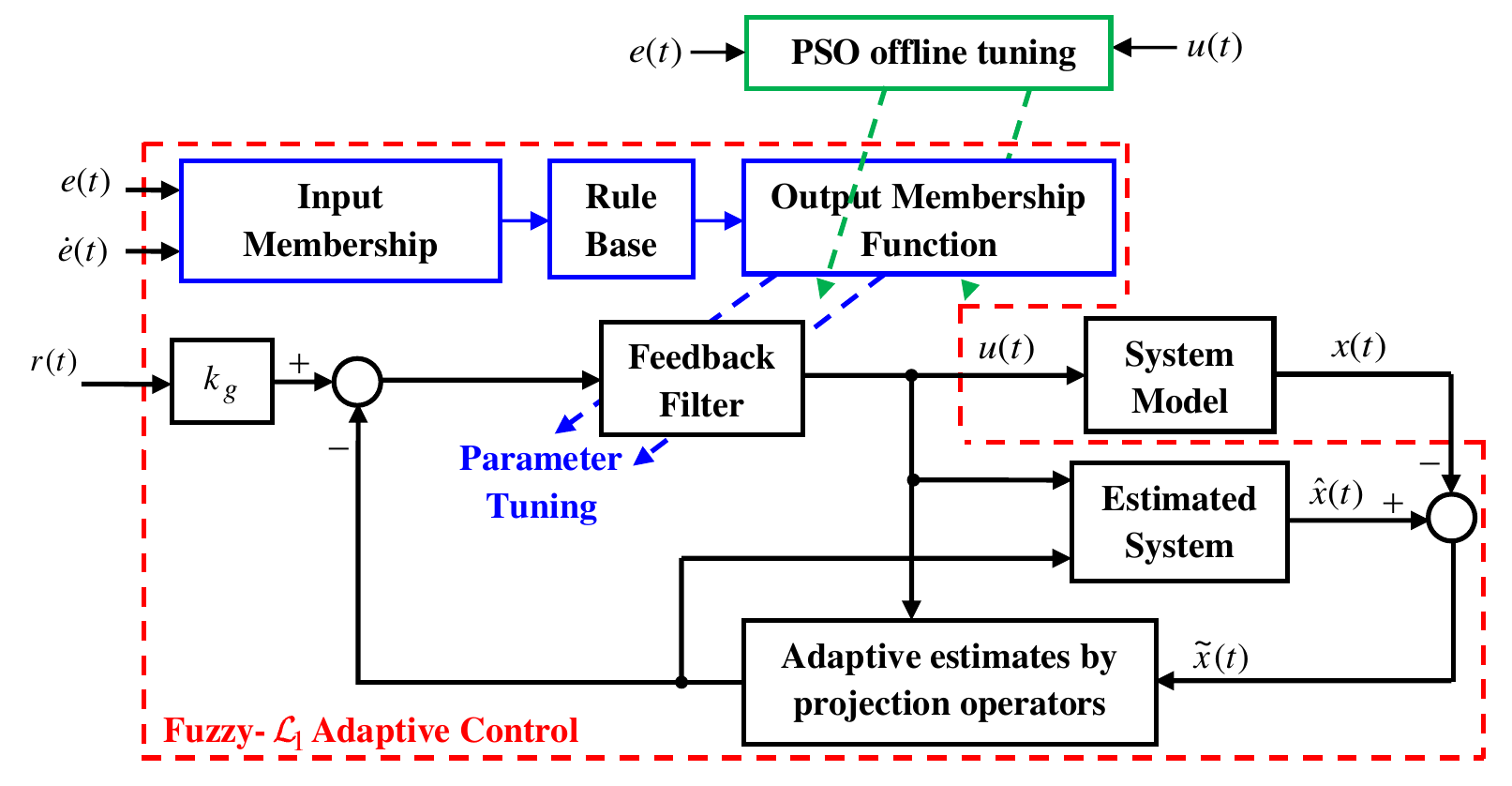}
               \caption{Proposed fuzzy- adaptive control structure.}
               \label{Fuzzy_L1_3PS}
            \end{figure*}

     One can observe that while the previous approaches to determine the optimal coefficients have different degrees of complexity, they agree on the fact that the selection of the appropriate coefficients is performed off-line; and once selected, these coefficients remain unchanged. This study claims that increasing the robustness while guaranteeing fast adaptation requires dynamic and on-line tuning of the feedback filter's coefficients and any proposed method should be relatively simple and easily implementable. To this end, this study proposes fuzzy tuning of the filter's coefficients optimized using PSO taking into account the rate and value of the tracking error between the model reference output and the system's output. The complete structure of fuzzy-\Lone adaptive controller is presented in Figure \ref{Fuzzy_L1_3PS}. The FLC-based tuning is performed on-line during operation. On the other hand, PSO identifies the optimal values of output membership functions through off-line tuning.
            \\
     Fuzzy logic controller (FLC) is classified as an intelligent technique and was first proposed in \cite{zadeh_fuzzy_1965}. FLC showed impressive results in control applications and it has been presented as a robustifying tool with adaptive controllers in \cite{tao_novel_2010},\cite{li_adaptive_2014}. It has been used to compensate unknown nonlinearities of twin rotor MIMO system with adaptive sliding mode control \cite{tao_novel_2010}.In \cite{li_adaptive_2014}, the authors suggested an observer-based adaptive backstepping control scheme and used FLC to approximate unknown uncertainties and to handle bounds of dead zone nonlinearity. On the other hand, evolutionary algorithms are introduced as potential optimization techniques in various control applications. They gained the interest of researchers and witnessed rapid developments over the past few decades. In particular, Particle swarm optimization (PSO) was introduced as a global search technique in \cite{eberhart_new_1995}. PSO has been applied successfully to optimize the structure and parameters of adaptive fuzzy controller in \cite{das_sharma_hybrid_2009} and optimize the variables of FLC membership functions in \cite{bingul_fuzzy_2011}, \cite{wong_pso-based_2008}. The need to tune controller systems with originally fixed coefficients has been widely recognized. In particular, fuzzy tuning has been investigated in several studies (see for instance \cite{precup2014novel}, \cite{valdez2014survey}, and \cite{kumar2014ann}) and controllers based on such approach have been implemented in many applications (see for instance \cite{kumar2014ann}, \cite{mendes2015indirect}, \cite{masumpoor2015adaptive}, \cite{zhang2015neuro}). This allows to conclude that the proposed approach is practical and can definitely be implemented with great benefit.\\
       
      To summarize, in this work, fuzzy-\Lone adaptive controller is proposed to  tune the filter's coefficients in order to improve the trade-off between robustness and fast adaptation.  The coefficients are dynamically tuned and not kept fixed as in the literature, thus allowing for better performance. In the proposed approach, FLC is in charge of online tuning of the filter coefficients taking into account the range and rate of the tracking error. The use of FLC to tune the coefficients improves the stability and the robustness of the system and allows faster closed loop dynamics. Input membership functions and other FLC parameters are assigned arbitrarily while PSO optimizes the optimal variables of the fuzzy output membership functions. The approach is validated using different nonlinear systems and the extensive simulation results are benchmarked to the \Lone adaptive controller with fixed constant gain. The method is simpler than those in the literature and easily implementable.  \\ This chapter is organized as follows: In section two, brief review of \Lone adaptive control including adaptation laws and the general structure is discussed.  Section three presents the idea of filter design and the structure of the proposed control. Section four states the optimization problem and presents the particle swarm optimization algorithm. Illustrative examples will be presented in section five in order to clarify and verify the proposed approach. Finally, last section contains the conclusion.

\section{Review of \Lone adaptive controller} \label{Sec_2}
{
Consider the following dynamics for nonlinear system
    \begin{equation}
     \label{eq:ch3L1sys}
     \begin{aligned}
       &\dot{x}(t) = A_{m}x(t) + b(\omega u(t) + f(x(t),u(t),t))\\
       & y(t) = cx(t)
     \end{aligned}
    \end{equation}

    where $x(t) \in \mathbb{R}^n$ is the system state vector (assumed measured); $u(t) \in \mathbb{R}$ is the control input; $y(t) \in \mathbb{R}$ is the system output; $b,c \in {R}^{n}$ are constant vectors (known);  $A_m $ is $\mathbb{R}^{n \times n}$ Hurwitz  matrix  (known) refers to the desired closed-loop dynamics; $\omega(t) \in \mathbb{R}$ is an unknown time variant parameter describes unmodeled input gain with known sign, and $f(x(t),u(t),t): \mathbb{R}^n\times \mathbb{R} \times \mathbb{R} \to \mathbb{R}$ is an unknown nonlinear continuous function.
    \begin{assum} 
    {\bf(Partially known with known sign control input)} Let the upper and the lower input gain bounds be defined by $\omega_l$ and $\omega_u$ respectively, where
     \begin{equation*}
       \omega \in \Omega \triangleq [\omega_l , \omega_u],\hspace{10pt} |\dot{\omega}| < \omega
     \end{equation*}
     $\Omega$ is assumed to be known convex compact set and $0<\omega_l<\omega_u$ are uniformly known conservative bounds.
    \end{assum}     
     
    \begin{assum} 
    {\bf (Uniform boundedness of $ f(0,u(t),t)) $} Let $B>0$ such that $ f(0,u(t),t)) \leq B $ for all  $ t \geq 0 $
    \end{assum} 
    \begin{assum} 
     {\bf (Partial derivatives are semiglobal uniform bounded)} For any $\delta >0$, there exist $d_{f_{x}} (\delta) >0 $ and $d_{f_{t}} (\delta) >0 $  such that for arbitrary $||x||_{\infty} \leq \delta$ and any $u$, the partial derivatives of $ f(x(t),u(t),t)) $ is piecewise-continuous and bounded,
      \begin{equation*}
        ||\frac{\partial f(x(t),u(t),t) }{\partial x}|| \leq d_{f_{x}} (\delta),\hspace{10pt} |\frac{\partial f(x(t),u(t),t) }{\partial t}| \leq d_{f_{t}} (\delta)
      \end{equation*}  
    \end{assum} 
    \begin{assum} 
      {\bf (Asymptotically stable of initial conditions)} The system assumed to start initially with $x_0$ inside an arbitrarily known set $\rho_0$ i.e., $||x_0||_{\infty} \leq \rho_0 < \infty$.
    \end{assum}    
        \begin{equation}
          \label{eq:chfuzL1Max}
          \theta_b \triangleq d_{f_{x}}(\delta),\hspace{10pt} \Delta \triangleq B + \epsilon
        \end{equation}
         {\bf Lemma:} If $||x||_{\mathcal{L}_{\infty}} \leq \rho$ and there exist $u(\tau)$, $\omega(\tau)$, $\theta(\tau)$and  $\sigma(\tau)$ over $[0,t]$ such that
         \begin{equation}
           \omega_l<\omega<\omega_u
         \end{equation} 
         \begin{equation}
          \vspace{-2mm}
           |\theta(\tau)| < \theta_b
         \end{equation}
         \begin{equation}
          \vspace{-2mm}
           |\sigma(\tau)| < \sigma_b
         \end{equation}
         \begin{equation*}
          \vspace{-2mm}
           f(x(\tau),u(\tau),\tau) = \omega u(\tau) + \theta(\tau) ||x(\tau)||_{\infty} +  \sigma(\tau)
         \end{equation*}
         If $\dot{x}(\tau)$ and $\dot{u}(\tau)$ are bounded then $\omega(\tau)$, $\theta(\tau)$and $\sigma(\tau)$ are differentiable with finite derivatives.\\

       
   The \Lone adaptive controller is composed of three parts defined as the state predictor, the adaption algorithm based on projection and the feedback filter (see Figure (\ref{Fuzzy_L1_Gen})). 
  The main function of the state predictor is developed based on the adaptation laws
           \begin{equation}
             \label{eq:ch3L1est}
             \begin{aligned}
               &\dot{\hat{x}}(t) = A_{m}\hat{x}(t) + b(\hat{\omega} u(t) + \hat{\theta} ||x(t)||_{\infty}+\hat{\sigma} )\\
               & \hat{y}(t) = c\hat{x}(t)
             \end{aligned}
           \end{equation}
          The adaptive estimates $\hat{\omega} \in \mathbb{R}$, $\hat{\theta} \in \mathbb{R}$ and $\hat{\sigma} \in \mathbb{R}$ are defined as follows
           \begin{equation}
             \label{eq:chFuzL1Proj}
             \begin{aligned}
               &\dot{\hat{\omega}} = \Gamma Proj(\hat{\omega},-\tilde{x}^{\top}Pbu(t)), \hspace{10pt} \hat{\omega}(0) = \hat{\omega}_0  \\
               &\dot{\hat{\theta}} = \Gamma Proj(\hat{\theta},-\tilde{x}^{\top}Pb||x(t)||_{\infty}) \hspace{10pt} \hat{\theta}(0) = \hat{\theta}_0\\
               &\dot{\hat{\sigma}} = \Gamma Proj(\hat{\sigma},-\tilde{x}^{\top}Pb) \hspace{10pt} \hat{\sigma}(0) = \hat{\sigma}_0
             \end{aligned}
           \end{equation}
          where $\tilde{x} \triangleq \hat{x} - x(t)$, $\Gamma \in \mathbb{R}^{+}$ is the adaptation gain, and the solution of Lyapunov equation $A_m^TP+PA_m=-Q$ with symmetric  $P>0$ and $Q>0$. The projection operator ensures that $\hat{\omega} \in \Omega \triangleq [\omega_l , \omega_u]$, $\hat{\theta} \in \Theta \triangleq [-\theta_b,\theta_b]$, $|\hat{\sigma}| \leq \Delta$  with $\theta_b$ and $\Delta$ being defined in \eqref{eq:chfuzL1Max}. Projection operators will be evaluated as defined in \cite{pomet_adaptive_1992}\\
          \\
       With special interest to this paper, the control law is defined as 
          \begin{equation}
            u(s) = -k\,D(s)(\hat{\eta}(s)-k_g\,r(s))
          \end{equation}
          
          where $k > 0$ is a feedback gain and $D(s)$ is a strictly proper transfer function leading to a strictly proper and stable transfer function. The Laplace transforms of $r(t)$ and $\hat{\eta}(t) = \hat{\omega} u(t) + \hat{\theta} x(t)+\hat{\sigma}$ are $r(s)$ and $\hat{\eta}(s)$. Finally, $k_g$ is a necessary feedforward gain ensuring a unity steady state gain where $k_g \triangleq -1/(cA_{m}^{-1}b)$ ; $k > 0$. Thus, after a certain transient determined by its bandwidth, the effect of the filter will vanish from the dynamic of the closed loop system.\\
          Thus, in this case, the filter
          \begin{equation}
            C(s) = \frac{\omega \,k\,D(s)}{1+\omega \,k\,D(s)}
          \end{equation}       
          With DC gain $C(0) = 1$. The general structure of \Lone adaptive controller is depicted in Figure \ref{Fuzzy_L1_Gen}.
          \begin{rem}
          The main objective of this work is to design a FLC in order to tune the feedback gain of \Lone adaptive controller and ensure that $y(t)$ tracks a continuous reference signal $r(t)$. In addition, it is aimed at improving the robustness and tracking capability and reducing the control signal range when compared to \Lone adaptive controller with constant parameters.
          \end{rem}

\section{Optimal Fuzzy-tuning of the feedback filter} \label{Sec_3}
      FLC has been used widely for various control applications. In this work, FLC is developed in order to tune the feedback filter gain of the \Lone adaptive controller. The importance of tuning this filter is crucial to improve the robustness and to reduce the control signal range.

\subsection{Structure of Fuzzy Logic Controller}
      The error $e(t)$ is the difference between reference input $r(t)$ and regulated output $y(t)$. $k_p$ and $k_d$ are proportional and differential weights respectively. These parameters will be assigned before designing the membership functions and their values rely on the expected range of both $e(t)$ and $\dot{e}(t)$ in order to normalize fuzzy input between 1 and 0.
          \begin{equation}
            k_p \leq \frac{1}{||e||_{\infty}}, \hspace{10pt} k_d \leq \frac{1}{||\dot{e}||_{\infty}}
          \end{equation}  
     The existence of these norms is guaranteed by \Lone adaptive controller in case of stable dynamics. In addition, they can also be dynamically assigned. The fuzzy filter has a triangular membership functions for both inputs and output. The fuzzy filter has two inputs represented by the error and its rate and one output which is the inverse of the feedback gain $k_f$. Fuzzy inputs are the absolute values of $e(t)$ and $\dot{e}(t)$ multiplied by weighted gains $k_p$ and $k_d$. \Lone adaptive controller will consider the fuzzy output $k_f$ as a feedback gain if the error is greater than $k_e$. Adversely, the controller will consider a constant feedback gain $k$ if the error is less than or equal $k_e$ as shown in figure \ref{Fuzzy_L1_4}.\\
       
\section{Particle Swarm Optimization} \label{Sec_4}
       Particle swarm optimization is an intelligent evolutionary computation algorithm.  PSO algorithm deploys a set of particles in the space as a population and each particle is a candidate solution. Each particle in the search space moves randomly in swarm of particles to find the optimal solution. Each solution is defined by a particle position in the space and the velocity of swarming is necessary to target the best position. The proper setting of the algorithm variables ensures swarming in the vicinity space of the optimal solution and increases the probability of fast convergence. The velocity and position of the particle are defined according to the following two equations \eqref{eq:chFuzL1v} and \eqref{eq:chFuzL1x} respectively
          \begin{equation}
          \label{eq:chFuzL1v}
          \begin{split}
           v_{i,j}&(t) = \alpha(t)v_{i,j}(t-1) + c_{1}r_{1}(x_{i,j}^{*}(t-1) \\
           &-x_{i,j}(t-1))+ c_{2}r_{2}(x_{i,j}^{**}(t-1)-x_{i,j}(t-1))
          \end{split}
          \end{equation}
          \begin{equation}
          \label{eq:chFuzL1x}
           x_{i,j}(t) = v_{i,j}(t) + x_{i,j}(t-1))
          \end{equation}
          
       where $i=1,2,\cdots, N_p$ and $N_p$ is the population size, $j=1,2,\cdots, P_s$ and $P_s$ is the number of parameters In each particle, $x_{i,j}^{*}$ and  $x_{i,j}^{**}$ represent the local and global solutions respectively, $\alpha(t)$ is an exponential decreasing inertia, $c_1$ and $c_2$ represent personal and social influence of parameters and finally $r_1$ and $r_2$ are random numbers where $r_1 ,r_2 \in [0,1]$. The objective function is defined to enhance the tracking capability and improve the control signal range as follows
             \begin{equation}
               \label{eq:chFuzL1Obj}
               \begin{aligned}
                Obj = \sum\limits_{t=0}^{t_{sim}}\big(\gamma_1 e^2(t)+\gamma_2 u^2(t)\big)
               \end{aligned}
             \end{equation}
       where $e(t) = r(t) - y(t)$, $e(t)$ and $u(t)$ are the system error and control signal respectively. $\gamma_1$ and $\gamma_2$ are weights that can be selected arbitrarily. Obviously, the output membership functions have 18 parameters and they should be optimized to minimize the objective function. Particle swarm optimization is developed to search for the optimal values of aforementioned parameters. It must be noted that triangular nodes of output membership function represent position $x_{i,j}$, each two triangular intersect on the horizontal axis on one node. The computational flow diagram of PSO algorithm is illustrated in Figure \ref{Fuzzy_L1_PSO}. The algorithm will be used with \Lone adaptive controller to define the optimal parameters of output membership functions for a specific number of generations as mentioned in \cite{eberhart_new_1995}, \cite{abido_optimal_2002}.
       \begin{figure}[!h]
          \centering
             \includegraphics[scale=0.6]{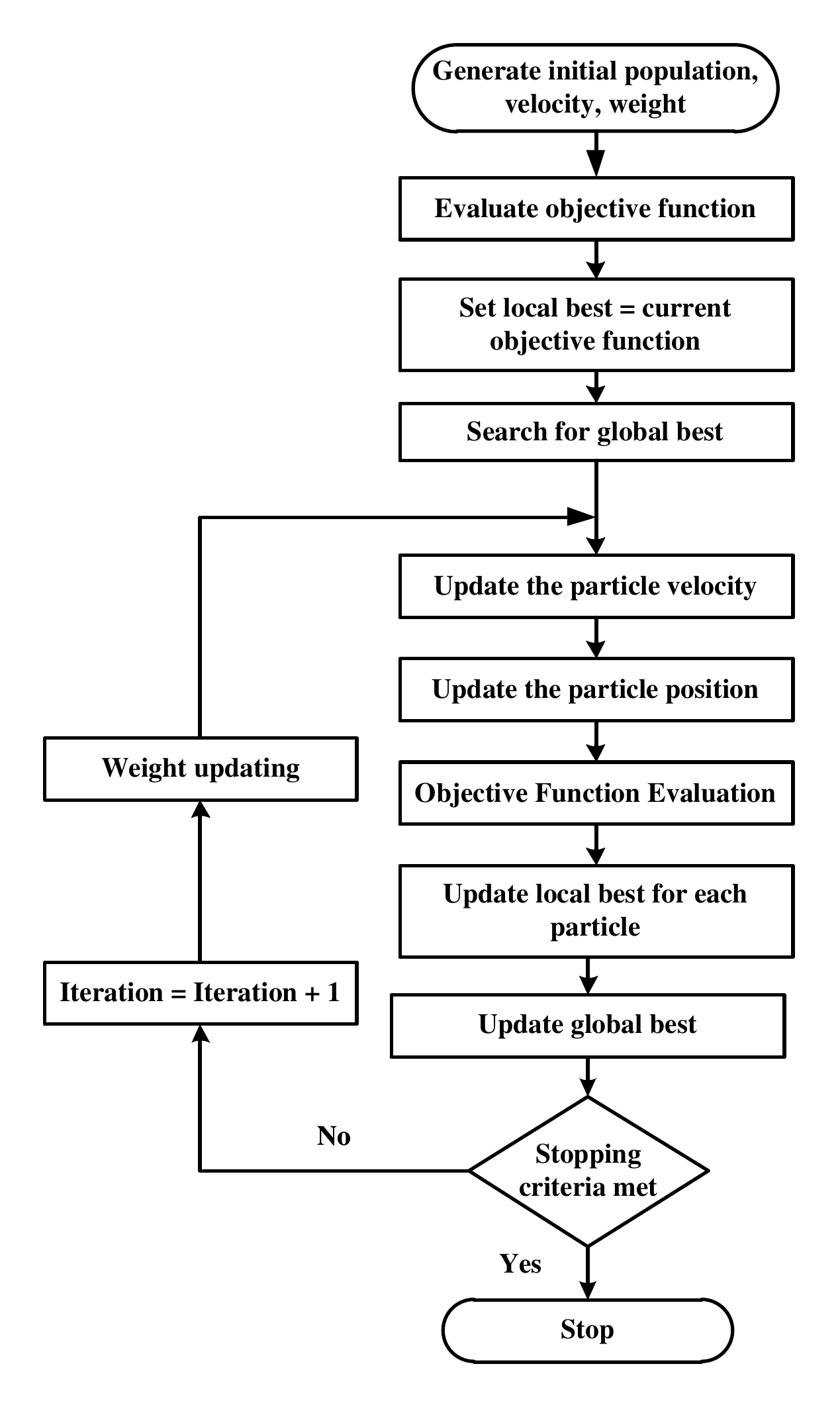}
          \caption{Flowchart of particle swarm Optimization \cite{abido_optimal_2002}.}\label{Fuzzy_L1_PSO}
       \end{figure}
       
\begin{rem}
In the proposed approach, the properties of the filter, such as strictly proper, low pass  with C(0)=1, are preserved. Consequently, stability of the Fuzzy-based-\Lone adaptive controller is guaranteed by the same analysis of stability done in \cite{cao_design_2006}.  
\end{rem}       
      
\section{Results and Discussions} \label{Sec_5}
\subsection{Fuzzy \Lone adaptive controller implementation:}
    Problem in \cite{hovakimyan_l1_2010} has been considered here with additive nonlinearities added to the system as follows
    \begin{equation*}
     \begin{aligned}
       &\dot{x}(t) = A_{m}x(t) + B(\omega u(t) + f(x(t),t))\\
       & y(t) = Cx(t)
     \end{aligned}
    \end{equation*}
    where $x(t) = [x_1(t),x_2(t)]^{\top}$ are the system states, $u(t)$ is the control input, $y(t)$ is the regulated output and $f(t,x(t))$ includes high nonlinearity assumed to be unknown. In addition,
     \begin{equation*}
      A = 
      \begin{bmatrix}
      0 & 1\\
      0 & 0
      \end{bmatrix}, \hspace{10pt}
      B = 
      \begin{bmatrix}
      0\\
      1
      \end{bmatrix}, \hspace{10pt}
      C = 
      \begin{bmatrix}
      0 & 1
      \end{bmatrix}
     \end{equation*}
     and
    \begin{equation*}
       f(x(t),t) = 2x_1^2(t)+2x_2^2(t)+x_1sin(x_1^2)+x_2cos(x_2^2)\\
    \end{equation*}
    \begin{equation*}
       \omega = \frac{75}{s+75}\\    
    \end{equation*}
    $\omega$ is a function with fast dynamic to ensure smoothness of the control signal. The compact sets of the projection operators for unmodeled input parameters, uncertainties and disturbances were assigned to $[\omega_{min},\omega_{max}] \in [0,10]$, $\Delta = 100$ and $\theta_b = 10$ . The control objective is to design a fuzzy-\Lone adaptive controller to enhance each of control signal range and tracking capability of a bounded reference input $r(t)$ for the output signal $y(t)$. Desired poles are set to = $-21 \pm j0.743$, the constant feedback gain($k$) = 20, the adaptation gain($\gamma$) = 1000000 and $Q = \big[ \begin{smallmatrix} 1 & 0 \\ 0 & 1 \end{smallmatrix}\big]$. Fuzzy control parameters are $k_p = 0.1$ , $k_d = 0.05$ and $k_e = 0.1$. Figure \ref{Fuzzy_L1_4} illustrates the FLC with \Lone adaptive controller.
          \begin{figure*}[ht]
          \centering
          \includegraphics[width=12cm, height=7cm]{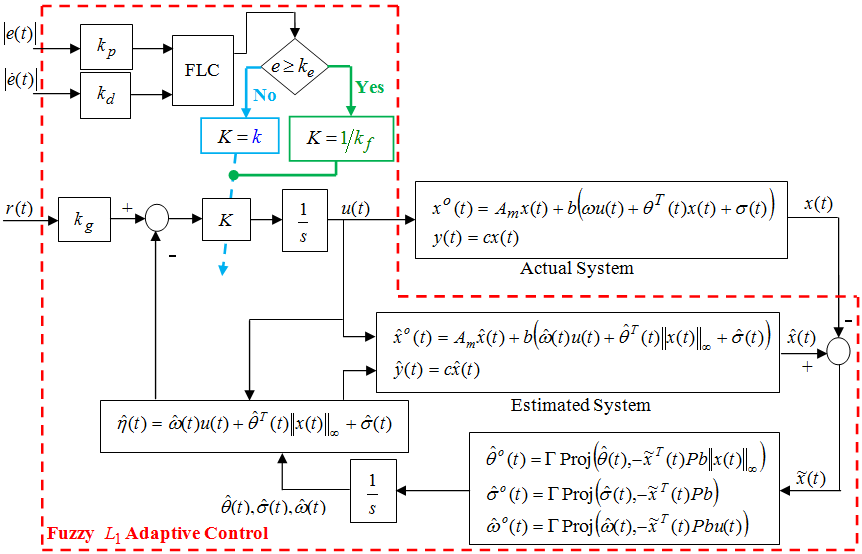}
          \caption{Fuzzy-\Lone adaptive controller for nonlinear SISO system.}\label{Fuzzy_L1_4}
       \end{figure*}

  \subsection{Membership Function Optimization}
      The objective of this work is to construct output membership function for FLC capable of reducing the error and the control signal. Values of input membership functions and constraints of the output membership functions were chosen based on trying different values by running a certain number of experiments. The range of input membership functions was adjusted between 0.08 and 1 and their values were selected as shown in Figure \ref{Fuzzy_L1_mem_e}. The fuzzy inputs and output have triangular membership functions with five linguistic variables. Linguistic variables are assigned as very large ($VL$), large ($L$), small ($S$), very small ($VS$) and zero ($Z$) where values of input membership function will be assigned arbitrarily. Values of output membership functions are  optimized using PSO. Rule base of the proposed filter is demonstrated in Table \ref{table:Tab_Rule}.
       \begin{table}[!t]
       \caption{Rule base of FLC.} 
       \centering 
       \small
       \begin{tabular}{|c| c| c| c| c| c|} 
       \hline\hline 
       {\bf $\Delta e/e$} & {\bf VL}  & {\bf L} & {\bf S} & {\bf VS} & {\bf Z} \\ [0.0ex]
       \hline\hline 
       {\bf VL} & $VL$ & $VL$ & $VL$ & $VL$ & $L$ \\[0ex]
       \hline                   
       {\bf L}  & $VL$ & $VL$ & $VL$ & $L$ & $S$ \\[0ex] 
       \hline
       {\bf S}  & $VL$ & $VL$ & $L$ & $S$ & $VS$ \\[0ex] 
       \hline
       {\bf VS} & $VL$ & $L$ & $S$ & $VS$ & $VS$ \\[0ex]
       \hline
       {\bf Z}  & $L$ & $S$ & $VS$ & $VS$ & $Z$ \\[0ex] 
       \hline\hline 
       \end{tabular}
       \label{table:Tab_Rule}
       \end{table}
       
        \begin{figure}[!h]
           \centering
             \includegraphics[scale=0.3]{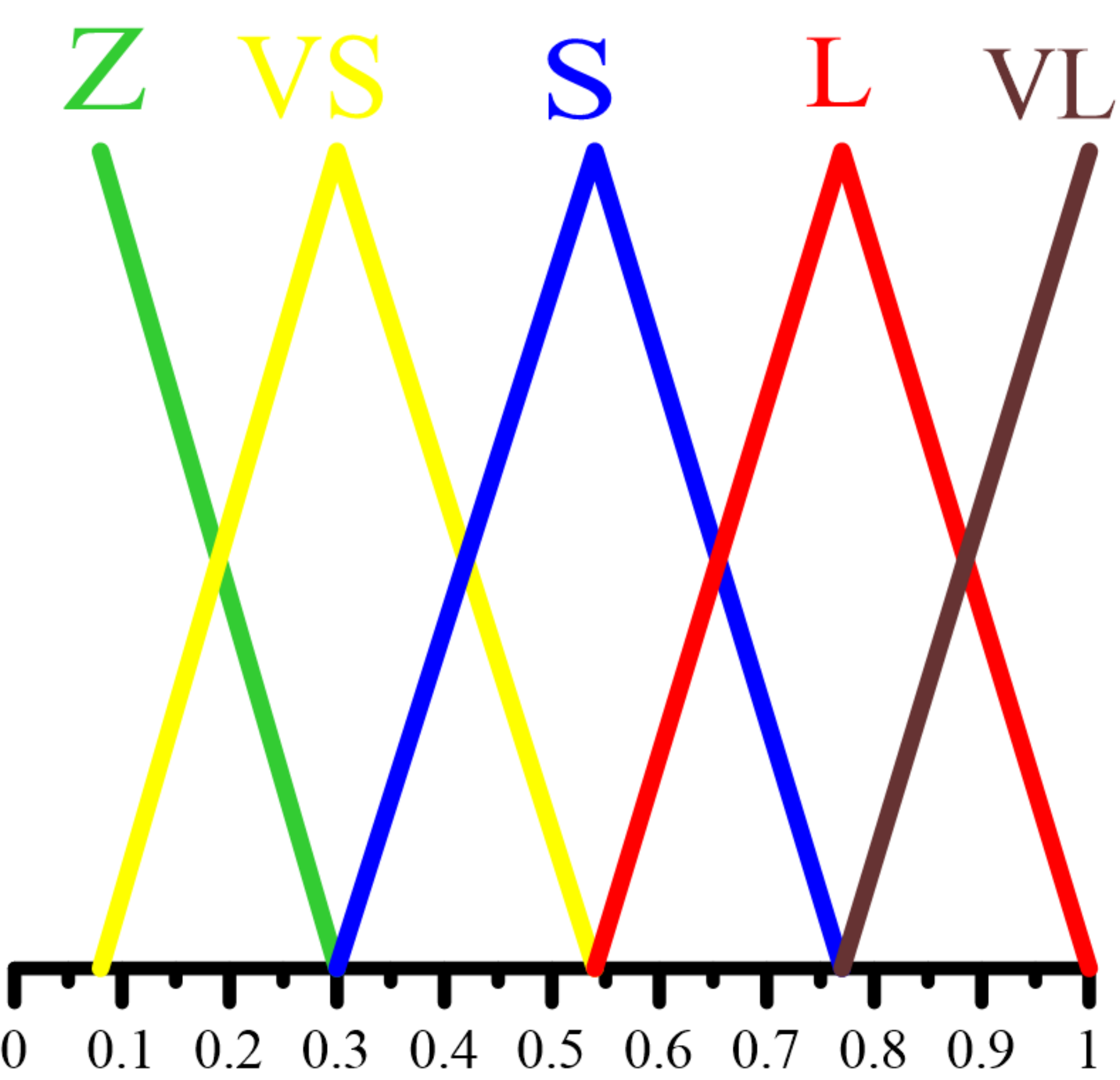}
           \caption{Error and rate of error membership functions.}\label{Fuzzy_L1_mem_e}
        \end{figure}
        
      Constraint values of output membership functions are represented by three parameters as lower ($l$), center ($c$) and higher ($h$) values. These three parameters of each triangular membership function will constrain between minimum and maximum bounds. Constraints bounds of the problem can be defined as follows 
      \begin{equation}
        \begin{aligned}
         &[4,8,8] \leq [VL_l,VL_c,VL_h] \leq [8,12,12]\\
         &[1.5,3,6] \leq [L_l,L_c,L_h] \leq [3,6,10]\\
         &[0.3,1.5,4] \leq [S_l,S_c,S_h] \leq [1.5,4,8]\\ 
         &[0,0.5,1.5] \leq [VS_l,VS_c,VS_h] \leq [0.5,1.5,3]\\  
         &[0.0,0.0,0.3] \leq [Z_l,Z_c,Z_h] \leq [0.0,0.0,1.5]\\ 
        \end{aligned}
      \end{equation}
     With $VL$, $ML$, $L$, $S$, $MS$, $VS$ and $Z$ were mentioned before as a linguistic variables. Also, we assigned $VL_c = VL_h$, $VL_l = S_h$, $L_l = VS_h$, $S_l = Z_h$, $VS_l = z_c$ and  $z_c = Z_l$ which means that we have only nine parameters to be optimized.
 \subsection{PSO Simulation results}
     The population size is set arbitrarily as 150 particles and each particle include 9 parameters will be optimized based on a minimization objective function and these parameter are $VL_c$, $VL_l$, $L_l$, $L_c$, $L_h$, $S_l$, $S_c$, $VS_l$ and $VS_c$ in (\ref{eq:chFuzL1Obj}). The initial settings of PSO algorithm are demonstrated in Table \ref{table:Tab_PSO} and the maximum numbers of generations is 100.
       \begin{table}[!h]
       \setlength{\tabcolsep}{5pt}
       \caption{Parameters setting for PSO.} 
       \centering 
       \small
       \begin{tabular}{|c| c| c| c| c| c|} 
       \hline\hline 
       {\bf Parameter} & $\lambda$  & $\alpha$ & $c_1$ & $c_2$ \\ [0.0ex]
       \hline\hline 
       {\bf Settings} & 10 & 0.99 & 2 & 2  \\[0ex]
       \hline\hline 
       \end{tabular}
       \label{table:Tab_PSO}
       \end{table}
  
    \subsection{PSO Results}
    The system was simulated for 8 seconds and the data was recorded every 0.01 seconds. The reference input was defined by $cos(0.5t)$ with zero initial conditions. The optimal variables of output triangular membership functions are illustrated in Figure \ref{Fuzzy_memb_u}. The fitness reduction during the search process is demonstrated in Figure \ref{Fuzzy_L1_PSOObj}. However, it is clear that objective function is reduced significantly and enormously to a suitable value which is reflected on the output performance as revealed in Figure \ref{Fuzzy_L1_out1}. Figure \ref{Fuzzy_L1_out1}.(a) demonstrates the optimal output performance and Figure \ref{Fuzzy_L1_out1}.(b) shows the control signal of the considered problem.\\
       \begin{figure}[!h]
          \centering
          \includegraphics[scale=0.3]{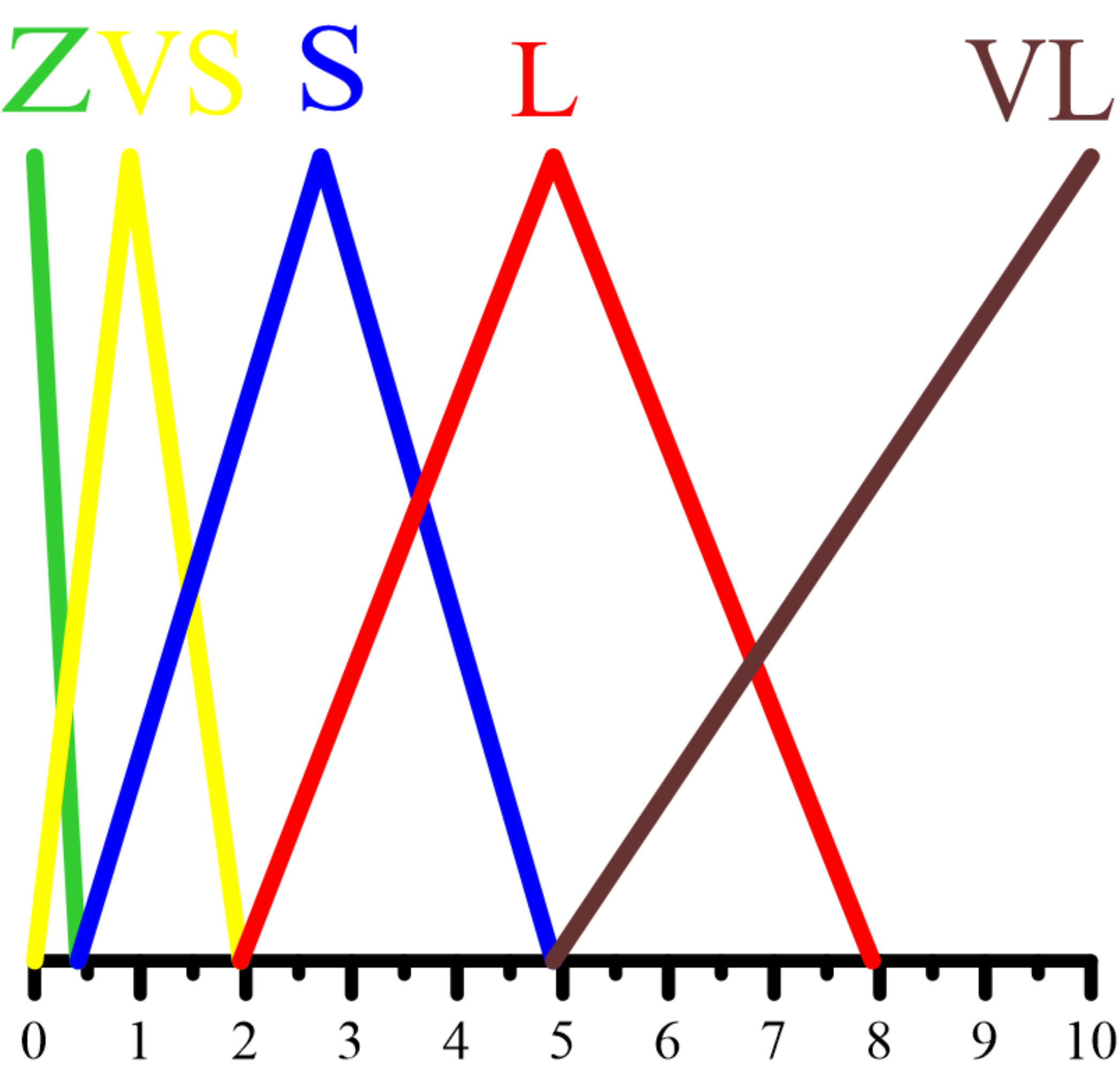}
          \caption{Graphical illustration of output membership functions.}\label{Fuzzy_memb_u}
       \end{figure}
          \begin{figure*}[ht]
          \centering
          \includegraphics[scale=0.6]{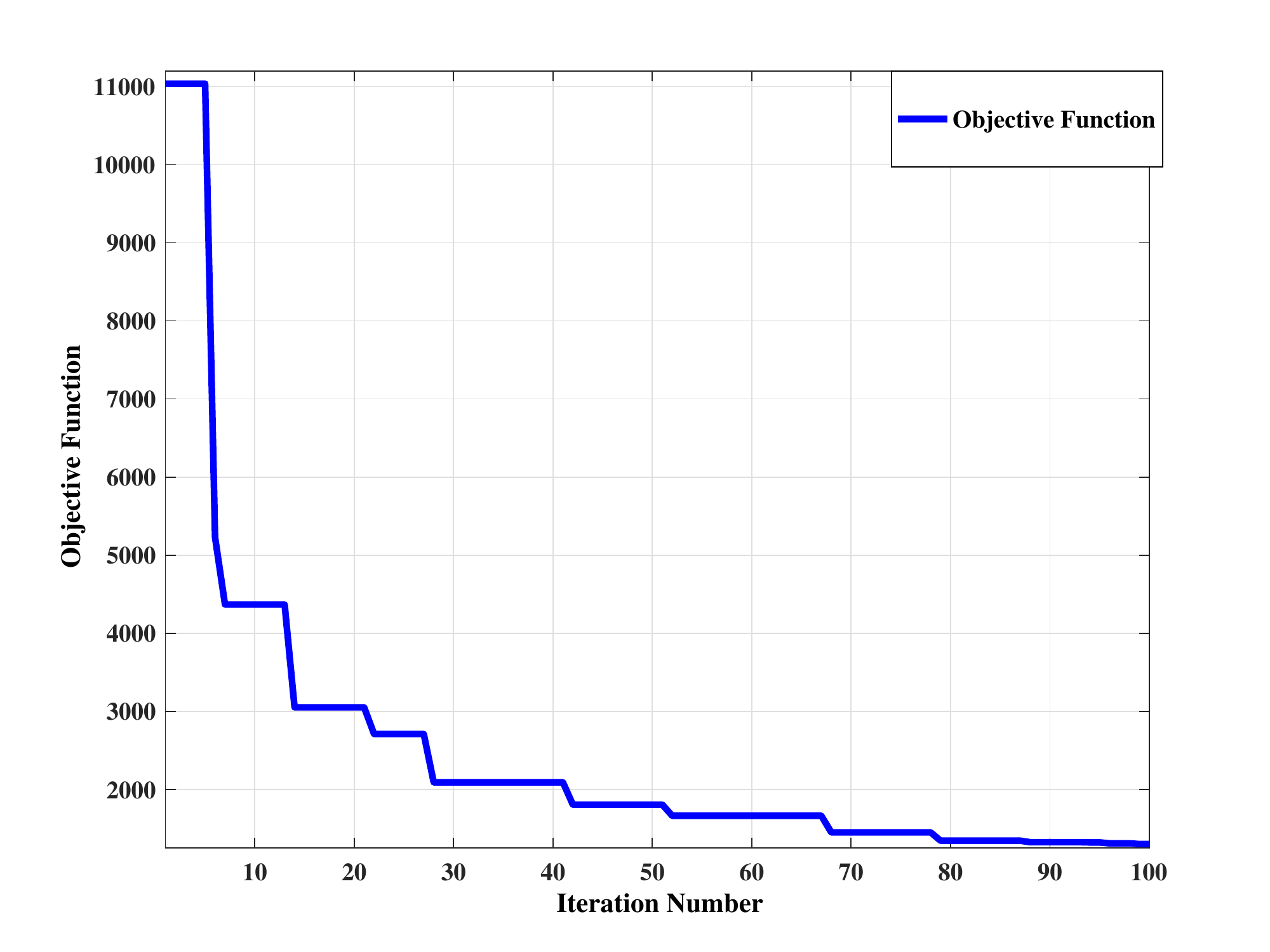}
          \caption{Objective function minimization with PSO search process.}\label{Fuzzy_L1_PSOObj}
       \end{figure*}
          \begin{figure*}[ht]
          \centering
                    \includegraphics[scale=0.6]{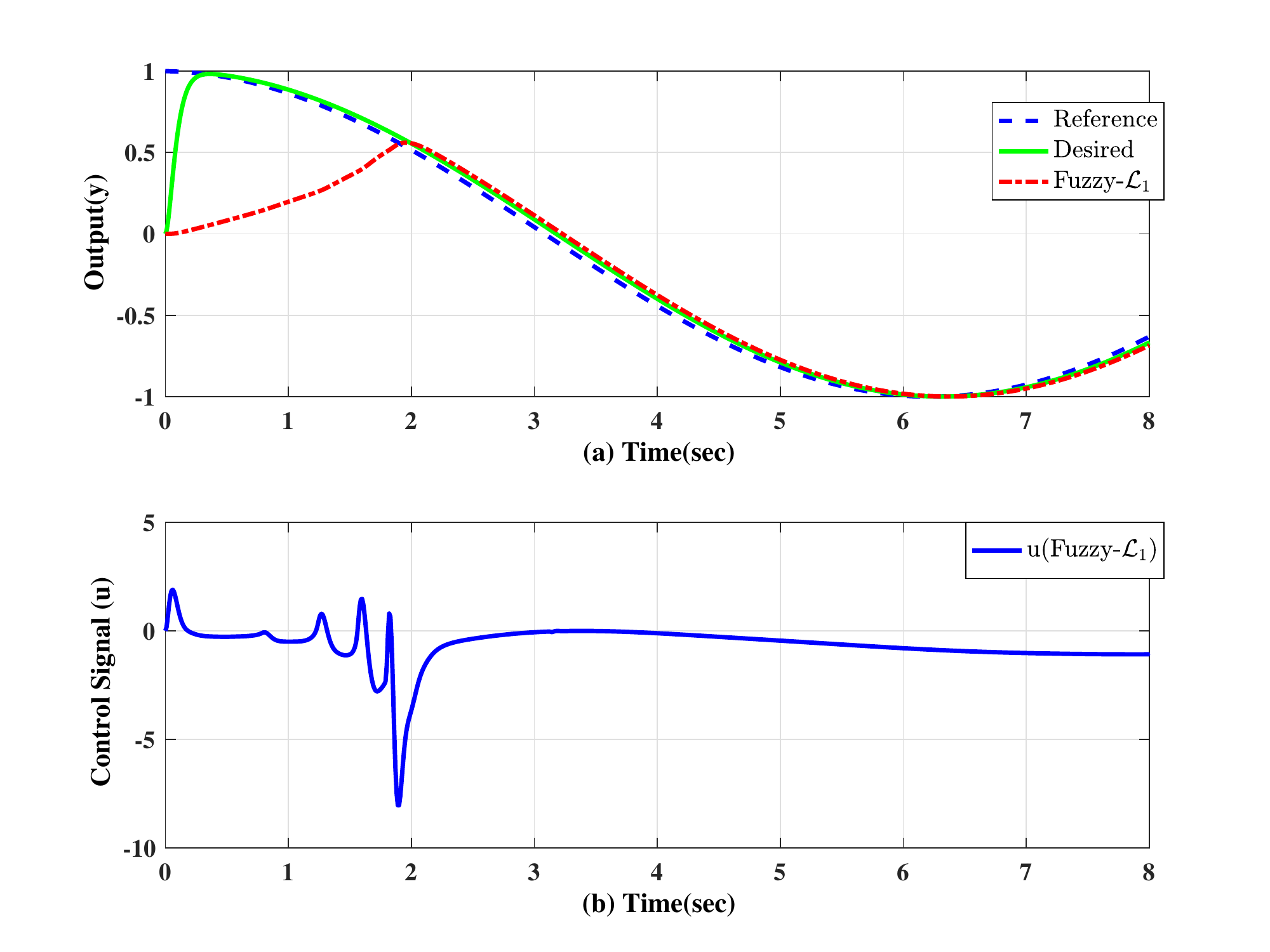}
          \caption{Performance of fuzzy-\Lone adaptive controller after 100 iterations search process.}\label{Fuzzy_L1_out1}
       \end{figure*}
     In this study, three different scenarios are considered to demonstrate the robustness of fuzzy-\Lone adaptive controller. All cases will be simulated for 40 seconds.
     The first case will discuss the nonlinear system included in the search process. Case 2 includes the nonlinear model with high uncertainties, unmodeled input parameters and adding some disturbances in order to validate the robustness of fuzzy filter with \Lone adaptive controller. Case 3 consider all assumptions in case 2 in addition to investigate the system with faster desired closed loop dynamics.\\
     {\bf Case 1:} Figure \ref{Fuzzy_L1_out2} presents the output performance of fuzzy-\Lone adaptive controller versus \Lone adaptive controller and their control signals. Fuzzy-\Lone adaptive controller guarantees uniform transient and smooth tracking performance. In addition, its major contribution lies in reducing the control signal range by tuning the feedback gain. Tuning feedback gain enhances the robustness of the system and reduces the control signal range. The correspondence difference of feedback gain between fuzzy-\Lone adaptive controller and \Lone adaptive controller is illustrated in Figure \ref{Fuzzy_L1_out2_2}.(a). The errors of both controllers are presented in Figure \ref{Fuzzy_L1_out2_2}.(b).
          \begin{figure*}[ht]
          \centering
                    \includegraphics[scale=0.6]{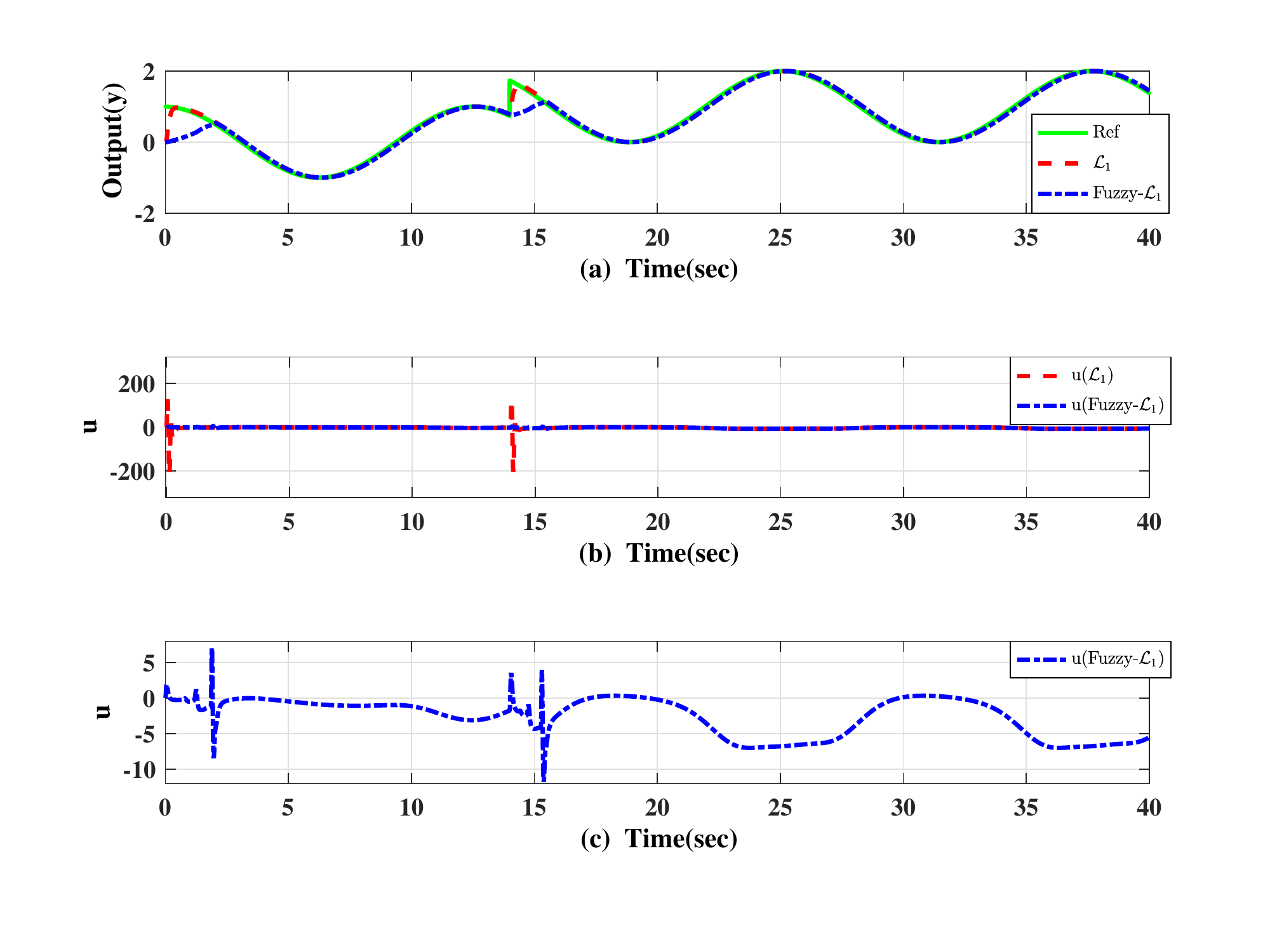}
          \caption{Performance of fuzzy-\Lone adaptive controller and \Lone adaptive controller for nonlinear system of case 1.}\label{Fuzzy_L1_out2}
       \end{figure*}
     
          \begin{figure*}[ht]
          \centering
                              \includegraphics[scale=0.5]{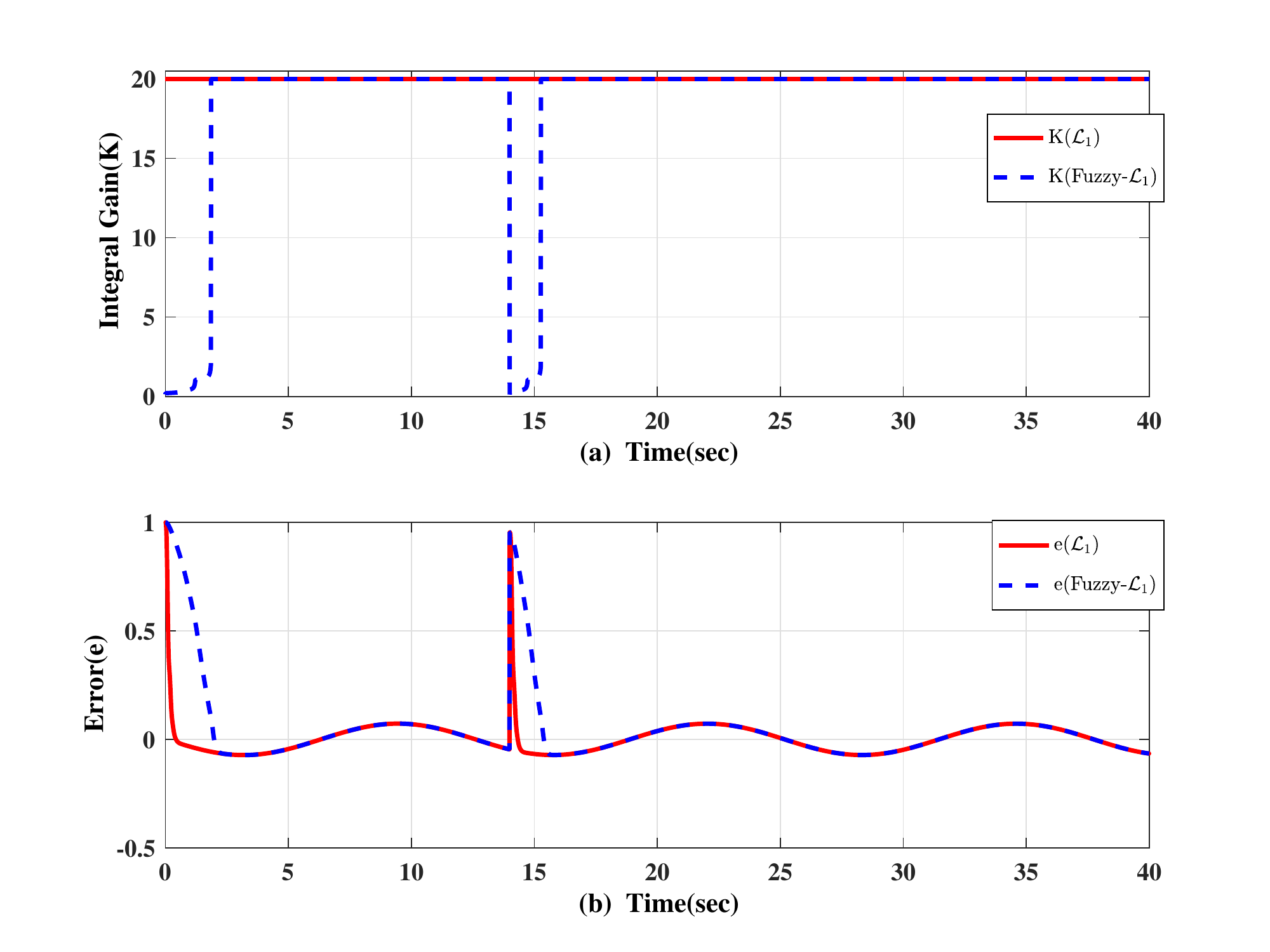}
          \caption{Feedback gain and output error of fuzzy-\Lone adaptive controller and \Lone adaptive controller of case 1.}\label{Fuzzy_L1_out2_2}
       \end{figure*}
       {\bf Case 2:}To illustrate the effectiveness of the proposed fuzzy filter with \Lone adaptive controller, robustness of fuzzy filter is examined against any existing of high uncertainties, unmodeled input parameters and disturbances.Here, the nonlinear model and other assumptions mentioned in case 1 are similar the nonlinear function; however, the nonlinearity includes high time variant uncertainties and disturbances and these changes will except be presented as follows
       \begin{equation*}
       \begin{split}
        f&(x(t),t) = \big(sin(0.4t)+1 \big)x_1^2(t) + \big(2cos(0.35t)+0.5 \big)x_2^2(t)\\
        &+ \big(sin(0.3t)+0.3 \big)x_1sin(x_1^2) + sin(0.35t)cos(0.4t)\\ 
        & + 0.5x_2cos(x_2^2+0.5cos(0.3t))+ sin(0.3t)cos(0.4t)z^2 
       \end{split}
       \end{equation*}
       where
       \begin{equation*}
         z(s) = \frac{s-1}{s^2+3s+2}v(s),\hspace{10pt} v(t) = x_1sin(0.2t)+x_2\\    
       \end{equation*}
       The robustness of fuzzy feedback filter gain with \Lone adaptive controller has been validated in Figure \ref{Fuzzy_L1_out3} and presented versus \Lone adaptive controller. The significant impact and the advantage of fuzzy-\Lone controller is revealed on control signals performance as shown Figure \ref{Fuzzy_L1_out3}. Figure \ref{Fuzzy_L1_out3_2}.(a) presents the performance of feedback gain for fuzzy-\Lone adaptive controller and \Lone adaptive controller. Finally, Figure \ref{Fuzzy_L1_out3_2}.(b) shows the error of both controllers. Uniform transient and tracking capability are validated as shown in Figure \ref{Fuzzy_L1_out2}  and \ref{Fuzzy_L1_out3}. The benefits of fuzzy-\Lone adaptive controller can be summarized in including fast desired dynamics and improving the tracking capability and robustness with less range of control signal.
          \begin{figure*}[ht]
          \centering
                                        \includegraphics[scale=0.5]{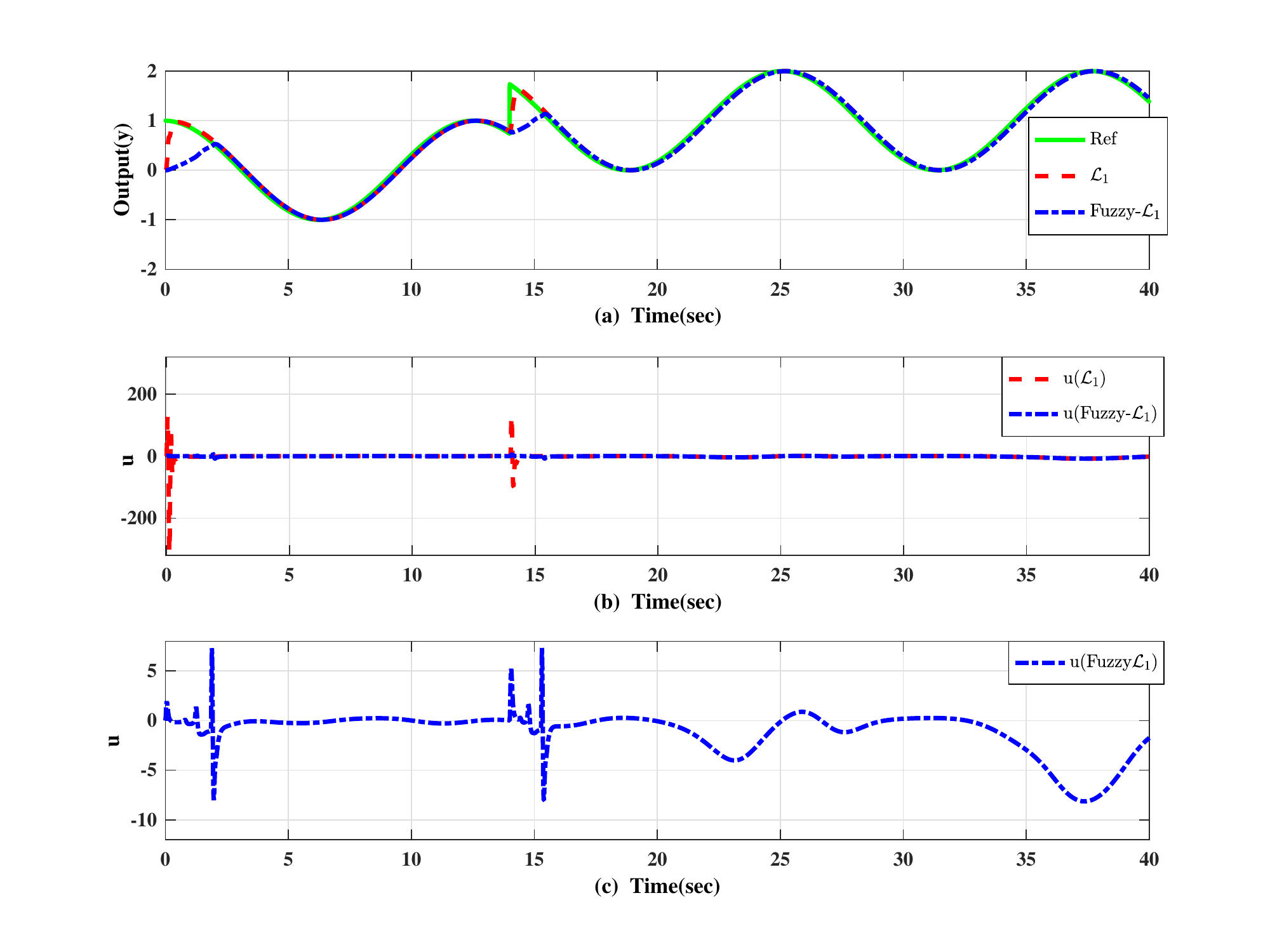}
          \caption{Performance of fuzzy-\Lone adaptive controller and \Lone adaptive controller for nonlinear system of case 2.}\label{Fuzzy_L1_out3}
       \end{figure*}
     
          \begin{figure*}[ht]
          \centering
                                                  \includegraphics[scale=0.5]{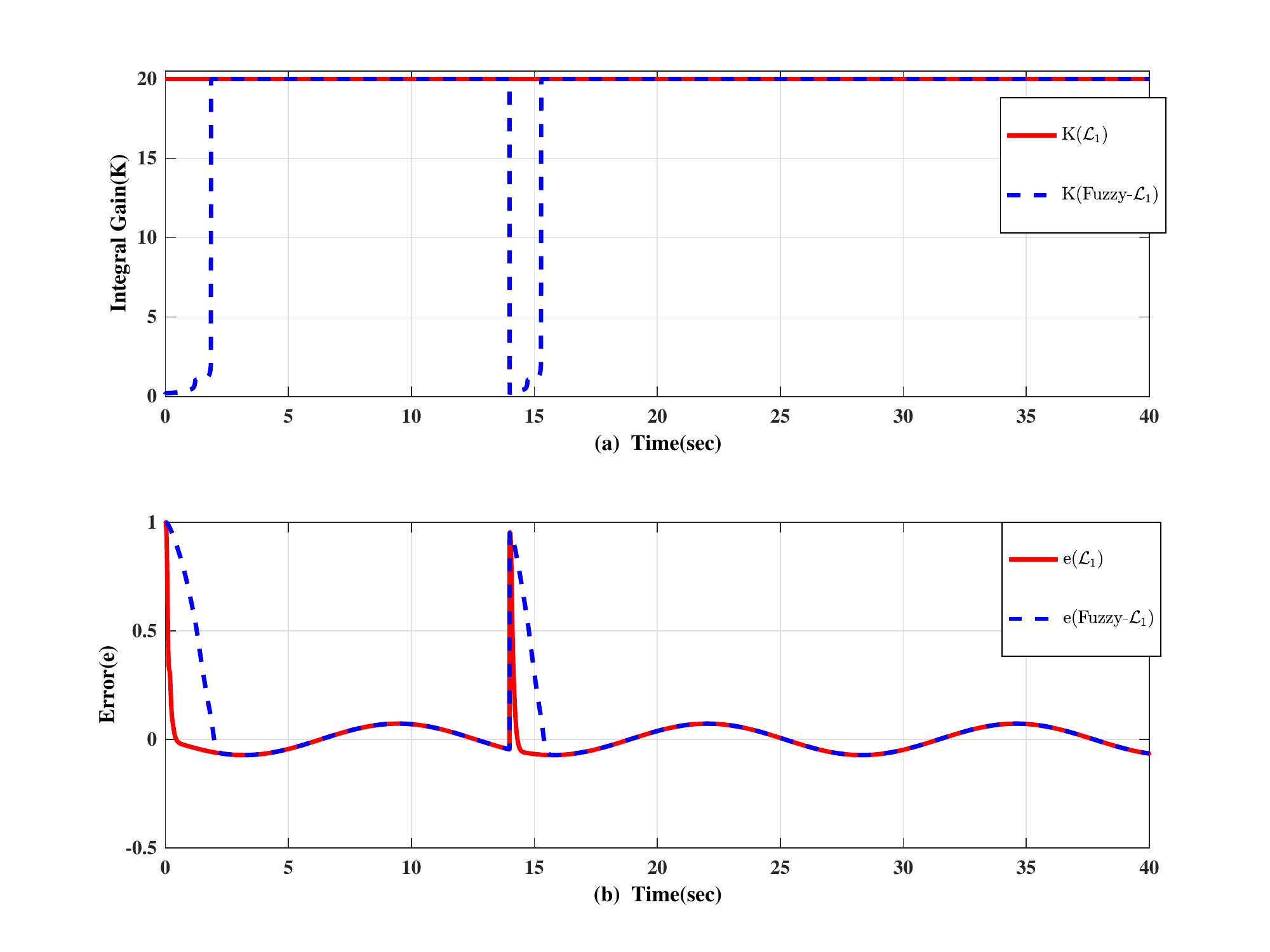}
          \caption{Feedback gain and output error of fuzzy-\Lone adaptive controller and \Lone adaptive controller of case 2.}\label{Fuzzy_L1_out3_2}
       \end{figure*}
       
       {\bf Case 3:} The robustness of fuzzy-\Lone adaptive controller and \Lone adaptive controller will reveal more in this case. All aforementioned assumptions in case 2 are similar here except the desired closed loop dynamics assumed to be faster than case 2. Desired poles are set to $p=-84 \pm j0.743$. According to this change in closed loop poles, the robustness of \Lone adaptive controller will be violated and the system will no longer be stable. However, fuzzy-\Lone adaptive controller will be able to track the output under this new condition with limitation in increasing the control signal range. Figure \ref{Fuzzy_L1_out4} illustrate the output performance of fuzzy-\Lone adaptive controller for case 3.\\
          \begin{figure*}[ht]
          \centering
                    \includegraphics[scale=0.6]{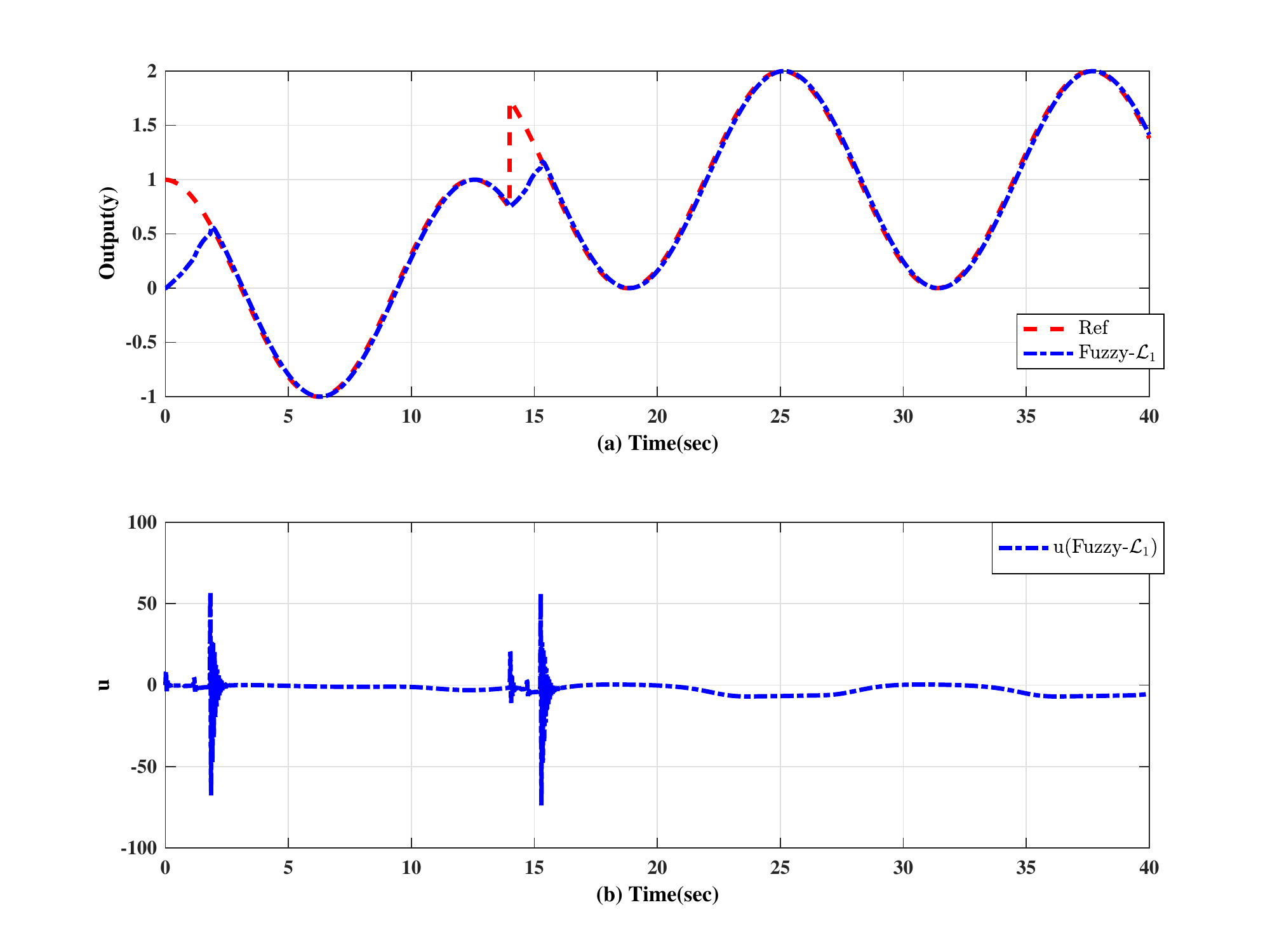}
          \caption{Performance of fuzzy-\Lone adaptive controller for nonlinear system of case 3.}\label{Fuzzy_L1_out4}
       \end{figure*}
        The robustness of this criterion has been simulated and validated with \Lone adaptive controller on high nonlinear system with different forms of nonlinearities and uncertainties in addition to fast closed loop dynamics compared to normal structure of \Lone adaptive controller.  It can be concluded based on the cases considered and results obtained that the proposed fuzzy-based approach to tune the feedback filter improves greatly the performance of \Lone adaptive controller. The proposed fuzzy-\Lone adaptive controller guarantee boundedness of the output and control signal and insures fast tracking and low range of control signal.
       
\section{Conclusion} \label{Sec_6}
     This paper presents a new FLC-PSO design of the feedback gain filter part of  \Lone adaptive controller. PSO determines the optimal variables of the output membership functions. The proposed algorithm tunes on-line the filter parameters, which in turn contributed to improving robustness and stability of \Lone adaptive controller. Moreover, owing to a smooth tuning of the filter the control signal range has been greatly reduced. Illustrative examples were developed and simulated to compare fuzzy-\Lone adaptive controller with \Lone adaptive controller with constant filter parameters and to validate the advantages of the proposed approach. The results show improved performance and robustness with high levels of time variant uncertainties and disturbances in addition to fast desired closed loop adaptation. 
     
     There are several directions for future work. One important area is to implement this approach on a real system and compare the performance with existing techniques. One should note that computation power is relatively cheap and the technology offers several hardware option over which this controller can be implemented. This study aimed at proposing an effective way of tuning the coefficients of the control filter. During this work, the structure of the filter was fixed. Extending this work to determine automatically the appropriate filter's structure has the potential to further improve the robustness and stability. Such an extension will take into account the health of the system and may lead to the design a failure tolerant robust controller.  In this study, PSO has been implemented off-line. To our knowledge, recursive PSO for online implementation is not explored. A recursive and online PSO will impose hard constraints on hardware capacity. The implementation of PSO as it is now is not feasible for online control implementation.  A comparison between PSO and other evolutionary algorithms can be established to define the most effective solution for the fuzzy-\Lone adaptive control problem.

\clearpage

\newpage


\chapter{NEURO-ADAPTIVE FOR STRICT FEEDBACK MIMO SYSTEMS WITH PPF}
\section{Introduction}
  This chapter is mainly concerned in reproducing recent study of robust neuro adaptive control with prescribed performance function on strict feedback MIMO system. The importance of this chapter relies on capturing prescribed performance idea on transient performance, tracking trajectory and smoothness of the control signal. This chapter consists of six sections with first section is an introduction. The second section presents introduction and necessary conditions of prescribed performance. The third section describes the problem formulation and the main idea of prescribed performance function. The fourth section presents neural network for nonlinearity approximation. Section five includes control law formulation and stability analysis. Section six presents simulation and controller benchmark results. The last section is a conclusion.

\section{Introduction of Prescribed Performance}
  Prescribed performance simply means tracking error into an arbitrarily small residual set and the convergence error should be within pre-assigned range. In addition, the convergence rate has to be less than a prescribed value and maximum overshoot should be less than a prescribed constant. Prescribed performance with robust adaptive control was mainly developed to provide a smooth control signal for soft tracking and to solve the problem of accurate computation of the upper bounds for systematic convergence. Due to nonexistence adaptive control for nonlinear systems with error convergence into a predefined small set, the controller with prescribed performance function is demanded. In this chapter, robust adaptive control with prescribed performance should have the ability to approximate the nonlinear model assuming completely unknown dynamics and provide smooth control signal to track the output into the desired trajectory smoothly and accurately.

  The main features of the prescribed performance is its ability of tracking the error into a defined small set. Prescribed performance should guarantee many factors
\begin{itemize}
    \item the convergence of the error within a prescribed bound, 
    \item a maximum overshot  less than a prescribed value,
    \item a uniform ultimate boundedness property for the transformed output error,
    \item adaptive and smooth tracking.
\end{itemize} 

   Neural network will be used to estimate the nonlinear model as an online estimation tool in the adaptive control problem. Adaptive control will be offered to stabilize the system by canceling undesired dynamics using neural network. Also, it will be used to provide robust tracking and forcing the error to be bounded in predefined set. The prescribed set will be reduced into a very small set according to a pre-assigned prescribed performance function. Number of neurons of the neural network and their types are defined based on try and error  which can be considered as a main drawback of this method.
   
  The work in this section is mainly based on reproducing \cite{bechlioulis_robust_2008} to catch the idea of prescribed performance function and to evaluate the function with adaptive control.

\section{Problem Formulation and Preliminaries}
For compactness and easy reading of the chapter, this section presents the concept of prescribed performance (for more details the reader is invited to consult \cite{bechlioulis_robust_2008}).\\
Consider the general case of nonlinear affine system as follows
       \begin{equation}
         \label{eq:ch3t}
         \begin{aligned}
           & x_1^{(n_1)}=f_1(x)+g_{11}(x)u_1+ \cdots +g_{1m}(x)u_m\\
           \vdots\\
           & x_m^{(n_m)}=f_m(x)+g_{m1}(x)u_1+ \cdots +g_{mm}(x)u_m
         \end{aligned}
       \end{equation}
 which can be adequately written in the form:\\
    \begin{equation*}
      x^{(n)}=f(x)+G(x)u\\
    \end{equation*}
   where\\
    \begin{equation*}
     x^{(n)}=
     \begin{bmatrix}
      x_1^{(n_1)} & \cdots & x_m^{(n_m)}
     \end{bmatrix}^{\top}
    \end{equation*}
    \begin{equation*}
     f(x)=
     \begin{bmatrix}
      f_1(x) & \cdots & f_m(x)
     \end{bmatrix}^{\top}
    \end{equation*}
    \begin{equation*}
     G(x)=
     \begin{bmatrix}
      g_{11}(x) & \cdots & g_{1m}(x)\\
      \vdots & \ddots & \vdots\\
      g_{m1}(x) & \cdots & g_{mm}(x)\\
     \end{bmatrix}
    \end{equation*}
The use of Prescribed performance with robust adaptive control demand considering four basic assumptions.
\begin{assumption}
 The matrix $\frac{G(x)G^{\top}(x)}{2)}$ has to be known with either uniformly positive definite or uniformly negative definite for all $x \in \Omega_{x}$  where $\Omega_{x} \subseteq \mathbb{R}^{n} $  is a compact set to guarantee system controllability.
 \label{assumption_51}
\end{assumption}
\begin{equation}
\underline{\sigma}(\frac{G(x)G^{\top}(x)}{2}) \geq g^* > 0 \hspace{10pt} \forall x \in \Omega_{x}
\end{equation}

where $\underline{\sigma}(W)$  is the smallest singular value of the matrix W and $g^{*}$ represents its lower bound. In addition, if $G(x)$ satisfies Assumption \ref{assumption_51} then system  \label{eq:ch3t} is uniformly strongly controllable \cite{bechlioulis_robust_2008}.
\begin{assumption}
The desired trajectories are known bounded functions of time with bounded known derivatives.
\end{assumption}
\begin{assumption}
The system states are available for measurement.
\end{assumption}
\begin{assumption}
 The functions $f_{i}(x)$ and $g_{ij}(x),i,j=1,\cdots,m$ are continuous but otherwise completely unknown.
\end{assumption}
Prescribed performance can be defined as the effort of tracking a generic error\\ $e\left(t\right)=[e_1\left(t\right), e_2\left(t\right),...,e_m\left(t\right)]\in \mathbb{R}^n $ such that each element of $e\left(t\right)$ evolves within PPB in a form of decaying functions of time that define the range of the residual error, the speed of convergence to the residual set, and the allowable overshoot or undershoot. In addition, prescribed performance with robust adaptive control was mainly developed to provide an adequate command signal for smooth tracking and solve the problem of accurate computation of the transient and steady state error bounds by guarantying uniform ultimate boundedness property of the error.  \\
A smooth function $\rho_i\left(t\right): \mathbb{R}_{+} \to \mathbb{R}_{+}$  is defined as a performance function associated with error component $e_i\left(t\right)$, $i=1, ..., m$, if $\rho_i\left(t\right)$ is positive, decreasing and $\lim\limits_{t \to \infty}\rho_i\left(t\right)=\rho_{i_\infty}>0$.

\subsection{Performance Functions}
A smooth function $\rho_i\left(t\right): \mathbb{R}_{+} \to \mathbb{R}_{+}$  is defined as a performance function associated with error component $e_i\left(t\right)$, $i=1, ..., m$, if $\rho_i\left(t\right)$ is positive, decreasing and $\lim\limits_{t \to \infty}\rho_i\left(t\right)=\rho_{i_\infty}>0$.
A possible choice of such function can be   
\begin{equation}
  \label{eq:ch3PPF}
  \rho_i\left(t\right)=(\rho_{i_0} - \rho_{i_\infty})\exp^{-\ell_i\, t}+\rho_{i_\infty}
\end{equation}
where $\rho_{i_0},\rho_{i_\infty}$ and $\ell_i$ are appropriately defined positive constants. The control objective is to guarantee that
\begin{equation}
  \label{eq:ch3onePPF}
  -\delta_i\rho_i\left(t\right)<e_i\left(t\right)<\rho_i\left(t\right),\hspace{10pt} if \: e_i(0)>0
\end{equation}
\begin{equation}
  \label{eq:ch3twoPPF}
  -\rho_i\left(t\right)<e_i\left(t\right)<\delta_i\rho_i\left(t\right),\hspace{10pt} if \: e_i(0)<0
\end{equation}
for all $ t \geq 0 $ and $ 0 \leq \delta_i \leq 1 $, and $i=1, ...,m$. Figure~\ref{fig:ch3PPF} illustrates the prescribed performance function and tracking error evolving from a large to a small set as per equations \eqref{eq:ch3onePPF} and \eqref{eq:ch3twoPPF}. 
\begin{figure}[h!]
 \centering
 \includegraphics[scale=0.5]{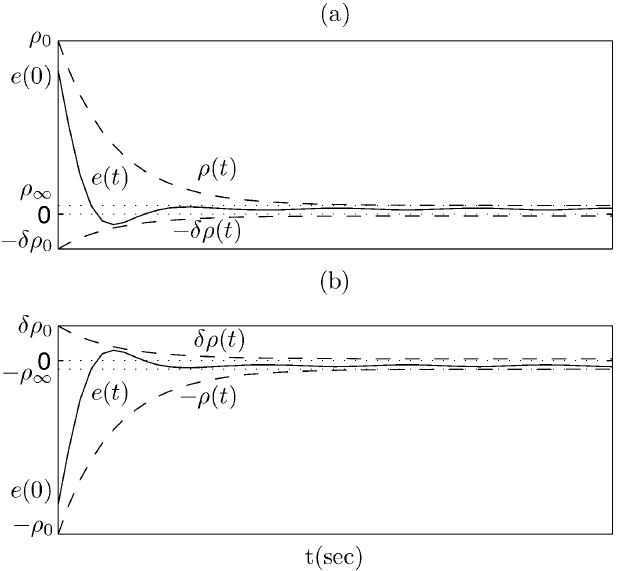}
 \caption{Graphical illustration of PPF for the tracking error behavior (a) graphical illustration of \eqref{eq:ch3onePPF}; (b) graphical illustration of \eqref{eq:ch3twoPPF}.}
 \label{fig:ch3PPF}
\end{figure}

\subsubsection{Error Transformation}
To implement the prescribing performance, one needs to solve a constrained control problem. To avoid such difficulty,
the following error transformation is used 
\begin{equation}
 \epsilon_i =T_i(\frac{e_i\left(t\right)}{\rho_i\left(t\right)})
\end{equation}
or equivalently,
\begin{equation}
  e_i\left(t\right)=\rho_i\left(t\right)S(\epsilon_i)
\end{equation}
where $\epsilon_i,\:i=1,2,...,m$ is the transformed error and $S_i(.)$ and $T_i^{-1}(.)$ are two smooth functions such that $S_i(.)=T_i^{-1}(.)$ and $S_i(.)$ satisfy the following properties:
\begin{enumerate}
  \item $S_i(\epsilon_i)$ is smooth and strictly increasing.
  \item $-\delta_i<S_i(\epsilon_i)<1,\hspace{10pt} if \: e_i(0)>0$\\
    $-1<S_i(\epsilon_i)<\delta_i,\hspace{10pt} if \: e_i(0)<0$
  \item
      $\left.
       \begin{aligned}
         lim_{\epsilon_i \rightarrow -\infty}S_i(\epsilon_i)=-\delta_i\\
         lim_{\epsilon_i \rightarrow +\infty}S_i(\epsilon_i)=1,
       \end{aligned}
       \right\}
       \qquad if \hspace{10pt} e_i(0)\geq 0$\\
      $\left.
       \begin{aligned}
         lim_{\epsilon_i \rightarrow -\infty}S_i(\epsilon_i)=-1\\
         lim_{\epsilon_i \rightarrow +\infty}S_i(\epsilon_i)=\delta_i,
       \end{aligned}
       \right\}
       \qquad if \hspace{10pt} e_i(0)< 0$
   \end{enumerate}
where\\
   \begin{equation}
    \label{eq:ch3PPFS}
     S(\epsilon)=
        \left\{
        \begin{aligned}
        \frac{\bar{\delta}e^{\epsilon}-\underline{\delta}e^{-\epsilon}}{e^{\epsilon}+e^{-\epsilon}},& \quad \bar{\delta}=1 \hspace{5pt} and \hspace{5pt} \underline{\delta}=0 \hspace{5pt} if \hspace{5pt} e(0)\geq 0\\
        \frac{\bar{\delta}e^{\epsilon}-\underline{\delta}e^{-\epsilon}}{e^{\epsilon}+e^{-\epsilon}},& \quad \underline{\delta}=1 \hspace{5pt} and \hspace{5pt} \bar{\delta}=0 \hspace{5pt} if \hspace{5pt} e(0)\geq 0\\
        \end{aligned}
        \right.
    \end{equation}
One should note that the overshoot/undershoot in equation \eqref{eq:ch3PPFS} is assumed to be zero.

where the overshot in equation \eqref{eq:ch3PPFS} assumed to be zero.
To continue, an error transformation that modulates $e_{i}\left(t\right)$ with respect to the corresponding performance bounds has to be defined. More specifically, we define the following transformed errors:
\begin{equation}
  \epsilon=S^{-1}\Big(\frac{\rho\left(t\right)}{e\left(t\right)}\Big)
\end{equation}
Next a metric error $E\left(t\right)$ will be defined to describe the system dynamics in a new form of system error.
\begin{equation}
  E_{i}\left(t\right)=(\frac{d}{dt}+\lambda_{i})^{n-1}\epsilon_i
\end{equation}
\begin{equation}
  \label{eq:ch3Edot}
  \dot{E}\left(t\right)= V+R\dot{x}
\end{equation}
where $\dot{E}\left(t\right)=\begin{bmatrix} E_1 & \cdots & E_n \end{bmatrix}^{\top}$ and $V=\begin{bmatrix} v_1 & \cdots & v_n \end{bmatrix}^{\top}$.
\begin{equation}
 \label{eq:ch3RR}
 R = 
    \begin{bmatrix}
     \frac{1}{2\rho_1\left(t\right)}(\partial S_1^{-1}/\partial(\frac{\rho_1\left(t\right)}{e_1\left(t\right)})) & \cdots & 0\\
     \vdots & \ddots & \vdots\\
     0 & \cdots & \frac{1}{2\rho_n\left(t\right)}(\partial S_n^{-1}/\partial(\frac{\rho_n\left(t\right)}{e_n\left(t\right)}))\\
   \end{bmatrix}
\end{equation}
  Equations \eqref{eq:ch3Edot} and \eqref{eq:ch3RR} can be driven easily. All foregoing equations in addition to approximated nonlinear model will be implemented in order to define the required control signal. Online training of linearly parameterized neural network is mainly implemented to estimate the nonlinear model as presented in the following subsection.
\section{Neural Approximations}
  Neural network with linear parameterization can be expressed by the following relation
  \begin{equation}
    y=Z^{\top}(x)\theta
  \end{equation}
 
  where $y \in \mathbb{R}^{m}$ is the neural net output, $x \in \mathbb{R}^{n}$  is the neural input,   $\theta \in \mathbb{R}^{p}$ is a p-dimensional vector of synaptic weights and $Z(x)$ is a p-dimensional vector of regressor terms. Regressor terms may include high order functions of radial basis function \cite{buhmann_radial_2000}, sigmoid functions \cite{ito_approximation_1992} and shifted sigmoids \cite{cybenko_approximation_1989} are defined as high order neural network.
  
  The nonlinear system is considered to be unknown functions and may be represented by one layer neural network structure with linear in weights plus modeling error term $\forall x \in \Omega_{x}$ obtaining:
  \begin{equation}
    f(x)=Z_f^{\top}(x)\theta^{*}+\omega_f(x)
  \end{equation}
  \begin{equation}
    G(x)=
        \begin{bmatrix}
         Z_{G_{11}}^{\top}(x)\theta^{*} & \cdots & Z_{G_{1m}}^{\top}(x)\theta^{*}\\
         \vdots & \ddots & \vdots\\
         Z_{G_{m1}}^{\top}(x)\theta^{*} & \cdots & Z_{G_{mm}}^{\top}(x)\theta^{*}\\
       \end{bmatrix}
       +\omega_G(x)
  \end{equation}
  where $Z_f(x)=\begin{bmatrix} Z_{f_1}(x) & \cdots & Z_{f_m}(x) \end{bmatrix}$ , $Z_{f_i}(x)$ and $Z_{G_{i,j}}(x)\in \mathbb{R}^{p}$,$i,j=1,\cdots,m$ are selected basis functions and $\theta^* \in \mathbb{R}^{p}$ are constants but unknown parameters which are used to minimize the approximation errors $\omega_f(x), \omega_G(x) \forall x \in \Omega_{x}$. Number of regressor $p$ should be chosen appropriately and sufficiently large in order to have a suitable representation of the nonlinear system. The approximated errors $\omega_f(x), \omega_G(x)$ should satisfy the following conditions
  \begin{equation}
    ||\omega_f(x)|| \leq W_f, \forall x \in \Omega_{x}
  \end{equation}
  \begin{equation}
    ||\omega_G(x)|| \leq W_G, \forall x \in \Omega_{x}
  \end{equation}
  where $W_f > 0$ and $W_G > 0$ and they are constants.\\
  Furthermore, if we define:
  \begin{equation}
    f(x,\theta)=Z_f^{\top}(x)\theta
  \end{equation}
  \begin{equation}
    G(x,\theta)=
      \begin{bmatrix}
       Z_{G_{11}}^{\top}(x)\theta & \cdots & Z_{G_{1m}}^{\top}(x)\theta\\
       \vdots & \ddots & \vdots\\
       Z_{G_{m1}}^{\top}(x)\theta & \cdots & Z_{G_{mm}}^{\top}(x)\theta\\
     \end{bmatrix}
  \end{equation}
  Then, defining the control law require the following variables
  \begin{equation}
    F_G(x,\theta)v=A_F(x,\theta)\theta
  \end{equation}
  where
  \begin{equation}
    A_F(x,\theta)=
      \begin{bmatrix}
       Z_{G_{11}}^{\top}(x) v_1+\cdots+ Z_{G_{1m}}^{\top}(x) v_m\\
       \vdots \\
       Z_{G_{m1}}^{\top}(x) v_1+\cdots+ Z_{G_{mm}}^{\top}(x) v_m\\
     \end{bmatrix}
  \end{equation}
  
\section{Robust Adaptive Control Design}
  The control law may be formulated as following 
  \begin{equation}
    \label{eq:ch3uu}
    u = \nu_a - (\eta_{G_a}|\nu_a|^{2}+\eta_{G_b}|\nu_b|^{2})\frac{R^{\top}E}{sign(G(x))}
  \end{equation}
  \begin{equation}
   \label{eq:ch3nua}
    \nu_a(x,\hat{\theta})=-\frac{Adj(F_G(x,\hat{\theta}))Det(F_G(x,\hat{\theta}))}{Det^2(F_G(x,\hat{\theta}))+\delta_d} \nu_b(x,\hat{\theta})
  \end{equation}
  \begin{equation}
   \label{eq:ch3nub}
    \nu_b(x,\hat{\theta})=F_f(x,\hat{\theta})+R^{-1}V+kR^{-1}E+n_fR^{\top}+E
  \end{equation}
  For $\eta_{G_a},\eta_{G_b},n_f,k$ and $\delta_d$ are positive constants and $F_f(x,\hat{\theta})$ and $F_G(x,\hat{\theta})$ are the approximations of $f(x)$ and $G(x)$. $\delta_d$ is necessary to make equation \eqref{eq:ch3nua} free of singularities. In order to validate equations \eqref{eq:ch3uu},\eqref{eq:ch3nua} and \eqref{eq:ch3nub}, let's formulate Lyapunov candidate function as
  \begin{equation}
    L = \frac{1}{2}E^{\top}E + \frac{1}{2} \tilde{\theta}^{\top}\Gamma^{-1}\tilde{\theta}
  \end{equation}
  \begin{equation*}
    \dot{L} = \frac{1}{2}\dot{E}^{\top}E + \frac{1}{2}E^{\top}\dot{E} + \frac{1}{2} \dot{\tilde{\theta}}^{\top}\Gamma^{-1}\tilde{\theta} + \frac{1}{2} \tilde{\theta}^{\top}\Gamma^{-1}\dot{\tilde{\theta}}
  \end{equation*}
  \begin{equation*}
    \dot{L} = \frac{1}{2}(V+R(f(x)+G(x)u))^{\top}E + \frac{1}{2}E^{\top}(V+R(f(x)+G(x)u)) + \frac{1}{2} \dot{\tilde{\theta}}^{\top}\Gamma^{-1}\tilde{\theta} + \frac{1}{2} \tilde{\theta}^{\top}\Gamma^{-1}\dot{\tilde{\theta}}
  \end{equation*}
  And after some manipulations, next equation will be chosen to validate global stability of the control law
  \begin{equation}
   \dot{\hat{\theta}} = \Gamma\big((Z_f^{\top}(x)+A_F(x,\nu_a))^{\top}R^{\top}E - \sigma(\theta-\theta_0)\big)
  \end{equation}
  Where $\sigma > 0$  and  $\theta_0$ a parameter vector used to incorporate a good guess of  $\theta$. Finally, $\dot{L}$ will be equivalent to
    \begin{equation*}
         \begin{split}
          \dot{L} \leq & -k|E^2| - \frac{\sigma}{2}|\tilde{\theta}|^2 - \eta_f|R^{\top}E|^2 + |R^{\top}E|W_f + |R^{\top}E|^2|\nu_a|W_G \\
          & - \eta_{G_a}g^*|\nu_a|^2|R^{\top}E|^2 + |R^{\top}E||\nu_b| - \eta_{G_b}g^*|\nu_b|^2|R^{\top}E|^2 + \frac{\sigma}{2}|\theta^*-\theta|^2
         \end{split}
    \end{equation*}
  Finally we will have
    \begin{equation*}
         \dot{L} \leq -k|E^2| - \frac{\sigma}{2}|\tilde{\theta}|^2 + \frac{W_f^2}{4\eta_f} + \frac{W_G^2}{4\eta_{G_a}g^*} + \frac{1}{4\eta_{G_b}g^*} + \frac{\sigma}{2}|\theta^*-\theta|^2
    \end{equation*}
    and if we choose $ d = \frac{W_f^2}{4\eta_f} + \frac{W_G^2}{4\eta_{G_a}g^*} + \frac{1}{4\eta_{G_b}g^*} + \frac{\sigma}{2}|\theta^*-\theta|^2 $, then the value of $d$ will be reflected on the value of $E$ or/and $\tilde{\theta}$.For more details look \cite{bechlioulis_robust_2008}.

\section{Problem Simulation and Results}
Consider equations of motion of 2 DOF planner robot in example 3.3.1, the nonlinear plant assumed to be completely unknown. Single layer neural network with 30 neurons sigmoid basis function were used to estimate the system nonlinearities $-M^{-1}(q)( C(\dot{q},q)\dot{q} + G_0(q) )$ and  $M^{-1}(q)$. The parameters of the sigmoid basis function $\zeta_j(x)=1/(1+e^{-\omega_j^{\top} -b_j})$ with $\omega_j \in \mathbb{R}^{4}$, $b_j \in \mathbb{R}^{4}$, $j=1,2,\cdots,30$ were chosen by off-line training try and error on the simulation then kept constant throughout the simulation. $\theta_0$ is a vector represents the good guess of the initial conditions of the parameter estimates and was taken to be a zero vector referring to completely unknown nonlinear dynamics. \\
  The robot assumed to start initially from the origin while the desired trajectory for both angles were chosen to be
  \begin{equation*}
     q_d = \begin{bmatrix} 0.5cos(0.7t) & -0.6cos(0.65t) \end{bmatrix}^{\top}
  \end{equation*}
  Prescribed performance function was chosen as
  \begin{equation*}
    \rho_i\left(t\right)=(\rho_{i0} - \rho_{i\infty})e^{-l_it}+\rho_{i\infty}, \hspace{5pt} i = 1,2
  \end{equation*}
  Prescribed performance parameters are demonstrated in table~\ref{tab:ch3_PPparm} and parameters of controller are defined table~\ref{tab:ch3_conparm}
  \begin{table}[h]
    \caption {Prescribed performance function parameters}
    \label{tab:ch3_PPparm} 
    \begin{center}
        \begin{tabular}{ l | l | l | l | l | l }
        \hline\hline
        $\rho_{10}$ & $\rho_{1\infty}$  & $l_1$ & $\rho_{20}$ & $\rho_{2\infty}$  & $l_2$  \\ \hline
        1.1 & 0.005 & 2.0 & 1.1 & 0.005 & 2.0  \\ \hline \hline
      \end{tabular}
    \end{center}
  \end{table}

  \begin{table}[h]
    \caption {Adaptive PPF Controller parameters}
    \label{tab:ch3_conparm} 
    \begin{center}
      \begin{tabular}{ l | l | l | l | l | l | l | l | l }
      \hline\hline
      $k$ & $n_f$  & $\eta_{G_a}$ & $\eta_{G_b}$ & $\delta_d$ & $\sigma$ & $\Gamma$ & $\lambda_1$ & $\lambda_2$  \\ \hline
      0.5 & 0.2 & 0.2 & 0.2 & 0.1 & 7.5 & 0.1I & 0.75 & 0.75  \\ \hline \hline
      \end{tabular}
    \end{center}
  \end{table}

   In figure~\ref{fig:PPFout}, angular positions of both actual and desired trajectory had verified the control efficacy. Figure~\ref{fig:PPFcont} demonstrates the smoothness of the control signal along the trajectory. Figure~\ref{fig:PPFTra} presents bounds of the prescribed performance function and verify that the error of each joint is bounded within a large set and ended within a small preassigned set. Finally, transformed errors both joints are demonstrated in figure ~\ref{fig:PPFTra}.
  \begin{figure}[h!]
   \centering
   \includegraphics[scale=0.45]{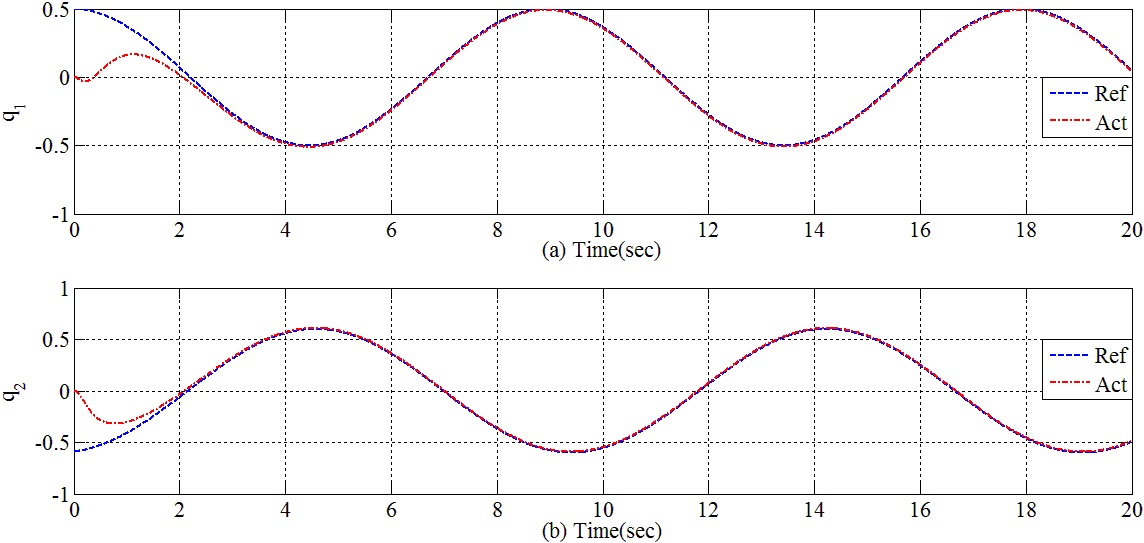}
   \caption{Output response of the robust adaptive control with PPF for $q_1$ and $q_2$ versus desired trajectory $q_{d1}$ and $q_{d2}$}
   \label{fig:PPFout}
  \end{figure}

  \begin{figure}[h!]
    \centering
    \includegraphics[scale=0.45]{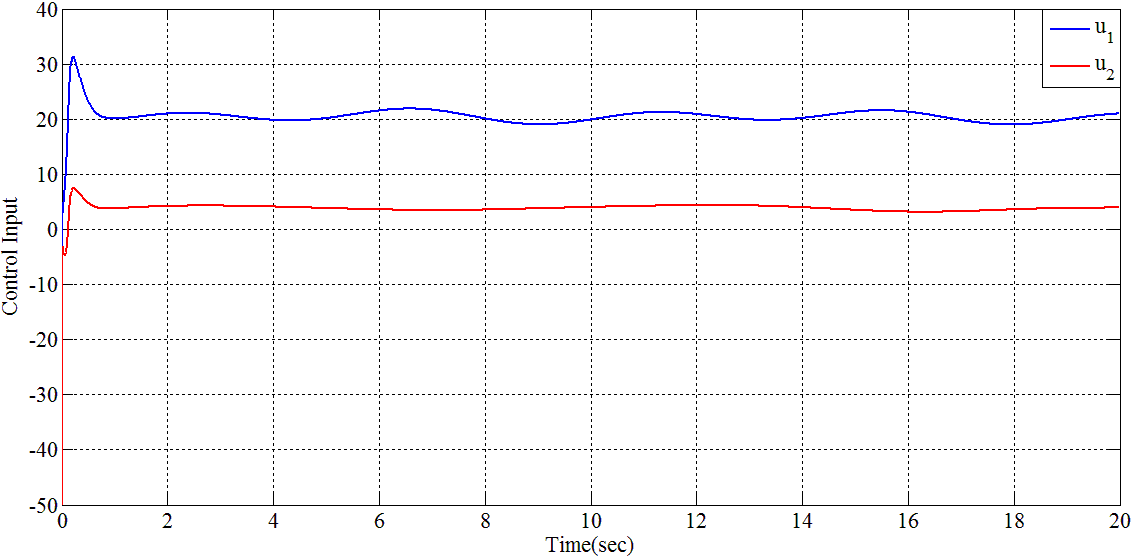}
    \caption{Control input provided by robust adaptive control with PPF where $u_1$ is $\tau_1$ and $u_2$ is $\tau_2$.}
    \label{fig:PPFcont}
  \end{figure}

\begin{figure}[h!]
  \centering
  \includegraphics[scale=0.45]{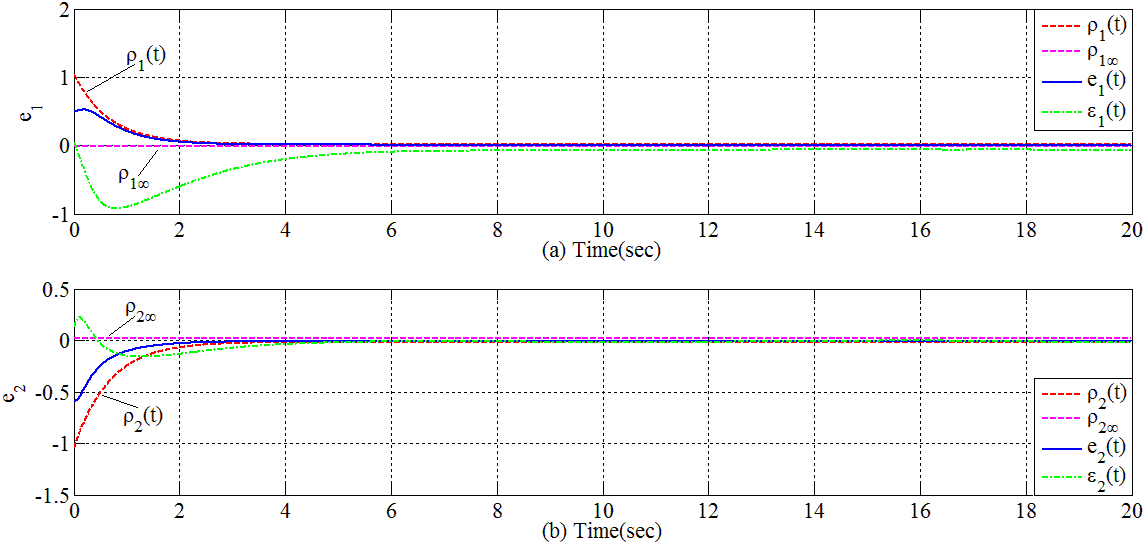}
  \caption{Prescribed error bounds between $\rho_0$ and $\rho_{\infty}$ and $\epsilon$ for both joints (a)$q_1$ and (b)$q_2$.}
  \label{fig:PPFTra}
\end{figure}

\section{Conclusion}
This chapter illustrated the significant role of prescribed performance function with robust adaptive control. The main idea of Prescribed performance has been gained and the controller showed smoothness in the control signal and impressive tracking performance. In a subsequent chapter, new controller stands on PPF will be developed relies on the result of this chapter.

\clearpage

\newpage


\chapter{ROBUST MRAC WITH PPF FOR NONLINEAR MIMO SYSTEMS}

\section{Introduction}           
In this work, we are motivated by the limitations of the studies presented in the literature and mentioned in chapter 3 and 5 to propose a robust MRAC with PPF. We show that the robust stabilization of the transformed error guaranties the stability and convergence of the constrained tracking error within the set of time varying constraints representing the performance limits. Simulation results benchmark the performance of the proposed approach with \Lone adaptive control and neuro-adaptive control with prescribed performance. The rest of the chapter is organized as follows. In section two, the problem formulation with important remark are presented. The design and analysis of the proposed robust MRAC-PPF, which represents the main contribution, is presented in section three. In section four, simulation results verify the effectiveness of the proposed control and show that the MRAC-PPF considerably improves the transient performance when compared to \Lone adaptive control and Neuro-Adaptive controller with PPF. We conclude the chapter in section five.

\section{Problem Formulation}
We consider the following uncertain system defined by
  \begin{equation}
  \label{system1}
    \begin{aligned}
      &\dot{x} = A x\left(t\right)+B u\left(t\right) + \theta^{\top}x\left(t\right) + \Delta f(x,u,t) + d\left(t\right),\hspace{10pt} x(0)=x_0.\\
      & y\left(t\right)=Cx\left(t\right).
    \end{aligned}
  \end{equation}
  where $ \Delta f $ is an unknown uncertainty and $d\left(t\right)$ is the system unknown but bounded disturbance. And Let the desired dynamics be defined as following
      \begin{equation}
      \label{refmodel}
        \begin{aligned}
          &\dot{x_m}\left(t\right) = A_m x\left(t\right)+B_mr\left(t\right),\hspace{10pt} B_m=Bk_g.
        \end{aligned}
      \end{equation}
 where $A_m$ is a Hurwitz matrix, and both pairs (A,B) and (Am,B) are controllable. Consider $ u\left(t\right) = u_m\left(t\right)+u_{ad}\left(t\right) $ where $ u_m\left(t\right)=-k_m x\left(t\right)$ and $k_m$ is a state feedback gain such that $A_m=A\,-\,k_m x\left(t\right)$.
    \begin{equation}
      \begin{aligned}
        &\dot{x}\left(t\right) = A_m x\left(t\right)+Bu_{ad}\left(t\right) + \theta^{\top}\,x\left(t\right) + \Delta f(x,u,t) + d\left(t\right)
      \end{aligned}
    \end{equation}  
  \begin{remark} $B$ is not necessary a square matrix but satisfies
  \begin{equation}
  \underline{\sigma}(\frac{B\,B^{\top}}{2}) \geq g^* > 0 
  \end{equation}
     $A\in \mathbb{R}^{n \times n} $, $B\in \mathbb{R}^{n \times m} $, $x\left(t\right)\in \mathbb{R}^n $, $u\left(t\right)$ and $r\left(t\right)$ are $ \in \mathbb{R}^m $, $\Delta f(x,u,t)\in \mathbb{R}^n $ and  $d\left(t\right)\in \mathbb{R}^n $.
 \end{remark}
 Let the error be $e=x-x_m$,  then    
    \begin{equation}
    \label{error1}
      \begin{aligned}
        &\dot{e} =\dot{x}-\dot{x}_m= A_m e+B(u_{ad}-k_g\,r)+ \theta^{\top}\,x\left(t\right) + \Delta f(x,u,t) + d\left(t\right)\\
      \end{aligned}
    \end{equation}  

\section{Controller Structure}
 Let 
   \begin{equation}
     e\left(t\right)=\rho\left(t\right)S(\epsilon)
   \end{equation}
   \begin{equation}
     \rho\left(t\right)=(\rho_0-\rho_\infty)e^{-lt}+\rho_\infty
   \end{equation}
the transformed error is then
      \begin{equation}
         \epsilon=S^{-1}(\frac{\rho\left(t\right)}{e\left(t\right)})
      \end{equation}
where
      \begin{equation}
        S(\epsilon)=
          \left\{
          \begin{aligned}
          \frac{\bar{\delta}e^{\epsilon}-\underline{\delta}e^{-\epsilon}}{e^{\epsilon}+e^{-\epsilon}},& \quad \bar{\delta}=1 and \underline{\delta}=0 if e(0)\geq 0\\
          \frac{\bar{\delta}e^{\epsilon}-\underline{\delta}e^{-\epsilon}}{e^{\epsilon}+e^{-\epsilon}},& \quad \underline{\delta}=1 and \bar{\delta}=0 if e(0)\geq 0\\
          \end{aligned}
          \right.
      \end{equation}  
 and
     \begin{equation}
        \epsilon=S^{-1}(\frac{\rho\left(t\right)}{e\left(t\right)})\,= \frac{1}{2}ln(\underline{\delta}+e\left(t\right)/\rho\left(t\right))-\frac{1}{2}ln(\bar{\delta}-e\left(t\right)/\rho\left(t\right))
     \end{equation}
Let
    \begin{equation}
        \frac{1}{2\rho\left(t\right)}\big(\partial S^{-1}(\epsilon)/\epsilon)= \frac{1}{2\rho\left(t\right)}\Big(\frac{1}{\underline{\delta}+e\left(t\right)/\rho\left(t\right)} - \frac{1}{e\left(t\right)/\rho\left(t\right) - \bar{\delta}}\Big)
    \end{equation}
which can be written in matrix form as    
    \begin{equation*}
      \Gamma =
      \begin{pmatrix}
       \frac{1}{2\rho_1\left(t\right)}\big(\partial S^{-1}(\epsilon_1)/\epsilon_1) & \cdots & 0\\
      \vdots & \ddots & \vdots\\
       0 & \cdots & \frac{1}{2\rho_n\left(t\right)}\big(\partial S^{-1}(\epsilon_n)/\epsilon_n)\\
      \end{pmatrix}
    \end{equation*} 
 Let   
    \begin{equation*}
      \Phi =-\Gamma
      \begin{pmatrix}
       e_1\left(t\right)/\rho_1\left(t\right) & \cdots & 0\\
      \vdots & \ddots & \vdots\\
       0 & \cdots & e_n\left(t\right)/\rho_n\left(t\right)\\
      \end{pmatrix}
    \end{equation*} 
    \begin{equation*}
      \epsilon =
      \begin{pmatrix}
       \epsilon_1\\
       \vdots\\
       \epsilon_n\\
      \end{pmatrix}
    \end{equation*} 
and
    \begin{equation*}
      \rho =
      \begin{pmatrix}
       \rho_1\\
       \vdots\\
       \rho_n\\
      \end{pmatrix}  
    \end{equation*}
 then
    \begin{equation}
       \label{transformed1}
      \dot{\epsilon} = \Gamma\dot{e}+\Phi\dot{\rho}=\Gamma \big(A_m e+B(u_ad -k_gr)+ \theta^{\top}x\left(t\right) + \Delta f(x,u,t) + d\left(t\right)\big)+\Phi\dot{\rho}
    \end{equation} 
 Let
        \begin{equation}
         \label{gamma}
          \gamma(x) = \theta^{\top}\,x\left(t\right) + \Delta f(x,u,t) + d\left(t\right)
        \end{equation}
Assume
    \begin{equation}
      \label{gamma2}
           \gamma(x) = \theta^{\top}\,x\left(t\right) + \sigma^{\top} \,\psi(x,u)+\alpha(x,u) 
        \end{equation}
where $\alpha(x,u)$ represents all the unknown nonlinear in parameters terms such that $\alpha(x) \leq \bar{\alpha}_i$.
Let         
        \begin{equation}
     \label{v1}
           V = \Gamma \big(A_m e-B k_gr)\big)+\Phi\dot{\rho}
        \end{equation}
and define    
  \begin{equation}
       V_n=\Gamma^{-1}\,V
       \label{vn} 
    \end{equation}  
Consider
\begin{equation}
      \label{gammaa2}
           \hat{\gamma}(x) = \hat{\theta}^{\top}\,x\left(t\right) + \hat{\sigma}^{\top} \,\psi(x,u) 
        \end{equation}
where $\hat{(.)}$ stands for the estimate. Then
    \begin{equation*}
       \gamma(x) - \hat{\gamma}(x,\hat{\theta},\hat{\sigma}) = \tilde{\theta}^{\top}x\left(t\right) + \tilde{\sigma}^{\top}\psi(x) + \alpha(x).   
    \end{equation*} 
    \begin{equation*}
       \tilde{\theta} = \hat{\theta} - \theta, \tilde{\sigma}=\hat{\sigma}-\sigma      
    \end{equation*} 
    \begin{equation*}
      \epsilon^{\top}\dot{\epsilon} = \epsilon^{\top}\Big(\Gamma \big( - \hat{\gamma}(x,\hat{\theta},\hat{\sigma})-V_n \big)+V\Big)
    \end{equation*}
    \begin{equation}
      \epsilon^{\top}\dot{\epsilon} = -\sum_{i=1}^{n}\epsilon_i \Gamma_{i,i}\tilde{\theta}_{:,i}^{\top} x\left(t\right)-\sum_{i=1}^{n}\epsilon_i \Gamma_{i,i}\tilde{\sigma}_{:,i}^{\top} \psi(x) +\sum_{i=1}^{n}\epsilon_i \Gamma_{i,i}\bar{\alpha}_i
    \end{equation} 
  It is important to notice that
    \begin{equation*}
      \tilde{\theta}_{:,i}^{\top}\hat{\theta}_{:,i} = \frac{1}{2} \tilde{\theta}_{:,i}^{\top}\tilde{\theta}_{:,i}+\frac{1}{2}\big(\hat{\theta}_{:,i}-\theta_{:,i}\big)^{\top}\big(\hat{\theta}_{:,i}+\theta_{:,i}\big)\geq \frac{1}{2} \tilde{\theta}_{:,i}^{\top}\tilde{\theta}_{:,i} -  \frac{1}{2} \theta_{:,i}^{\top}\theta_{:,i}
    \end{equation*}
    \begin{equation*}
      -\tilde{\theta}_{:,i}^{\top}\hat{\theta}_{:,i} \leq -\frac{1}{2} \tilde{\theta}_{:,i}^{\top}\tilde{\theta}_{:,i} +  \frac{1}{2} \theta_{:,i}^{\top}\theta_{:,i}
    \end{equation*}
The control signal can be selected as
\begin{equation}
\label{uad}
      u_{ad}\left(t\right) = B^{-1}\big(-\hat{\theta}^{\top}x\left(t\right) - \hat{\sigma}^{\top} \, \psi(x)-V_n \big)+u_{r}\left(t\right)
 \end{equation} 
where $B^{-1}$ can be replaced by its Moore$-$Penrose inverse when it is not square owing to Assumption. Let the adaption rules for $\hat{\theta}$ and  $\hat{\sigma}$ be defined as follows respectively
\begin{equation}
      \label{thetahat}
      \dot{\hat{\theta}}_{:,i}= -\gamma_{1i}\epsilon_i \Gamma_{i,i} x\left(t\right)
    \end{equation}
    \begin{equation}
     \label{sigmahat}
      \dot{\hat{\sigma}}_{:,i}= \int_0^\infty \Gamma_{i,i}\big(-\gamma_{2i}|\epsilon_i|\upsilon_i\hat{\sigma}_{:,i}+\gamma_{2i}\epsilon_i\psi(x)\big)\mathrm{d}\tau -\beta_i\delta_i
    \end{equation}   
    \begin{equation}
  \label{delta}
      \delta_i = \gamma_{2i}|\epsilon_i|\upsilon_i\hat{\sigma}_{:,i}+\gamma_{2i}\epsilon_i\psi(x)
    \end{equation}
    \begin{equation}
     \label{epsilonbar}
      \hat{\bar{\alpha}} \geq \bar{\alpha}_i + \frac{1}{2}||\sigma_{:,i}||^2\Gamma_{i,i}\upsilon_i
    \end{equation}
 and the robustifying term
       \begin{equation}
       \label{robust1}
         u_{r}=\left[u_{ri}\right] = \left[-sign(\epsilon)_i \cdot \hat{\bar{\alpha}}\right]
       \end{equation}
We are now ready to announce the following theorem. 
\begin{theorem}
Under Assumption 1 with the 
prescribed performance defined by ( \ref{eq:ch3PPF}),  the MRAC of 
System (\ref{system1}) with reference model (\ref{refmodel}) having the error dynamic (\ref{error1}) and the transformed error dynamic (\ref{transformed1}), the control input defined by (\ref{uad}), equations (\ref{gammaa2})-(\ref{vn}), and the adaption rule (\ref{thetahat})-(\ref{delta}) and the robustifying term (\ref{epsilonbar})-(\ref{robust1}), forces the transformed error to asymptotically reach zero and therefore the tracking error to satisfy the prescribed performance.
\end{theorem}

\section{Stability Analysis}
The proof is similar to the one in \cite{sun_fuzzy_2014}. We adapted it to our case. Let us consider the Lyapunov candidate
  The Lyapunov candidate may be chosen as
      \begin{equation}
        W=W_1+W_2
      \end{equation}
      \begin{equation*}
        W_1=\epsilon^{\top}\epsilon
      \end{equation*}      
      \begin{equation*}
        W_2=\sum_{i=1}^{n}\frac{1}{2\gamma_{1i}}\tilde{\theta}_{:,i}^{\top}\tilde{\theta}_{:,i} + \sum_{i=1}^{n}\frac{1}{2\gamma_{2i}}\big(\tilde{\sigma}_{:,i}+\beta_i\delta_i\big)^{\top}\big(\tilde{\sigma}_{:,i}+\beta_i\delta_i\big)
      \end{equation*}
     \begin{equation*}
            \dot{W} = \dot{W_1}+\dot{W_2} 
      \end{equation*}
      \begin{equation*}
        \dot{W_1} = \dot{\epsilon}^{\top}\epsilon + \epsilon^{\top}\dot{\epsilon}
      \end{equation*}
      \begin{eqnarray}
        \dot{W_2}=\sum_{i=1}^{n}\frac{1}{2\gamma_{1i}}\tilde{\theta}_{:,i}^{\top}\dot{\hat{\theta}}_{:,i} + \sum_{i=1}^{n}\frac{1}{\gamma_{2i}}\big(\tilde{\sigma}_{:,i}+ \beta_i\delta_i\big)^{\top}\big(\dot{\hat{\sigma}}_{:,i}+\beta_i\dot{\delta}_i\big)
      \end{eqnarray}
      \begin{eqnarray}
        \dot{W_1}&\leq -\sum_{i=1}^{n}\epsilon_i \Gamma_{i,i}\tilde{\theta}_{:,i}^{\top} x\left(t\right)-\sum_{i=1}^{n}\epsilon_i \Gamma_{i,i}\tilde{\sigma}_{:,i}^{\top} \psi(x) +\sum_{i=1}^{n}|\epsilon_i| \Gamma_{i,i}\bar{\alpha}_i \nonumber \\&- \sum_{i=1}^{n}\epsilon_i \Gamma_{i,i}K_{i,i}\epsilon_i +  \sum_{i=1}^{n}\epsilon_i \Gamma_{i,i}u_{ri}
      \end{eqnarray}
      \begin{equation}
      \begin{split}
        \sum_{i=1}^{n}\frac{1}{\gamma_{2i}}& \big(\tilde{\sigma}_{:,i}+ \beta_i\delta_i\big)^{\top}\big(\dot{\hat{\sigma}}_{:,i}+\beta_i\dot{\delta}_i\big) \leq - \sum_{i=1}^{n}\frac{1}{2}||\tilde{\sigma}_{:,i}||^{2}|\epsilon_i|\Gamma_{i,i}\upsilon_i +\sum_{i=1}^{n}\frac{1}{2}||\sigma_{:,i}||^{2}|\epsilon_i|\Gamma_{i,i}\upsilon_i + \nonumber \\
         & \sum_{i=1}^{n}\frac{1}{2} \tilde{\sigma}_{:,i}\epsilon_i\Gamma_{i,i}\psi(x) - \sum_{i=1}^{n}\beta_{i}||\delta_i||
      \end{split}
      \end{equation} 
      \begin{eqnarray}
       \begin{split}
        \dot{W} & = \dot{W_1}+\dot{W_2} \leq -\sum_{i=1}^{n}\epsilon_i \Gamma_{i,i}\tilde{\theta}_{:,i}^{\top} x\left(t\right)-\sum_{i=1}^{n}\epsilon_i \Gamma_{i,i}\tilde{\sigma}_{:,i}^{\top} \psi(x) +\sum_{i=1}^{n}|\epsilon_i| \Gamma_{i,i}\bar{\alpha}_i - \sum_{i=1}^{n}\epsilon_i \Gamma_{i,i}K_{i,i}\epsilon_i \nonumber \\
        & +  \sum_{i=1}^{n}\epsilon_i \Gamma_{i,i}u_{ri} +\sum_{i=1}^{n}\frac{1}{2\gamma_{1i}}\tilde{\theta}_{:,i}^{\top}\dot{\hat{\theta}}_{:,i} + \sum_{i=1}^{n}\frac{1}{\gamma_{2i}}\big(\tilde{\sigma}_{:,i}+ \beta_i\delta_i\big)^{\top}\big(\dot{\hat{\sigma}}_{:,i}+\beta_i\dot{\delta}_i\big)\leq 0
       \end{split}
      \end{eqnarray}
  Then by choosing
      \begin{equation*}
        \dot{\hat{\sigma}}_{:,i}= \int_0^\infty \Gamma_{i,i}\big(-\gamma_{2i}|\epsilon_i|\upsilon_i\hat{\sigma}_{:,i}+\gamma_{2i}\epsilon_i\psi(x)\big)\mathrm{d}\tau -\beta_i\delta_i
      \end{equation*}   
      \begin{equation*}
        \delta_i = \gamma_{2i}|\epsilon_i|\upsilon_i\hat{\sigma}_{:,i}+\gamma_{2i}\epsilon_i\psi(x)
      \end{equation*}
 and
      \begin{equation*}
        u_{ri} = -sign(\epsilon)_i \cdot \hat{\bar{\alpha}}
      \end{equation*}
 one gets
      \begin{eqnarray}
       \begin{split}
        \dot{W} \leq &-\sum_{i=1}^{n}\epsilon_i \Gamma_{i,i}\tilde{\theta}_{:,i}^{\top} x\left(t\right)-\sum_{i=1}^{n}\epsilon_i \Gamma_{i,i}\tilde{\sigma}_{:,i}^{\top} \psi(x) +\sum_{i=1}^{n}|\epsilon_i| \Gamma_{i,i}\bar{\alpha}_i - \sum_{i=1}^{n}\epsilon_i \Gamma_{i,i}K_{i,i}\epsilon_i \\
        & -  \sum_{i=1}^{n}|\epsilon_i| \Gamma_{i,i} \hat{\bar{\alpha}} +\sum_{i=1}^{n}\frac{1}{2\gamma_{1i}}\tilde{\theta}_{:,i}^{\top}\dot{\hat{\theta}}_{:,i}- \sum_{i=1}^{n}\frac{1}{2}||\tilde{\sigma}_{:,i}||^{2}|\epsilon_i|\Gamma_{i,i}\upsilon_i\\
        & + \sum_{i=1}^{n}\frac{1}{2}\tilde{\sigma}_{:,i}\epsilon_i\Gamma_{i,i}\psi(x) - \sum_{i=1}^{n}\beta_{i}||\delta_i|| + \sum_{i=1}^{n}\frac{1}{2}||\sigma_{:,i}||^{2}|\epsilon_i|\Gamma_{i,i}\upsilon_i\leq 0
       \end{split}
      \end{eqnarray}
  Using the adaption rule   
      \begin{equation*}
        \dot{\hat{\theta}}_{:,i}= -\gamma_{1i}\epsilon_i \Gamma_{i,i} x\left(t\right)
      \end{equation*}    
and leads to
      \begin{equation}
       \begin{split}
        & +\sum_{i=1}^{n}|\epsilon_i| \Gamma_{i,i}\bar{\alpha}_i - \sum_{i=1}^{n}\epsilon_i \Gamma_{i,i}K_{i,i}\epsilon_i -  \sum_{i=1}^{n}|\epsilon_i| \Gamma_{i,i} \hat{\bar{\alpha}} \\
        &  - \sum_{i=1}^{n}\frac{1}{2}||\tilde{\sigma}_{:,i}||^{2}|\epsilon_i|\Gamma_{i,i}\upsilon_i + \sum_{i=1}^{n}\frac{1}{2}||\sigma_{:,i}||^{2}|\epsilon_i|\Gamma_{i,i}\upsilon_i - \sum_{i=1}^{n}\beta_{i}||\delta_i|| \leq 0
       \end{split}
      \end{equation}
      The following terms are negative  $ - \sum_{i=1}^{n}\epsilon_i \Gamma_{i,i}K_{i,i}\epsilon_i$, $- \sum_{i=1}^{n}\frac{1}{2}||\tilde{\sigma}_{:,i}||^{2}|\epsilon_i|\Gamma_{i,i}\upsilon_i $ and $- \sum_{i=1}^{n}\beta_{i}||\delta_i||$, Therefore one can  select
      \begin{equation}
        \sum_{i=1}^{n}|\epsilon_i| \Gamma_{i,i}\bar{\alpha}_i
          - \sum_{i=1}^{n}|\epsilon_i| \Gamma_{i,i} \hat{\bar{\alpha}} + \sum_{i=1}^{n}\frac{1}{2}||\sigma_{:,i}||^{2}|\epsilon_i|\Gamma_{i,i}\upsilon_i   \leq 0
      \end{equation}
      Which leads to
      \begin{equation}
        \sum_{i=1}^{n}|\epsilon_i| \Gamma_{i,i} \hat{\bar{\alpha}}
         \geq  \sum_{i=1}^{n}|\epsilon_i| \Gamma_{i,i}\bar{\alpha}_i  + \sum_{i=1}^{n}\frac{1}{2}||\sigma_{:,i}||^{2}|\epsilon_i|\Gamma_{i,i}\upsilon_i
      \end{equation}
     which is satisfied if $\hat{\bar{\alpha}}$ is selected as
       \begin{equation*}
          \hat{\bar{\alpha}}
          \geq  \big( \bar{\alpha}_i  + \frac{1}{2}||\sigma_{:,i}||^{2}\upsilon_i\big)
       \end{equation*}   
  In the next section, several simulation results to validate the approach ad assess its stability will be presented.    
\section{Simulation Examples}\label{simEx}
the performance of the proposed robust MRAC control design is demonstrated using two different cases. In each case, the control performance and its ability to guarantee the desired performance are benchmarked to first $\mathcal{L}_{1}$ adaptive controller and Neuro-adaptive controller.\\
   {\bf Example 6.5.1}
     \begin{equation*}
       \begin{aligned}
          &\dot{x} = A x\left(t\right)+ Bu\left(t\right) + \theta^{\top}x\left(t\right) + \Delta f + d\left(t\right),\hspace{10pt} x(0)=x_0.\\
          & y\left(t\right)=Cx\left(t\right).
       \end{aligned}
      \end{equation*}
     \begin{equation*}
       A =
       \begin{bmatrix}
        -36 & 36 & 0\\
        0 & 20 & 0\\
        0 & 0 & -3
       \end{bmatrix}
       ,B = 
       \begin{bmatrix}
        1 & 0 & 0\\
        0 & 1 & 0\\
        0 & 0 & 1
       \end{bmatrix}
       ,C = 
       \begin{bmatrix}
        1 & 0 & 0\\
        0 & 1 & 0\\
        0 & 0 & 1
       \end{bmatrix}
     \end{equation*} 
     \begin{equation*}
       \Delta f =
       \begin{bmatrix}
        x_3^2+0.2sin(x_1)\\
        -x_1x_3-0.2xos(x_3)x_1\\
        x_1x_2
       \end{bmatrix}
       ,d\left(t\right) = 
       \begin{bmatrix}
        1+sin\left(t\right)\\
        1.2+cos\left(t\right)\\
        sin\left(t\right)+cos\left(t\right)-1
       \end{bmatrix}
     \end{equation*}
     \begin{equation*}
       ,\theta\left(t\right) = 
       \begin{bmatrix}
        3sin(0.5t) & 2sin(0.4t)cos(0.3t) & 0.7sin(0.2t)\\
        0.9sin(0.2t) & 2.5sin(0.3t)+0.3cos\left(t\right) & sin(0.1t)\\
        0.5sin(0.13t) & 0.6cos(0.15t) & 1.5cos(0.7t)+1.6sin(0.3t)
       \end{bmatrix}
     \end{equation*}
   {\bf Example 6.5.2}
   \begin{equation*}
     \begin{aligned}
        &\dot{x} = A x\left(t\right)+ Bu\left(t\right) + \theta^{\top}x\left(t\right) + \Delta f + d\left(t\right),\hspace{10pt} x(0)=x_0.\\
        & y\left(t\right)=Cx\left(t\right).
     \end{aligned}
    \end{equation*}
    \begin{equation*}
      A =
      \begin{bmatrix}
       -36 & 36 & 0\\
       0 & 20 & 0\\
       0 & 0 & -3
      \end{bmatrix}
      ,B = 
      \begin{bmatrix}
       1 & 0 & 0\\
       0 & 1 & 0\\
       0 & 0 & 1
      \end{bmatrix}
      ,C = 
      \begin{bmatrix}
       1 & 0 & 0\\
       0 & 1 & 0\\
       0 & 0 & 1
      \end{bmatrix}
    \end{equation*} 
    \begin{equation*}
      \Delta f =
      \begin{bmatrix}
       x_3^2+0.2sin(x_1)-2.5u_3cos(u_1)\\
       -x_1x_3-0.2xos(x_3)x_1+0.7u_{3}^2\\
       x_1x_2
      \end{bmatrix}
      ,d\left(t\right) = 
      \begin{bmatrix}
       1+sin\left(t\right)\\
       1.2+cos\left(t\right)\\
       sin\left(t\right)+cos\left(t\right)-1
       \end{bmatrix}
     \end{equation*}
     \begin{equation*}
      ,\theta\left(t\right) = 
       \begin{bmatrix}
        3sin(0.5t) & 2sin(0.4t)cos(0.3t) & 0.7sin(0.2t)\\
        0.9sin(0.2t) & 2.5sin(0.3t)+0.3cos\left(t\right) & sin(0.1t)\\
        0.5sin(0.13t) & 0.6cos(0.15t) & 1.5cos(0.7t)+1.6sin(0.3t)
       \end{bmatrix}
     \end{equation*}      
 Desired poles are selected as $p=-70,-60\pm i$. \\
 {\bf Robust Adaptive Prescribed Performance Parameters Parameters} 
 $\rho_{i0}=2$,$\rho_{i\infty}=0.05$,$l_{i}=1.5$,and estimator parameters $\beta_{i}=2$,$\gamma_{1i}=50000$,$\gamma_{2i}=50000$,$\upsilon_{i}=0.05$ where ,$i=1,2,3$ and finally $\psi(x) = \begin{bmatrix} 2 & 2 & 2 \end{bmatrix}^{\top}$,$\hat{\bar{\epsilon}} = \begin{bmatrix} 10 & 10 & 10 \end{bmatrix}^{\top}$,$K=0.1diag(3)$\\
 Reference input assigned to be $r\left(t\right) = \begin{bmatrix} cos(0.75t) & cos(0.8t) & cos(0.7t) \end{bmatrix}^{\top}$,\\ 
 {\bf \Lone Adaptive Controller Parameters}
 $\theta_b \in [-5,5]$, $\Delta \in 20$, $\hat{\omega} \in [0.3,10]$, \\
 The parameters of the sigmoid basis function $\zeta_j(x)=1/(1+e^{-\omega_j^{\top} -b_j})$ with $\omega_j \in \mathbb{R}^{3}$, $b_j \in \mathbb{R}^{3}$, $j=1,2,\cdots,80$ were chosen by off-line training try and error on the simulation then kept constant throughout the simulation. $\theta_0$ is a vector represents the good guess of the initial conditions of the parameter estimates and was taken to be a zero vector referring to completely unknown nonlinear dynamics. Prescribed performance parameters are $\rho_{i0}=2$, $\rho_{i\infty}=0.05$, $l_{i}=1.5$.
   \begin{table}[htbp]
     \caption {Robust Neuro Adaptive Control with PPF parameters}
     \label{tab:ch3_conparm} 
     \begin{center}
       \begin{tabular}{ l | l | l | l | l | l | l  }
       \hline\hline
       $k$ & $n_f$  & $\eta_{G_a}$ & $\eta_{G_b}$ & $\delta_d$ & $\sigma$ & $\Gamma$ \\ \hline
       0.5 & 0.2 & 0.2 & 0.2 & 0.1 & 7.5 & 0.1$I$  \\ \hline \hline
       \end{tabular}
     \end{center}
   \end{table}
   Figure \ref{Outc1} shows the output performance of the proposed approach versus \Lone adaptive controller, the control signal of these two controllers are presented in figure \ref{Contc1}. Figure \ref{Transc1} and \ref{Trans1c1} reveal the idea of prescribed performance and demonstrates the error of these three controllers with respect to pre-assigned prescribed values with high nonlinear uncertainties and nonlinearities as mentioned in case 1.\\      
   The following figures of case 2 overlay the simulation results of the proposed approach as well as two controllers from the literature. we can see in figure \ref{Fig47} the output performance of three controllers, their control signal is presented in \ref{Fig46}, error and transformed error are presented in figure \ref{Fig48}, and finally figure \ref{Fig49} highlights the advantage of the proposed controller. In all, the performance of the proposed approach, its efficiency, and robustness compete with \Lone and $Neuro-Adaptive$.

\begin{figure}[htbp]
           \centering
           \includegraphics[width=13cm, height=7cm]{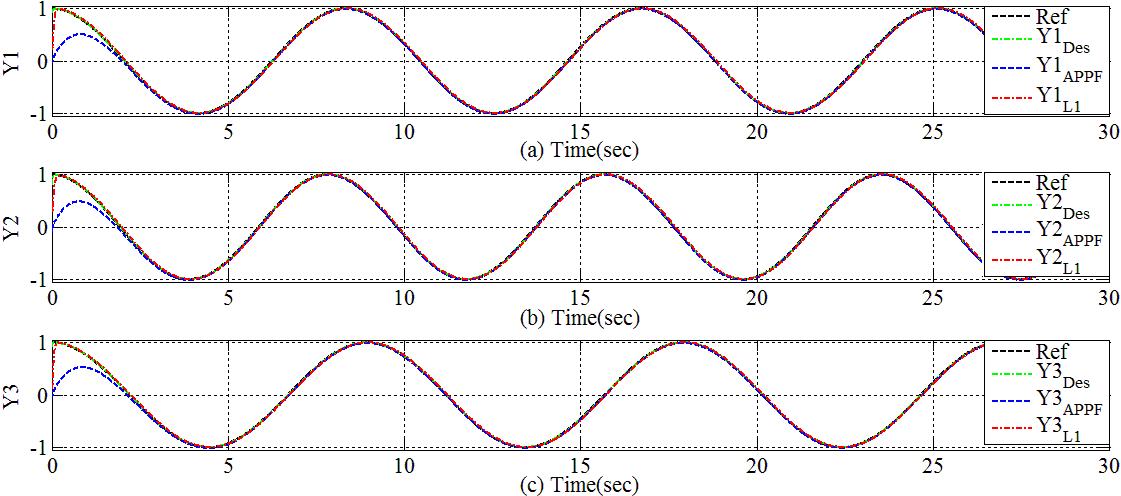}
           \caption{Output Performance of robust MRAC-PPF and \Lone adaptive controller for case 1.}\label{Outc1}
        \end{figure}        
        \begin{figure}[htbp]
           \centering
           \includegraphics[width=13cm, height=7cm]{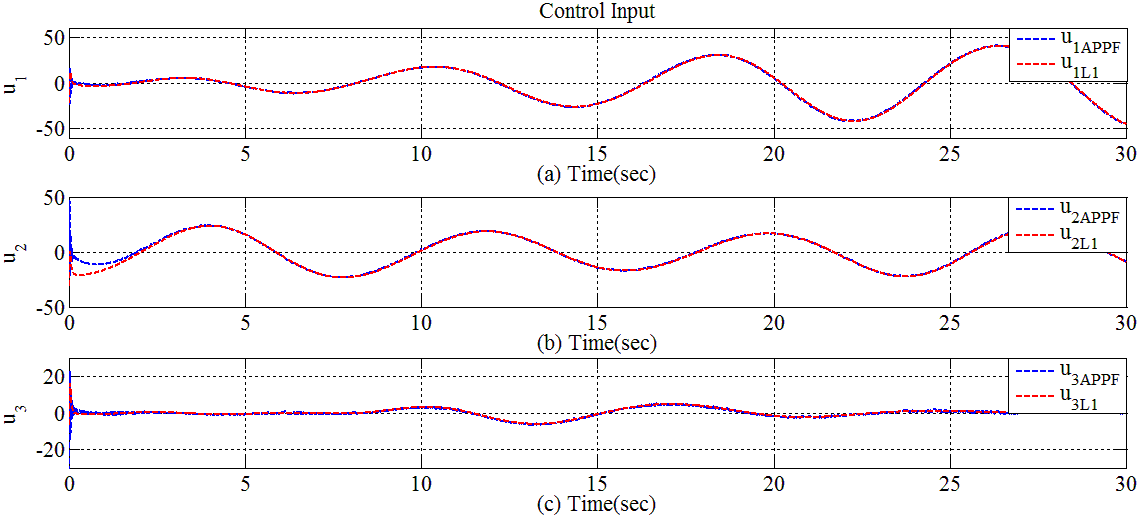}
           \caption{Control Signal of robust MRAC-PPF and \Lone adaptive controller for case 1.}\label{Contc1}
        \end{figure}        
        \begin{figure}[htbp]
           \centering
           \includegraphics[width=13cm, height=7cm]{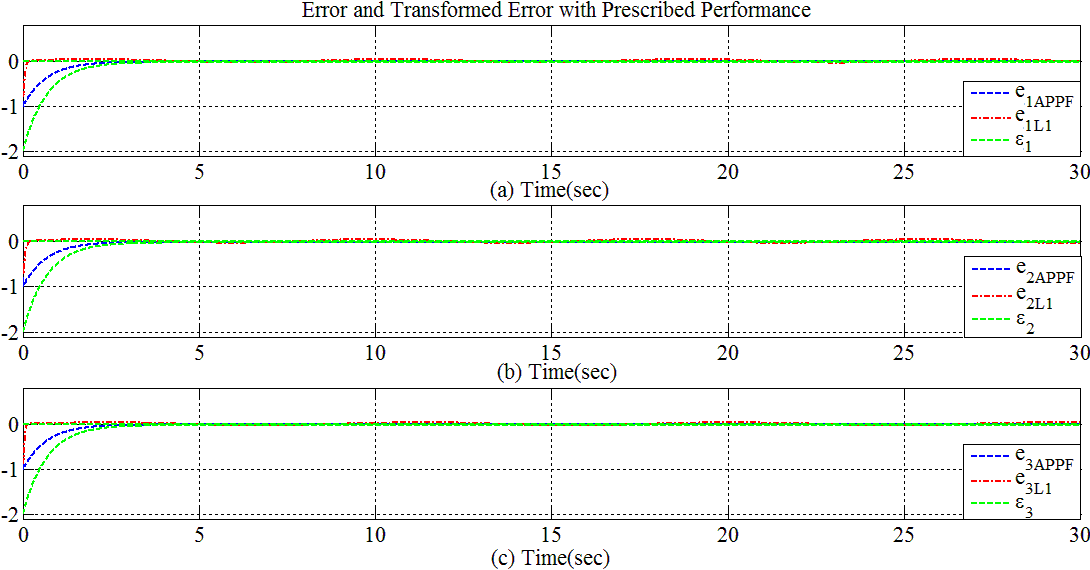}
           \caption{$e_2$ and $\epsilon_2$ of robust MRAC-PPF and \Lone adaptive controller for case 1.}\label{Transc1}
        \end{figure}        
        \begin{figure}[htbp]
           \centering
           \includegraphics[width=13cm, height=7cm]{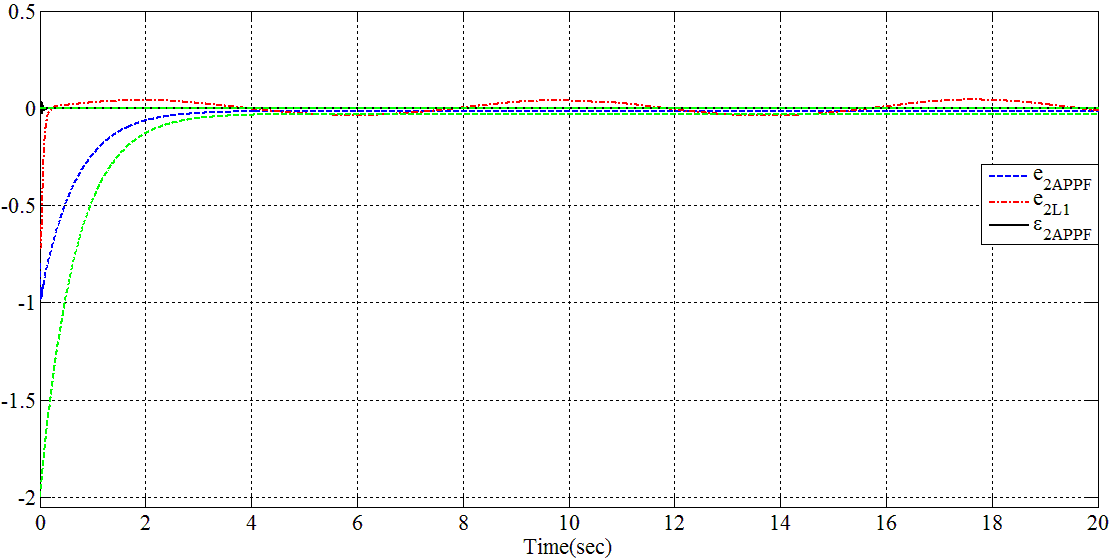}
           \caption{$e_2$ and $\epsilon_2$ of robust MRAC-PPF and \Lone adaptive controller for case 1.}\label{Trans1c1}
        \end{figure}
         \begin{figure}[htbp]
                   \centering
                   \includegraphics[width=13cm, height=7cm]{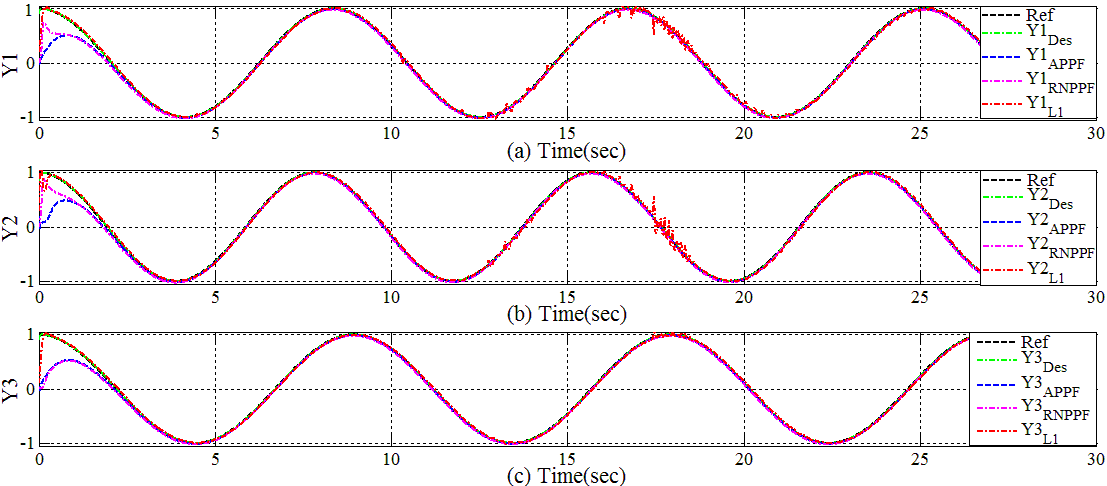}
                   \caption{Output Performance of robust MRAC-PPF, \Lone adaptive controller and $Neuro-Adaptive$ with PPF for case 2.}\label{Fig47}
                \end{figure}        
                \begin{figure}[htbp]
                   \centering
                   \includegraphics[width=13cm, height=7cm]{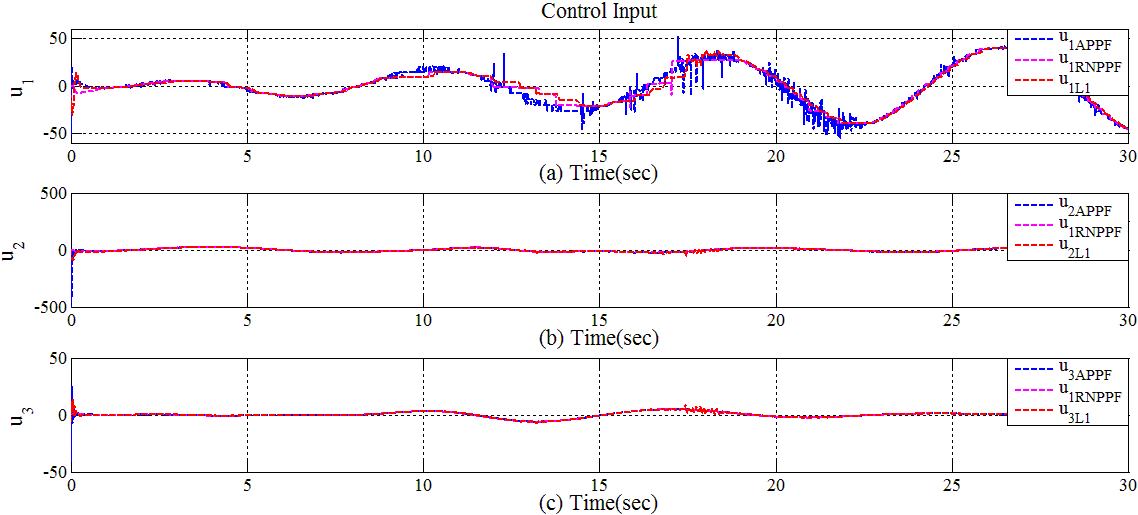}
                   \caption{Control Signal of robust MRAC-PPF, \Lone adaptive controller and $Neuro-Adaptive$ controller with PPF for case 2.}\label{Fig46}
                \end{figure}        
                \begin{figure}[htbp]
                   \centering
                   \includegraphics[width=13cm, height=7cm]{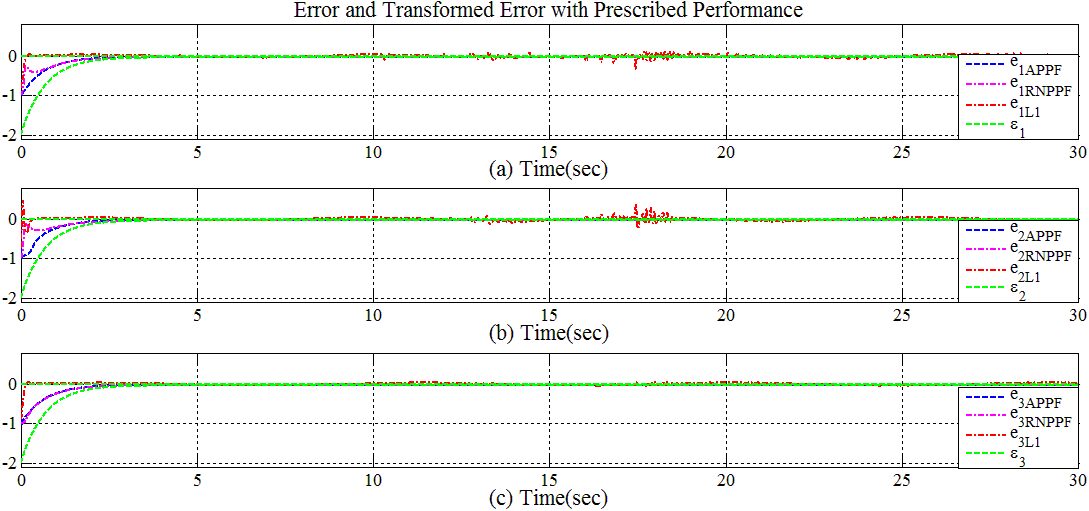}
                   \caption{$e_2$ and $\epsilon_2$ of robust MRAC-PPF, \Lone adaptive controller and $Neuro-Adaptive$ controller with PPF for case 2.}\label{Fig48}
                \end{figure}        
                \begin{figure}[htbp]
                   \centering
                   \includegraphics[width=13cm, height=7cm]{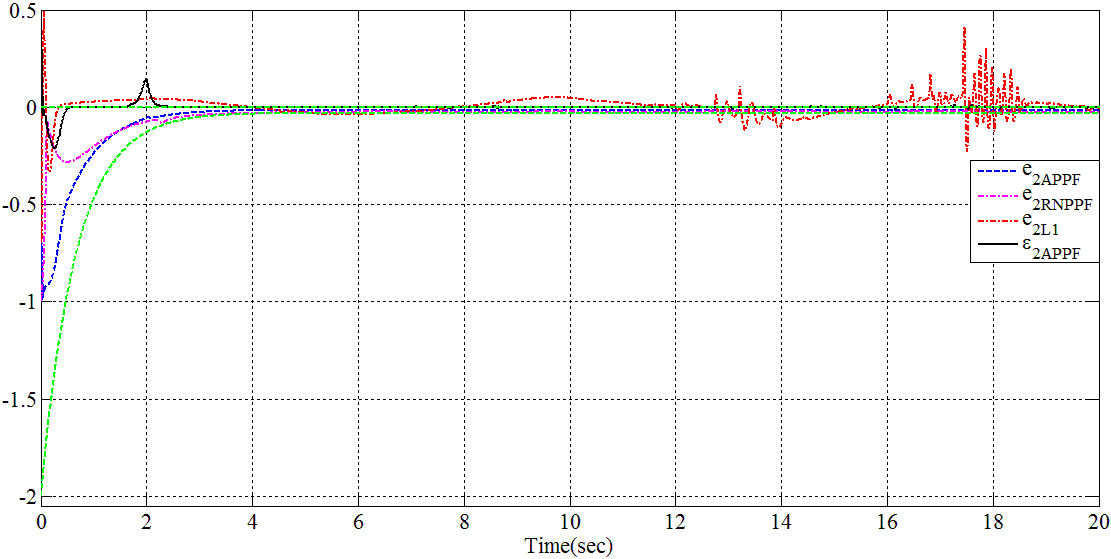}
                   \caption{$e_2$ and $\epsilon_2$ of robust MRAC-PPF, \Lone adaptive controller and $Neuro-Adaptive$ controller with PPF for case 2.}\label{Fig49}
                \end{figure}
                       
\section{Conclusion}
In this chapter, we proposed an adaptive control of multi-input multi-output uncertain high-order nonlinear system capable of guaranteeing a predetermined prescribed performance. The robust stabilization of the transformed error, guaranties the stability and convergence of the constrained tracking error within the set of time varying constraints representing the performance limits. Simulation results demonstrated the efficiency of the proposed approach when compared to \Lone adaptive control and to the neuro-adaptive approach with similar requirement.

\clearpage

\newpage


\chapter{ROBUST ADAPTIVE OBSERVER FOR \Lone ADAPTIVE CONTROLLER}
\section{Introduction}
   Designing a robust adaptive observer for nonlinear systems could be headed in order to estimate inaccessible states from the measured output but can be challenging due to unmodeled dynamics, presence high nonlinearities and time varying uncertainties. In this chapter, robust adaptive observer design for \Lone adaptive controller is mainly adopted from \cite{liu_robust_2009}. The work in \cite{liu_robust_2009} was designed to deal with SISO and MIMO systems with high level of nonlinearities that are assumed to be completely unknown in addition to the presence of structured uncertainties. The chapter is organized as following: section one is an introduction. Problem formulation is presented in section two. The observer design and stability analysis are presented in section three. In section four, discussion of illustrative examples validate the effectiveness of the observer design with \Lone adaptive controller. Finally, the chapter is concluded.

\section{Problem formulationn}
   Consider the following problem:
   \begin{equation}
    \begin{aligned}
     & \dot{x}(t) = Ax(t)+Bf(x,u,t)+g(y,u)\\
     & y = Cx(t)
    \end{aligned}
   \end{equation}
   where $x \in \mathbb{R}^{n}$, $u \in \mathbb{R}^{m}$ and $y \in \mathbb{R}^{p}$  are the system stats (unmeasured),  the control input (unmeasured) and the system output (measured) respectively. $g(y,u)$ is nonlinear function with known parameters and $f(x,u,t)$ is an unknown nonlinear function. Finally, $A$, $B$ and $C$ are constant matrices (known) with appropriate sizes.
   
   The objective of this chapter is to design an adaptive observer for uncertain nonlinear system with unknown dynamics in order to estimate states values for \Lone adaptive controller from the regulated output value. Four basic assumptions will be considered
 \begin{assumption}
    The pair $(A,B)$ is controllable and the pair $(A,C)$ is detectable.
 \end{assumption}
 \begin{assumption}
    Lyapunov function of the system $V(\omega)$ is uniformly bounded and satisfies
   \begin{equation}
    \begin{aligned}
     & \alpha_{1}(||\nu||) \le V_{\nu}(\nu) \le \alpha_{2}(||\nu||)
    \end{aligned}
   \end{equation}
   \begin{equation}
    \begin{aligned}
     \label{eq:ch6non}
     & \frac{\partial \le V_{\nu}(\nu)}{\partial \nu} S(y,\nu) \le -\alpha_{3}(||\nu||)
    \end{aligned}
   \end{equation}
   \begin{equation}
    \begin{aligned}
     & \alpha_{3}(||\nu||) = \tau_{0}V_{\nu}(\nu) - \gamma(||y||) - d_{0}
    \end{aligned}
   \end{equation}
   where $\alpha_{1}$ , $\alpha_{2}$ and $\alpha_{3}$ are positive definite class $K_{\infty}$ functions \cite{khalil_nonlinear_2002}, and $\tau_{0} > 0$ ; $ d_{0} > 0$ are positive constants. $\gamma_{0}$ is a smooth nonnegative function  and has the form of $\gamma(s)=s^{2}\gamma_{0}(s^{2})$ which will be equivalent to $y^{2}\gamma_{0}(y^{2})$ as mentioned in \cite{jiang_design_1998} and $\bar{\epsilon}_0$ is a small positive number.
 \end{assumption}
 \begin{assumption}
    The nonlinear function can be written in the form of
   \begin{equation}
    \begin{aligned}
     & ||f(x,u,t)|| \leq \lambda_{1} + \lambda_{2}||x||\xi(y,u) + \lambda_{3}\zeta(y,u) + \lambda_{4}\alpha(||\nu||)
    \end{aligned}
   \end{equation}
   with $\lambda_{i}\ge 0$, $i=1,2,3,4$ are unknown nonnegative constants, $\alpha(||\cdot||)$ is a class $K_{\infty}$ function and both of $\xi(y,u)$ and $\zeta(y,u)$ are functions assigned arbitrarily nonnegative.
 \end{assumption}
 \begin{assumption}
    $Q$, $P$ are positive definite matrices satisfying
   \begin{equation}
    \begin{aligned}   
     & (A-k_lC)^{T}P + P(A-k_lC) + Q \leq 0\\
     &  PB = C^T
    \end{aligned}  
   \end{equation}
 \end{assumption}
   \section{Robust adaptive observer}
   The observer design is given by
   \begin{equation}
    \begin{aligned}   
     & \dot{\breve{x}} = A\breve{x} - k_l(\breve{y}-y) - \breve{\beta}B(\breve{y}-y)\beta_l + g(y,u)
    \end{aligned}
   \end{equation}
   where $\breve{y} = C\breve{x}$.\\
   Let $e = \breve{x}-x$,$A_l = A-k_lC$ and $\bar{e}_l = \breve{y}-y$
   \begin{equation}
    \begin{aligned}   
     & \dot{e}_l = A_le_l - \breve{\beta}B\bar{e}_l\beta_l - Bf(x,u,t)
    \end{aligned}
   \end{equation}
   The adaptation law of parameter $\breve{\beta}$given by
   \begin{equation}
    \begin{aligned} 
     \label{eq:ch6betadot}  
     & \dot{\breve{\beta}} = \Gamma_l ||\bar{e}_l||\beta_l -  \Gamma_l\sigma_l\breve{\beta}
    \end{aligned}
   \end{equation}
   Where $\Gamma_l$ and $\sigma_l$ are positive constants and $\beta_l$ can be defined by
   \begin{equation}
    \label{eq:ch6betal}
    \begin{aligned}   
     & \beta_l = 1 + \xi^{2}(y,u) + ||\breve{x}||\xi^{2}(y,u) + \eta^{2}(y,u) + \big[\alpha\big(\alpha_1^{-1}(2\delta)\big)\big]^2
    \end{aligned}
   \end{equation}
   $\delta$ is a dynamic signal include unmodeled dynamics and it has the following form
   \begin{equation}
    \begin{aligned}   
     & \dot{\delta} = -\lambda_{0}\delta + \delta_{l},\delta(0)>0
    \end{aligned}
   \end{equation}
   where $\tau_0 \in [0,\tau_0]$ and $\delta_l(y)$ is a smooth nonnegative function $\delta_l(y)=||y||^{2}\gamma_{0}(||y||^{2}) + d_0$ .  As mentioned in \cite{jiang_design_1998,liu_robust_2009} that the relation between dynamic signal and Lyapunov function is
   \begin{equation}
    \begin{aligned}
     \label{eq:ch6delta} 
     & V_{\nu}(\nu) \le \delta + D \\
     & D = max{0,e^{-\tau_0t}V_{\nu}(\nu)-e^{-\lambda_0\delta_0t}}
    \end{aligned}
   \end{equation}
   where $\tau_0>\lambda_0>0$.\\
   The full illustration of L1 adaptive controller with robust adaptive observer is depicted in figure \ref{fig:Chap6_Observer_1}.
   \begin{figure}[h!]
      \centering
      \includegraphics[height=7cm, width=13.5cm]{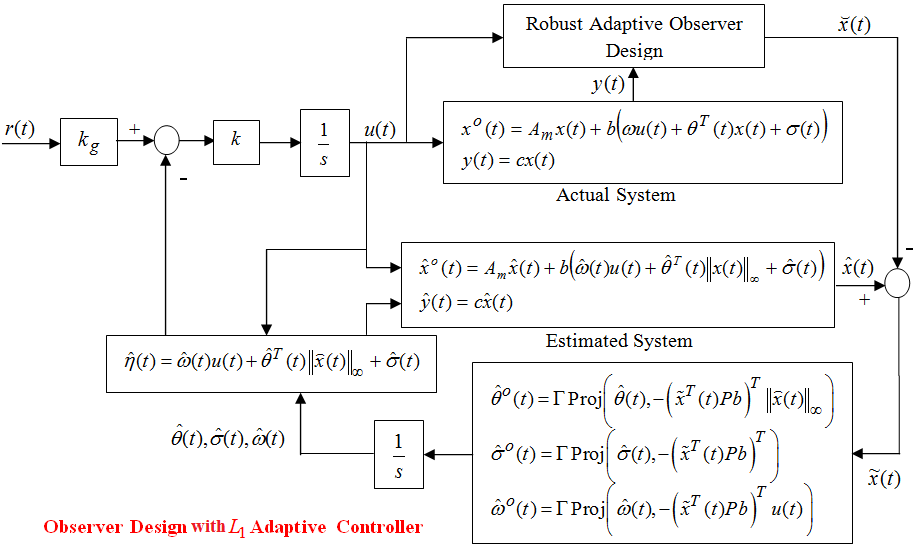}
      \caption{Robust adaptive observer design with \Lone adaptive controller.}
      \label{fig:Chap6_Observer_1}
   \end{figure}
   \subsection{Lyapunov function}
   consider the following Lyapunov candidate
   \begin{equation}
    \begin{aligned}
    \label{eq:ch6Lyap1}
     & V_{\nu} = \frac{1}{2}\big[e_l^TPe_l + \Gamma_l\tilde{\beta}^2]
    \end{aligned}
   \end{equation}
   where $\tilde{\beta} = \beta - \beta^{\ast}$ and $\beta^{\ast}>0$ is a constant representing the desired value of $\beta$. The derivative of \eqref{eq:ch6Lyap1} in addition to the use of \eqref{eq:ch6betadot}, \eqref{eq:ch6betal} and assumption 7.3.
   \begin{equation*}
    \begin{aligned}
     & \dot{V}_{\nu} = \frac{1}{2}e_l(A_l^TP + PA_l)e_l - e_l^TP\hat{\beta}B\bar{e}_l\beta_l - e_l^TPBf(x,t) + \Gamma_l^{-1} \tilde{\beta}\dot{\hat{\beta}}
    \end{aligned}
   \end{equation*}
   \begin{equation*}
    \begin{aligned}
     & \dot{V}_{\nu} = \frac{1}{2}e_lQe_l - e_l^TP\hat{\beta}B\bar{e}_l\beta_l - e_l^TPBf(x,t) + \Gamma_l^{-1} \tilde{\beta}\dot{\hat{\beta}}
    \end{aligned}
   \end{equation*}
   \begin{equation}
    \begin{split}
      \dot{V}_{\nu} = & \frac{1}{2}e_lQe_l - \hat{\beta}||\bar{e}_l||^2\big[1 + \xi^{2}(y,u) + ||\breve{x}||\xi^{2}(y,u) + \eta^{2}(y,u) + \big[\alpha\big(\alpha_1^{-1}(2\delta)\big)\big]^2\big]\\
      & - ||\bar{e}_l||\big[\lambda_{1} + \lambda_{2}||x||\xi(y,u) + \lambda_{3}\zeta(y,u) + \lambda_{4}\alpha(||\nu||)\big] + \Gamma_l^{-1} \tilde{\beta}\dot{\hat{\beta}}
    \end{split}
   \end{equation}
   \begin{equation}
    \begin{split}
      \dot{V}_{\nu} = & \frac{1}{2}e_lQe_l - ||\bar{e}_l||\big[\lambda_{1} + \lambda_{2}||x||\xi(y,u) + \lambda_{3}\zeta(y,u) + \lambda_{4}\alpha(||\nu||)\big] + \sigma_l \tilde{\beta}\beta\\
      & - \beta^{*}||\bar{e}_l||^2\big[1 + \xi^{2}(y,u) + ||\breve{x}||\xi^{2}(y,u) + \eta^{2}(y,u) + \big[\alpha\big(\alpha_1^{-1}(2\delta)\big)\big]^2\big] 
    \end{split}
   \end{equation}
   From \eqref{eq:ch6delta}
   \begin{equation}
    \begin{split}
      \alpha(||\nu||) \le \alpha\big(\alpha_1^{-1}(2\delta)) + \alpha\big(\alpha_1^{-1}(2D))
    \end{split}
   \end{equation}
   \begin{equation}
    \begin{split}
      \dot{V}_{\nu} = & \frac{1}{2}e_lQe_l - ||\bar{e}_l||\lambda_{1} + ||\bar{e}_l||\lambda_{2}||x||\xi(y,u) + ||\bar{e}_l||\lambda_{3}\zeta(y,u)\\ 
      & + ||\bar{e}_l||\lambda_{4}\alpha\big(\alpha_1^{-1}(2\delta)) + ||\bar{e}_l||\lambda_{4}\alpha\big(\alpha_1^{-1}(2D)) - \sigma_l \tilde{\beta}\beta\\
      & - \beta^{*}||\bar{e}_l||^2\big[1 + \xi^{2}(y,u) + ||\breve{x}||\xi^{2}(y,u) + \eta^{2}(y,u) + \big[\alpha\big(\alpha_1^{-1}(2\delta)\big)\big]^2\big] 
    \end{split}
   \end{equation}
   Choosing $\bar{\lambda}_{1} = \lambda_{1} + \lambda_{4}\alpha\big(\alpha_1^{-1}(2D))$ and  $||x||\le ||e_l||+||\breve{x}||$
   \begin{equation}
    \begin{split}
      \dot{V}_{\nu} = & \frac{1}{2}e_lQe_l - ||\bar{e}_l||\bar{\lambda}_{1}+ ||\bar{e}_l||\lambda_{2}||e_l||\xi(y,u)+ ||\bar{e}_l||\lambda_{2}||\breve{x}||\xi(y,u)\\ 
      &  + ||\bar{e}_l||\lambda_{3}\zeta(y,u) + ||\bar{e}_l||\lambda_{4}\alpha\big(\alpha_1^{-1}(2\delta)) - \sigma_l^{-1} \tilde{\beta}\beta\\
      & - \beta^{*}||\bar{e}_l||^2\big[1 + \xi^{2}(y,u) + ||\breve{x}||\xi^{2}(y,u) + \eta^{2}(y,u) + \big[\alpha\big(\alpha_1^{-1}(2\delta)\big)\big]^2\big] 
    \end{split}
   \end{equation}
   \begin{equation}
    \begin{split}
      \dot{V}_{\nu} = & \frac{1}{2}e_lQe_l -  \sigma_l^{-1} \tilde{\beta}\beta - \beta^{*}||\bar{e}_l||^2\beta_l + M  
    \end{split}
   \end{equation}
   Where $M$ includes the rest terms which is equivalent to equation \eqref{eq:ch6non}.
   \section{Results and Discussions}
   Two cases will validate the robustness of robust adaptive observer design with \Lone adaptive controller. The first case represent the observer with high nonlinear SISO system and in the second case and the observer is designed for high nonlinear MIMO system. The nonlinearity, states  and control input are assumed to be completely unknown for previous two cases.\\
   {\bf Example 7.5.1} Consider the following nonlinear SISO system
   \begin{equation*}
    \begin{split}
      & \dot{x} = Ax+B(\omega u + f(x,t))\\
      & y = Cx
    \end{split}
   \end{equation*}
   where $x=[x_1,x_2]^T$ are system states (unmeasured), $u$ is the control input (unmeasured), $y$ is the output (measured). $A$, $B$ and $C$ are known matrices and they indicate that the system is controllable and detectable. The unknown nonlinearity is $f(x,t)$.
   \begin{equation*}
    A = 
    \begin{bmatrix}
     0 & 1\\
     0 & 0 
    \end{bmatrix},
    B =
    \begin{bmatrix}
     0\\
     1 
    \end{bmatrix},
    C =
    \begin{bmatrix}
     1 & 1 
    \end{bmatrix}
   \end{equation*}
   and
   \begin{equation*}
    \begin{split}
      \omega = \frac{75}{s+75}, z(s) = \frac{s-1}{s^2+3s+2}v(s), v(t)= x_1sin(0.2t)+x_2
    \end{split}
   \end{equation*}
   \begin{equation*}
    \begin{split}
      f(x,t) = 2x_1^2 + 2x_2^2 + x_1sin(x_1^2) + x_2cos(x_2^2) + z^2
    \end{split}
   \end{equation*}
   Each of the unmodeled input parameters, uncertainties in the states and disturbances were assigned in compact sets $[\omega_{min},\omega_{max}]\in [0,10]$, $\Delta = 100$ and $\theta_b = 10$. The desired closed loop poles are chosen to be $-1.4\pm j0.743$, the feedback gain = $20$, the adaptation gain$(\Gamma) = 1000000$ and $Q=[\begin{smallmatrix} 1 & 0 \\ 0 & 1\end{smallmatrix}]$. The observer design parameters were selected as $\Gamma_l = 10$, $\sigma_l = 0.0001$, $\lambda_0 = 2.5$, $d_0 = 0.625$ and finally $k_l = [8,64]^T$ . The parameter of the adaptive law $\beta_l$ is defined by $\beta_l = 1 + ||y||^4 + ||\breve{x}||^2||y||^4 + 2\delta$ with $\delta(0) = 1$ and $\breve{\beta}(0) = 1$. The reference input was chosen to $r = cos(0.5t)$ with step change by $+1$ and $-1$ at $14$ and $35$ second respectively in order to validate the robustness of the observer with \Lone adaptive controller.\\
   Figure (\ref{fig:L1_chap6_case1_out}) illustrates the output performance and the control signal of \Lone adaptive controller with the observer design. The actual and estimated states are demonstrated in figure (\ref{fig:L1_chap6_case1_state}). The change in the adaptive estimate  during the control process is revealed in figure (\ref{fig:L1_chap6_case1_beta}).\\
   \begin{figure}[h!]
      \centering
      \includegraphics[height=7cm, width=13.5cm]{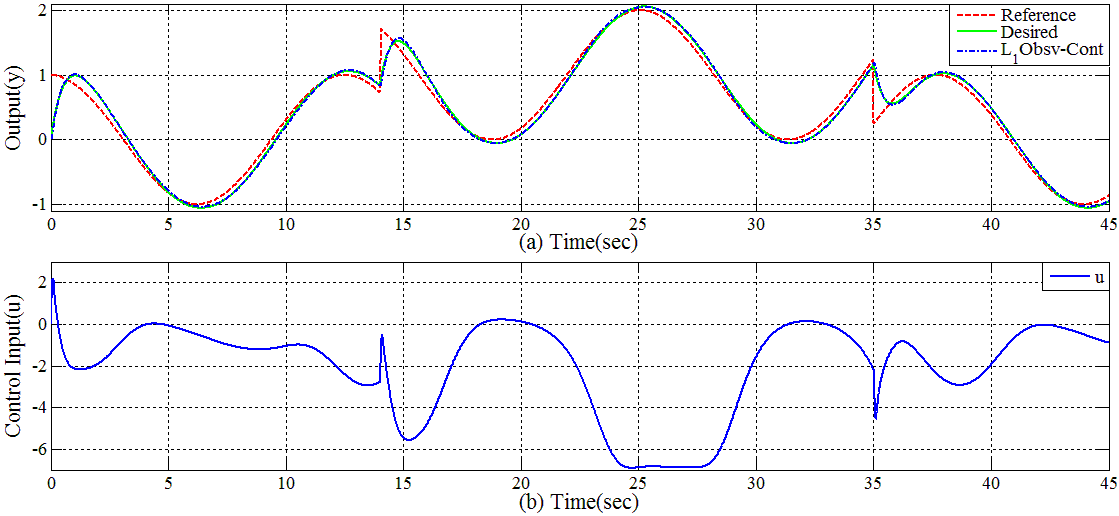}
      \caption{Output performance of \Lone adaptive controller with robust adaptive observer.}
      \label{fig:L1_chap6_case1_out}
   \end{figure}
   \begin{figure}[h!]
      \centering
      \includegraphics[height=6cm, width=13.5cm]{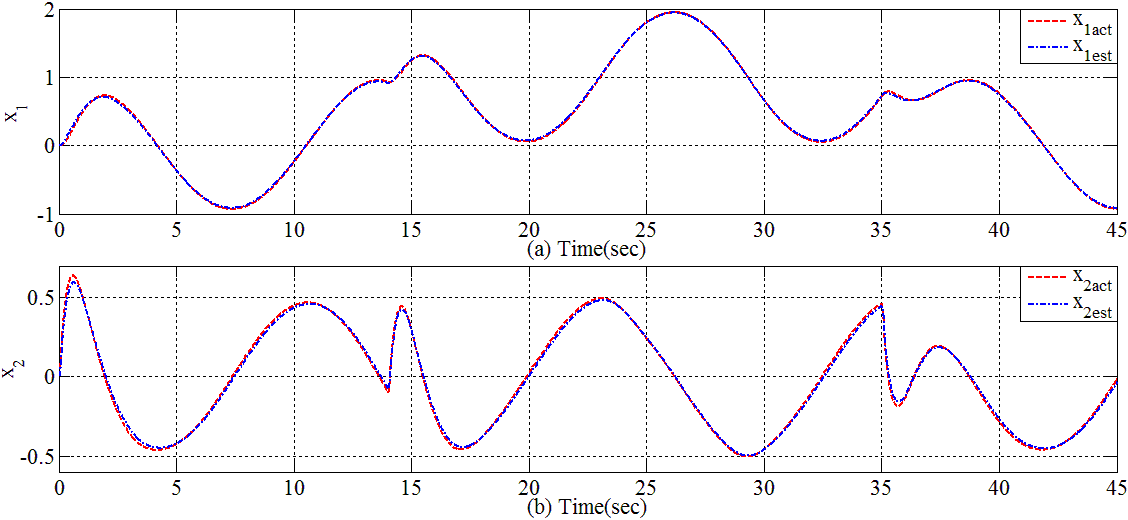}
      \caption{$x$ and $\breve{x}$ of robust observer with \Lone adaptive controller.}
      \label{fig:L1_chap6_case1_state}
   \end{figure}
   \begin{figure}[h!]
      \centering
      \includegraphics[height=6cm, width=13.5cm]{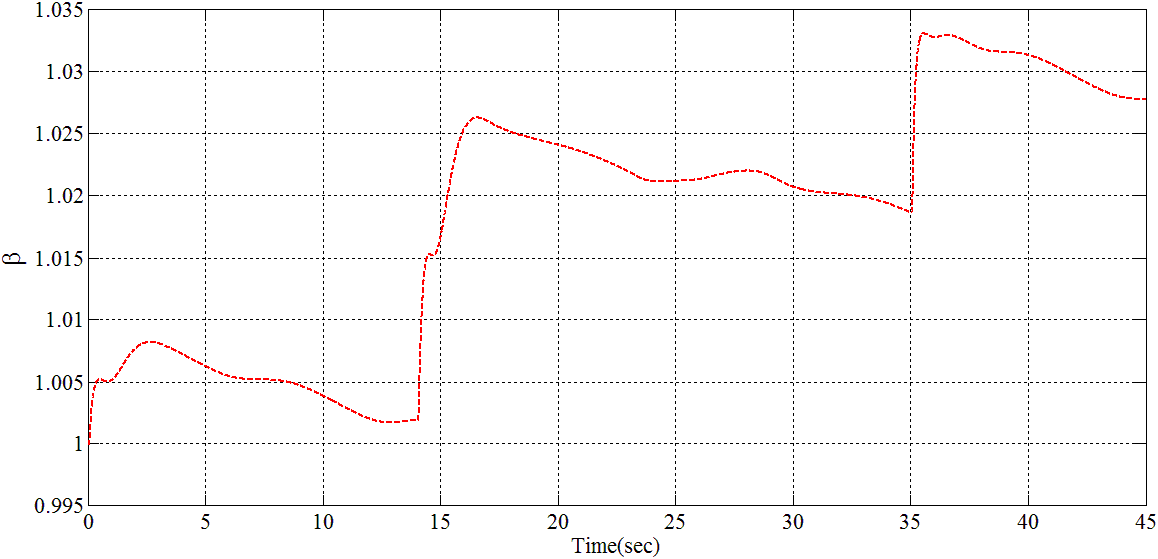}
      \caption{$\breve{\beta}$ of robust observer with \Lone adaptive controller.}
      \label{fig:L1_chap6_case1_beta}
   \end{figure}
   {\bf Example 7.5.2} Consider the following 2-DOF planner robot example 3.3.1 which is similar to our case with some time variant uncertainties in the inertia matrix to be
   \begin{equation*}
     \begin{split}
       M(q) = 
       \begin{bmatrix}
       M_{11}+d_1(t) & M_{12}+d_2(t)\\
       M_{21}+d_2(t) & M_{22}+d_3(t)
       \end{bmatrix}
    \end{split}
   \end{equation*}
   where $d_1(t) = 0.6sin(0.3t)$, $d_2(t) = 0.7sin(0.25t)$ and $d_3(t) = |0.5sin(0.35t)|$ are time varying uncertain parameters included in the model.  Projection operator bounds are  $\hat{\omega} \in [\begin{smallmatrix}[0.3,9.0] & [0.0,0.3] \\ [0.0,0.3] & [0.3,4] \end{smallmatrix}]$, $\Delta = 100$ and $\theta_b = 10$. The desired closed loop poles were chosen to $-10\pm j0.5, -15\pm j0.5$, the feedback gain = $K = [\begin{smallmatrix}20 & 0 \\ 0 & 20 \end{smallmatrix}]$, the adaptation gain$(\Gamma) = 100000$ and $Q=eye(4,4)$. The observer deign parameters were selected as $\Gamma_l = 100$, $\sigma_l = 0.0001$, $\lambda_0 = 2.5$, $d_0 = 0.625$. The parameter of adaptive law $\beta_l$ is defined by $\beta_l = 1 + ||y||^4 + ||\breve{x}||^2||y||^4 + 2\delta$ with $\delta(0) = 1$ and $\breve{\beta}(0) = 1$ and the desired closed loop poles of observers are $-60\pm j0.5, -50\pm j0.5$.\\
   Figure \ref{fig:L1_chap6_case2_out} shows the output performance and the control signal of \Lone adaptive controller with the observer design for joints $q_1$ and $q_2$. In figure (\ref{fig:L1_chap6_case2_state}), actual and observed states are plotted. Finally, figure (\ref{fig:L1_chap6_case2_beta}) illustrates the change in adaptive estimate $\breve{\beta}$ during the control process.
   \begin{figure}[h!]
      \centering
      \includegraphics[height=7cm, width=13.5cm]{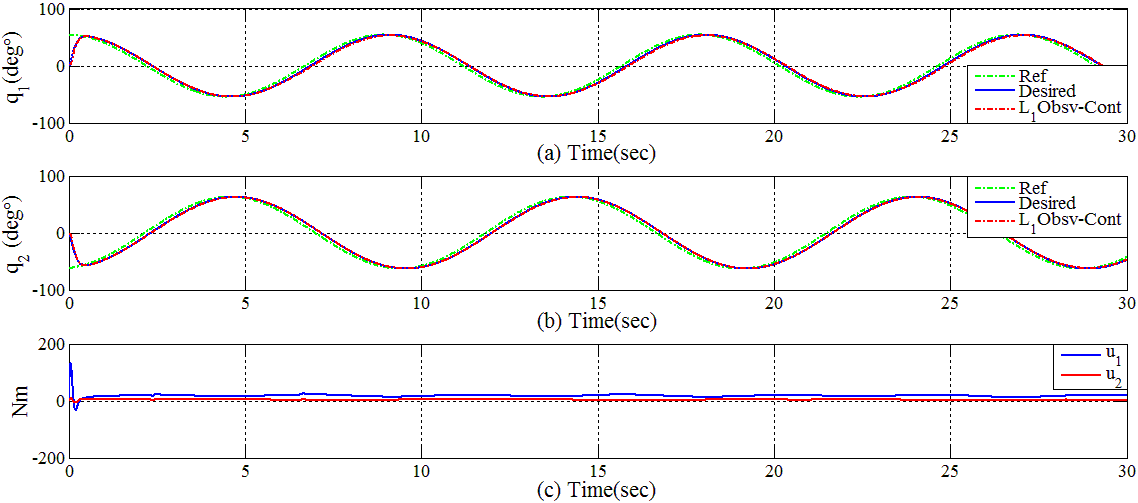}
      \caption{Output performance of \Lone adaptive controller with robust observer for 2 DOF planner robot.}
      \label{fig:L1_chap6_case2_out}
   \end{figure}
   \begin{figure}[h!]
      \centering
      \includegraphics[height=6cm, width=13.5cm]{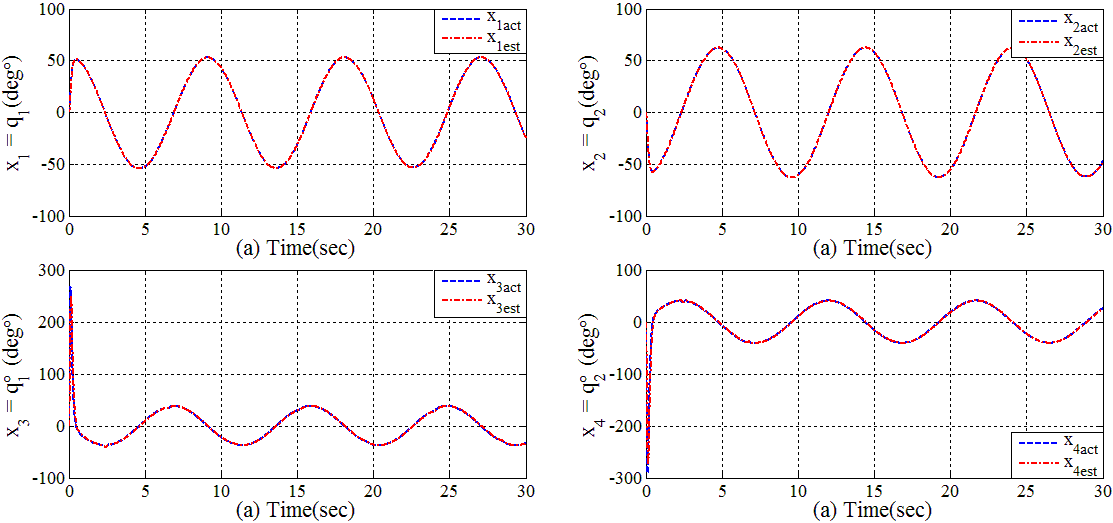}
      \caption{$x$ and $\breve{x}$ of robust observer with \Lone adaptive controller for 2-DOF planer robot.}
      \label{fig:L1_chap6_case2_state}
   \end{figure}
   \begin{figure}[h!]
      \centering
      \includegraphics[height=6cm, width=13.5cm]{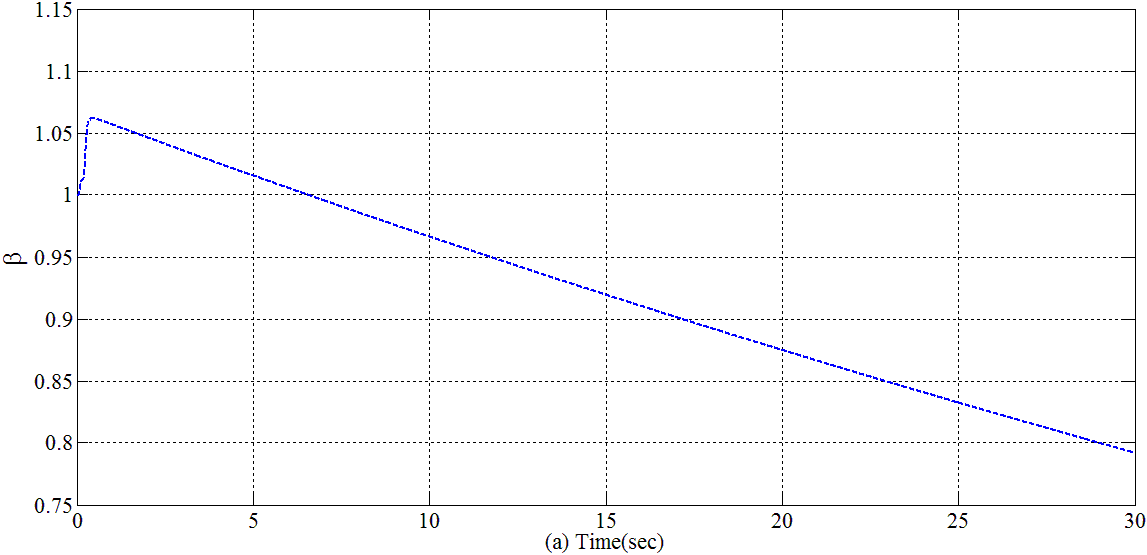}
      \caption{$\breve{\beta}$ in the estimate robust observer with \Lone adaptive controller.}
      \label{fig:L1_chap6_case2_beta}
   \end{figure}
   {\bf Example 7.5.3} Consider simulation problem of quadrotor in example (3.3.2), 
     The observer deign parameters were selected as $\Gamma_l = 100$, $\sigma_l = 0.0001$, $\lambda_0 = 2.5$, $d_0 = 0.625$. The parameter of adaptive law $\beta_l$ is defined by $\beta_l = 1 + ||y||^4 + ||\breve{x}||^2||y||^4 + 2\delta$ with $\delta(0) = 1$ and $\breve{\beta}(0) = 1$ and the desired closed loop poles of observers are $-70\pm j0.5, 75\pm j0.5 and -85\pm j0.5$.\\
     Figure (\ref{fig:L1_chap6_case3_position}) shows the output performance for positions of $x$, $y$ and $z$ of quadrotor. The angles performance and control signal are illustrated in figure (\ref{fig:L1_chap6_case3_Ang}) and (\ref{fig:L1_chap6_case3_Cont}) respectively. Figure (\ref{fig:L1_chap6_case3_3D}) shows the output position in 3D-frame. Finally, figure (\ref{fig:L1_chap6_case3_state}) benchmark the estimated states and actual states. The figure illustrate the robustness of the observer design. 
   \begin{figure}[h!]
      \centering
      \includegraphics[height=7cm, width=13.5cm]{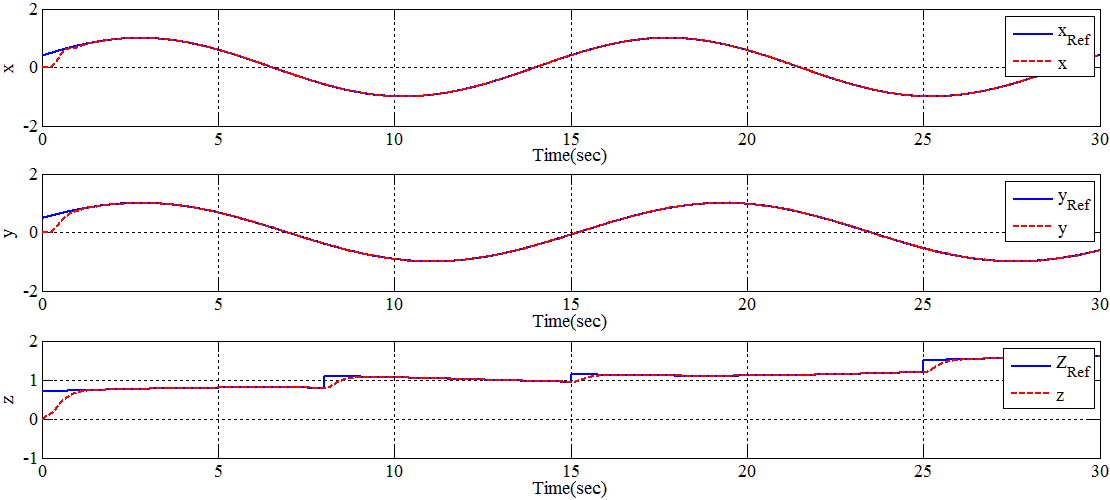}
      \caption{Position performance of \Lone adaptive controller with robust observer for quadrotor.}
      \label{fig:L1_chap6_case3_position}
   \end{figure}
   \begin{figure}[h!]
      \centering
      \includegraphics[height=7cm, width=13.5cm]{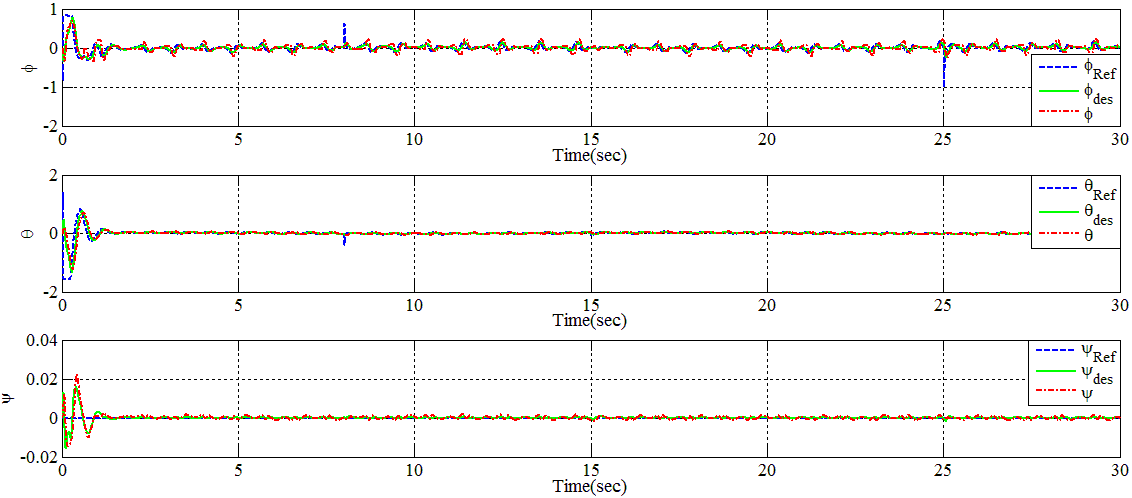}
      \caption{Angles performance of \Lone adaptive controller with robust observer for quadrotor.}
      \label{fig:L1_chap6_case3_Ang}
   \end{figure}
   \begin{figure}[h!]
      \centering
      \includegraphics[height=7cm, width=13.5cm]{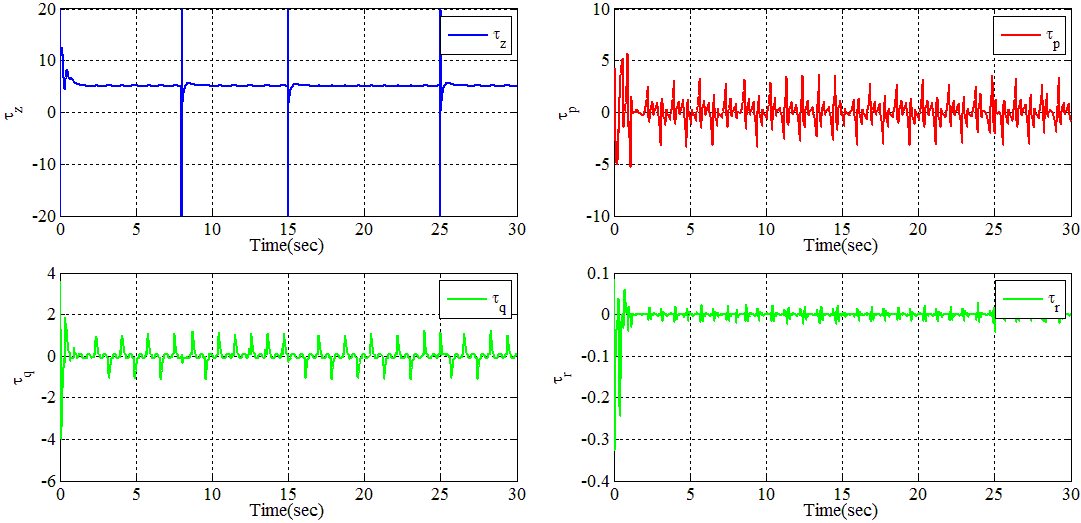}
      \caption{Control signal of \Lone adaptive controller with robust observer for quadrotor.}
      \label{fig:L1_chap6_case3_Cont}
   \end{figure}
   \begin{figure}[h!]
      \centering
      \includegraphics[height=7cm, width=13.5cm]{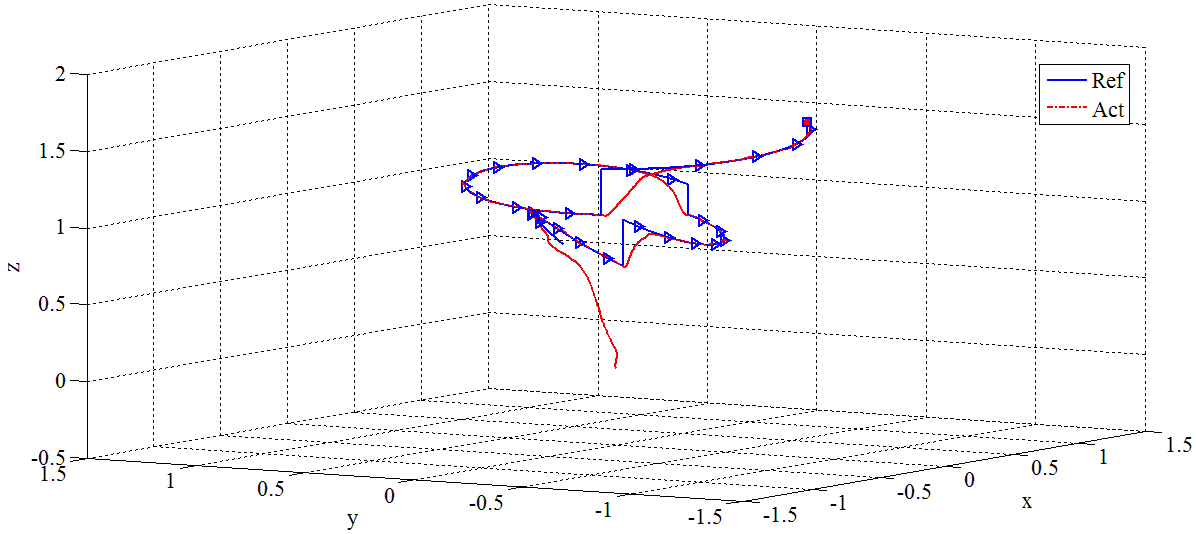}
      \caption{Angles performance of \Lone adaptive controller with robust observer for quadrotor.}
      \label{fig:L1_chap6_case3_3D}
   \end{figure}
   \begin{figure}[h!]
      \centering
      \includegraphics[height=6cm, width=13.5cm]{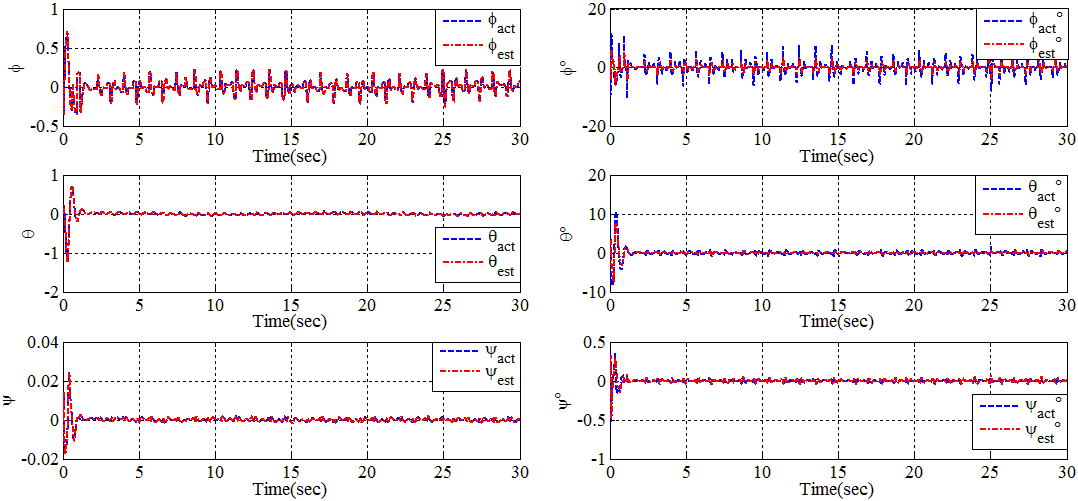}
      \caption{Actual and estimated angles of robust observer with \Lone adaptive controller for quadrotor.}
      \label{fig:L1_chap6_case3_state}
   \end{figure}
   
   \section{Conclusion}
   In this work, robust adaptive observer has been examined with \Lone adaptive controller for nonlinear systems. Nonlinearities are assumed to be completely unknown in addition to unmodeled input parameters and uncertainties. System outputs were available for measurements while states were unmeasurable and control inputs were not used in the observer design. Two illustrative simulations were developed including SISO and MIMO systems to prove the robustness of the observer design with \Lone adaptive controller and to validate the tracking performance. The output performance was impressive and both observed and actual states were very close in their values which validate the efficacy of the observer design with \Lone adaptive controller.

\clearpage

\newpage


\chapter{CONCLUSIONS AND FUTURE WORK}
\section{Summary of Conclusions and Contributions}
   \Lone adaptive controller was applied on different structures of nonlinear systems. In addition, the proposed controllers fuzzy-\Lone adaptive controller and robust MRAC with PPF have been implemented on different nonlinear systems.
   In this thesis, the following problems and results have been presented\\
   {\bf Chapter 3}
   \begin{enumerate}
      \item \Lone adaptive controller has been presented for high nonlinear SISO and MIMO systems with matched and unmatched uncertainties.
      \item High nonlinear systems include UVS such as twin rotor, quadrotor and UAV. Also, two degree of freedom planar robot and other nonlinear systems from recent papers have been simulated.
   \end{enumerate}
   {\bf Chapter 4}
   \begin{enumerate}
      \item Fuzzy filter for \Lone adaptive controller has been proposed for high  nonlinear uncertain systems.
      \item Stability analysis and robustness of the controller has been validated.
      \item The proposed controller showed better results in terms of control signal, robustness margin and tracking capability compared to \Lone adaptive controller.
   \end{enumerate}
   
   {\bf Chapter 5}
   \begin{enumerate}
      \item The work of neuro adaptive control with PPF has been developed successfully.
   \end{enumerate}
    
   {\bf Chapter 6}
   \begin{enumerate}
      \item Robust MRAC with PPF for high nonlinear uncertain systems has been proposed.
      \item Stability analysis and robustness of the controller has been validated.
      \item The proposed controller showed better results from \Lone adaptive controller in case of not-affine systems and it solved the limitations of neuro adaptive control with PPf.
   \end{enumerate}
   
   {\bf Chapter 7}
   \begin{enumerate}
      \item Developed and implemented a robust adaptive observer with \Lone adaptive controller.
      \item The observer showed impressive results with the controller applied to different systems.
   \end{enumerate}
   
\section{Future Work}
      \begin{enumerate}
         \item Optimizing fuzzy membership functions on scale of MIMO systems for fuzzy \Lone-adaptive controller.
         \item Propose MRAC with PPF for high nonlinear systems with unmatched uncertainties.
         \item Propose \Lone adaptive controller with PPF for nonlinear systems.
      \end{enumerate}

\addcontentsline{toc}{chapter}{Bibliography}
\bibliographystyle{IEEEtran}  
\bibliography{bib_ALL} 


\end{document}